\documentclass[11pt]{article}
\usepackage{amsmath, upgreek}
\usepackage{amssymb}
\usepackage{color}
\usepackage{amscd}
\usepackage{xspace}
\usepackage{verbatim}
\usepackage{graphicx}
\newcommand{\mysection}[1]{
\section{#1}\setcounter{equation}{0}}
\title{\bf Radial solutions of scaling invariant
nonlinear elliptic equations with mixed reaction terms}
\author{{\bf Marie-Fran\c{c}oise Bidaut-V\'eron\footnote{\noindent Laboratoire de Math\'{e}matiques et Physique Th\'{e}orique,
Universit\'e de Tours, 37200 Tours, France. E-mail: veronmf@univ-tours.fr},} \\{\bf Marta Garcia-Huidobro \footnote{\noindent
Departamento de Matematicas, Pontifica Universidad Catolica de Chile
Casilla 307, Correo 2, Santiago de Chile. E-mail: mgarcia@mat.puc.cl}}\\
 {\bf Laurent V\'eron \footnote{\noindent
Laboratoire de Math\'{e}matiques et Physique Th\'{e}orique, Universit\'e de Tours, 37200 Tours, France. E-mail: veronl@univ-tours.fr}}\\[2mm]
}

\date{}
\begin{document}
 \maketitle


\newcommand{\txt}[1]{\;\text{ #1 }\;}
\newcommand{\tbf}{\textbf}
\newcommand{\tit}{\textit}
\newcommand{\tsc}{\textsc}
\newcommand{\trm}{\textrm}
\newcommand{\mbf}{\mathbf}
\newcommand{\mrm}{\mathrm}
\newcommand{\bsym}{\boldsymbol}
\newcommand{\scs}{\scriptstyle}
\newcommand{\sss}{\scriptscriptstyle}
\newcommand{\txts}{\textstyle}
\newcommand{\dsps}{\displaystyle}
\newcommand{\fnz}{\footnotesize}
\newcommand{\scz}{\scriptsize}
\newcommand{\be}{\begin{equation}}
\newcommand{\bel}[1]{\begin{equation}\label{#1}}
\newcommand{\ee}{\end{equation}}
\newcommand{\eqnl}[2]{\begin{equation}\label{#1}{#2}\end{equation}}
\newcommand{\barr}{\begin{eqnarray}}
\newcommand{\earr}{\end{eqnarray}}
\newcommand{\bars}{\begin{eqnarray*}}
\newcommand{\ears}{\end{eqnarray*}}
\newcommand{\nnu}{\nonumber \\}
\newtheorem{subn}{\name}
\renewcommand{\thesubn}{}
\newcommand{\bsn}[1]{\def\name{#1}\begin{subn}}
\newcommand{\esn}{\end{subn}}
\newtheorem{sub}{\name}[section]
\newcommand{\dn}[1]{\def\name{#1}}   
\newcommand{\bs}{\begin{sub}}
\newcommand{\es}{\end{sub}}
\newcommand{\bsl}[1]{\begin{sub}\label{#1}}
\newcommand{\bth}[1]{\def\name{Theorem}
\begin{sub}\label{t:#1}}
\newcommand{\blemma}[1]{\def\name{Lemma}
\begin{sub}\label{l:#1}}
\newcommand{\bcor}[1]{\def\name{Corollary}
\begin{sub}\label{c:#1}}
\newcommand{\bdef}[1]{\def\name{Definition}
\begin{sub}\label{d:#1}}
\newcommand{\bprop}[1]{\def\name{Proposition}
\begin{sub}\label{p:#1}}

\newcommand{\R}{\eqref}
\newcommand{\rth}[1]{Theorem~\ref{t:#1}}
\newcommand{\rlemma}[1]{Lemma~\ref{l:#1}}
\newcommand{\rcor}[1]{Corollary~\ref{c:#1}}
\newcommand{\rdef}[1]{Definition~\ref{d:#1}}
\newcommand{\rprop}[1]{Proposition~\ref{p:#1}}
\newcommand{\BA}{\begin{array}}
\newcommand{\EA}{\end{array}}
\newcommand{\BAN}{\renewcommand{\arraystretch}{1.2}
\setlength{\arraycolsep}{2pt}\begin{array}}
\newcommand{\BAV}[2]{\renewcommand{\arraystretch}{#1}
\setlength{\arraycolsep}{#2}\begin{array}}
\newcommand{\BSA}{\begin{subarray}}
\newcommand{\ESA}{\end{subarray}}
\newcommand{\BAL}{\begin{aligned}}
\newcommand{\EAL}{\end{aligned}}
\newcommand{\BALG}{\begin{alignat}}
\newcommand{\EALG}{\end{alignat}}
\newcommand{\BALGN}{\begin{alignat*}}
\newcommand{\EALGN}{\end{alignat*}}
\newcommand{\note}[1]{\textit{#1.}\hspace{2mm}}
\newcommand{\Proof}{\note{Proof}}
\newcommand{\qeda}{\hspace{10mm}\hfill $\square$}
\newcommand{\qed}{\\
${}$ \hfill $\square$}
\newcommand{\Remark}{\note{Remark}}
\newcommand{\modin}{$\,$\\[-4mm] \indent}
\newcommand{\forevery}{\quad \forall}
\newcommand{\set}[1]{\{#1\}}
\newcommand{\setdef}[2]{\{\,#1:\,#2\,\}}
\newcommand{\setm}[2]{\{\,#1\mid #2\,\}}
\newcommand{\mt}{\mapsto}
\newcommand{\lra}{\longrightarrow}
\newcommand{\lla}{\longleftarrow}
\newcommand{\llra}{\longleftrightarrow}
\newcommand{\Lra}{\Longrightarrow}
\newcommand{\Lla}{\Longleftarrow}
\newcommand{\Llra}{\Longleftrightarrow}
\newcommand{\warrow}{\rightharpoonup}
\newcommand{
\paran}[1]{\left (#1 \right )}
\newcommand{\sqbr}[1]{\left [#1 \right ]}
\newcommand{\curlybr}[1]{\left \{#1 \right \}}
\newcommand{\abs}[1]{\left |#1\right |}
\newcommand{\norm}[1]{\left \|#1\right \|}
\newcommand{
\paranb}[1]{\big (#1 \big )}
\newcommand{\lsqbrb}[1]{\big [#1 \big ]}
\newcommand{\lcurlybrb}[1]{\big \{#1 \big \}}
\newcommand{\absb}[1]{\big |#1\big |}
\newcommand{\normb}[1]{\big \|#1\big \|}
\newcommand{
\paranB}[1]{\Big (#1 \Big )}
\newcommand{\absB}[1]{\Big |#1\Big |}
\newcommand{\normB}[1]{\Big \|#1\Big \|}
\newcommand{\produal}[1]{\langle #1 \rangle}

\newcommand{\thkl}{\rule[-.5mm]{.3mm}{3mm}}
\newcommand{\thknorm}[1]{\thkl #1 \thkl\,}
\newcommand{\trinorm}[1]{|\!|\!| #1 |\!|\!|\,}
\newcommand{\bang}[1]{\langle #1 \rangle}
\def\angb<#1>{\langle #1 \rangle}
\newcommand{\vstrut}[1]{\rule{0mm}{#1}}
\newcommand{\rec}[1]{\frac{1}{#1}}
\newcommand{\opname}[1]{\mbox{\rm #1}\,}
\newcommand{\supp}{\opname{supp}}
\newcommand{\dist}{\opname{dist}}
\newcommand{\myfrac}[2]{{\displaystyle \frac{#1}{#2} }}
\newcommand{\myint}[2]{{\displaystyle \int_{#1}^{#2}}}
\newcommand{\mysum}[2]{{\displaystyle \sum_{#1}^{#2}}}
\newcommand {\dint}{{\displaystyle \myint\!\!\myint}}
\newcommand{\q}{\quad}
\newcommand{\qq}{\qquad}
\newcommand{\hsp}[1]{\hspace{#1mm}}
\newcommand{\vsp}[1]{\vspace{#1mm}}
\newcommand{\ity}{\infty}
\newcommand{\prt}{\partial}
\newcommand{\sms}{\setminus}
\newcommand{\ems}{\emptyset}
\newcommand{\ti}{\times}
\newcommand{\pr}{^\prime}
\newcommand{\ppr}{^{\prime\prime}}
\newcommand{\tl}{\tilde}
\newcommand{\sbs}{\subset}
\newcommand{\sbeq}{\subseteq}
\newcommand{\nind}{\noindent}
\newcommand{\ind}{\indent}
\newcommand{\ovl}{\overline}
\newcommand{\unl}{\underline}
\newcommand{\nin}{\not\in}
\newcommand{\pfrac}[2]{\genfrac{(}{)}{}{}{#1}{#2}}

\def\ga{\alpha}     \def\gb{\beta}       \def\gg{\gamma}
\def\gc{\chi}       \def\gd{\delta}      \def\ge{\epsilon}
\def\gth{\theta}                         \def\vge{\varepsilon}
\def\gf{\phi}       \def\vgf{\varphi}    \def\gh{\eta}
\def\gi{\iota}      \def\gk{\kappa}      \def\gl{\lambda}
\def\gm{\mu}        \def\gn{\nu}         \def\gp{\pi}
\def\vgp{\varpi}    \def\gr{\rho}        \def\vgr{\varrho}
\def\gs{\sigma}     \def\vgs{\varsigma}  \def\gt{\tau}
\def\gu{\upsilon}   \def\gv{\vartheta}   \def\gw{\omega}
\def\gx{\xi}        \def\gy{\psi}        \def\gz{\zeta}
\def\Gg{\Gamma}     \def\Gd{\Delta}      \def\Gf{\Phi}
\def\Gth{\Theta}
\def\Gl{\Lambda}    \def\Gs{\Sigma}      \def\Gp{\Pi}
\def\Gw{\Omega}     \def\Gx{\Xi}         \def\Gy{\Psi}

\def\CS{{\mathcal S}}   \def\CM{{\mathcal M}}   \def\CN{{\mathcal N}}
\def\CR{{\mathcal R}}   \def\CO{{\mathcal O}}   \def\CP{{\mathcal P}}
\def\CA{{\mathcal A}}   \def\CB{{\mathcal B}}   \def\CC{{\mathcal C}}
\def\CD{{\mathcal D}}   \def\CE{{\mathcal E}}   \def\CF{{\mathcal F}}
\def\CG{{\mathcal G}}   \def\CH{{\mathcal H}}   \def\CI{{\mathcal I}}
\def\CJ{{\mathcal J}}   \def\CK{{\mathcal K}}   \def\CL{{\mathcal L}}
\def\CT{{\mathcal T}}   \def\CU{{\mathcal U}}   \def\CV{{\mathcal V}}
\def\CZ{{\mathcal Z}}   \def\CX{{\mathcal X}}   \def\CY{{\mathcal Y}}
\def\CW{{\mathcal W}} \def\CQ{{\mathcal Q}}
\def\BBA {\mathbb A}   \def\BBb {\mathbb B}    \def\BBC {\mathbb C}
\def\BBD {\mathbb D}   \def\BBE {\mathbb E}    \def\BBF {\mathbb F}
\def\BBG {\mathbb G}   \def\BBH {\mathbb H}    \def\BBI {\mathbb I}
\def\BBJ {\mathbb J}   \def\BBK {\mathbb K}    \def\BBL {\mathbb L}
\def\BBM {\mathbb M}   \def\BBN {\mathbb N}    \def\BBO {\mathbb O}
\def\BBP {\mathbb P}   \def\BBR {\mathbb R}    \def\BBS {\mathbb S}
\def\BBT {\mathbb T}   \def\BBU {\mathbb U}    \def\BBV {\mathbb V}
\def\BBW {\mathbb W}   \def\BBX {\mathbb X}    \def\BBY {\mathbb Y}
\def\BBZ {\mathbb Z}

\def\GTA {\mathfrak A}   \def\GTB {\mathfrak B}    \def\GTC {\mathfrak C}
\def\GTD {\mathfrak D}   \def\GTE {\mathfrak E}    \def\GTF {\mathfrak F}
\def\GTG {\mathfrak G}   \def\GTH {\mathfrak H}    \def\GTI {\mathfrak I}
\def\GTJ {\mathfrak J}   \def\GTK {\mathfrak K}    \def\GTL {\mathfrak L}
\def\GTM {\mathfrak M}   \def\GTN {\mathfrak N}    \def\GTO {\mathfrak O}
\def\GTP {\mathfrak P}   \def\GTR {\mathfrak R}    \def\GTS {\mathfrak S}
\def\GTT {\mathfrak T}   \def\GTU {\mathfrak U}    \def\GTV {\mathfrak V}
\def\GTW {\mathfrak W}   \def\GTX {\mathfrak X}    \def\GTY {\mathfrak Y}
\def\GTZ {\mathfrak Z}   \def\GTQ {\mathfrak Q}

\font\Sym= msam10 
\def\SYM#1{\hbox{\Sym #1}}
\newcommand{\bdw}{\prt\Gw\xspace}
\date{}
\maketitle\medskip

\noindent{\small {\bf Abstract} We study global properties of positive radial solutions of $-\Gd u=u^p+M\abs{\nabla u}^{\frac{2p}{p+1}}$ in $\BBR^N$ where $p>1$ and $M$ is a real number. We prove the existence or the non-existence of ground states and of solutions with singularity at $0$ according to the values of $M$ and $p$.
}\smallskip

\tableofcontents
\vspace{1mm}
\hspace{.05in}
\medskip
\mysection{Introduction}
The aim of this article is  to study local and global properties of positive radial solutions of the equation
\bel{I-0}
-\Gd u=|u|^{p-1}u+M\abs{\nabla u}^{\frac{2p}{p+1}},
\ee
in $\BBR^N$ or $\BBR^N\setminus\{0\}$ where  $p>1$ and $M$ is a real parameter. This is a particular case of the following class of equations
\bel{I-0-1}
-\Gd u=|u|^{p-1}u+M\abs{\nabla u}^{q},
\ee
where $q>1$ which has been the subject or many works  in the radial case when $M<0$, where a basic observation is that the two terms $|u|^{p-1}u$ and $M\abs{\nabla u}^{q}$ are in competition. The first  work in that case is due to Chipot and Weissler \cite{ChWe} who, in particular, solved completely the case $N=1$, then Serrin and Zou \cite{SeZo} performed a very detailed analysis. Much less is known in the case $M>0$. Under the scaling transformation $T_k$ defined for $k>0$ by
\bel{I-2}
u_k:=T_k[u](x)=k^{\frac{2}{p-1}}u(kx),
\ee
$(\ref{I-0-1})$ becomes
\bel{I-0-1a}
-\Gd u_k=|u_k|^{p-1}u_k+k^{\frac{2p-q(p+1)}{p-1}}M\abs{\nabla u_k}^{q},
\ee
Therefore, if $q\neq\frac{2p}{p+1}$, $(\ref{I-0-1})$ can be reduced to
\bel{I-0-1b}
-\Gd u=|u|^{p-1}u\pm\abs{\nabla u}^{\frac{2p}{p+1}}.
\ee
Moreover, when $q<\frac{2p}{p+1}$, the limit equation of $(\ref{I-0-1a})$ when $k\to 0$ is the Lane-Emden equation
\bel{I-0-1c}
-\Gd u=|u|^{p-1}u,
\ee
and thus the exponent $p$ is dominant. The other scaling transformation
\bel{I-2a}
v_k:=S_k[u](x)=k^{\frac{2-q}{q-1}}u(kx),
\ee
transforms $(\ref{I-0-1})$ into
\bel{I-0-1d}
-\Gd v_k=k^{\frac{q-p(2-q)}{q-1}}|v_k|^{p-1}v_k+M\abs{\nabla v_k}^{q},
\ee
and if $q>\frac{2p}{p+1}$, the limit equation of $(\ref{I-0-1d})$ when $k\to 0$ is the Riccati equation
\bel{I-0-1e}
-\Gd v=M\abs{\nabla v}^{q},
\ee
therefore the exponent $q$ is dominant. In  \cite{ChWe} and \cite{SeZo} most of the study deals with the case $q\neq\frac{2p}{p+1}$. In the {\it critical case} i.e. when
 \bel{I-x}
q=\frac{2p}{p+1},
\ee
then not only the sign of $M$  but also its value plays a fundamental role, with a delicate interaction with the exponent $p$. Notice that an equivalent form of $(\ref{I-0})$ is
\bel{I-0-1f}
-\Gd v=\gl \abs{v}^{p-1}v\pm \abs{\nabla v}^{q}
\ee
with $\gl>0$.  In the critical case first studies in the case $M<0$ are due to Chipot and Weissler \cite{ChWe} for $N=1$. The case $N\geq 2$ was left open by Serrin and Zou \cite{SeZo} and the first partial results are due to Fila and Quittner \cite{FiQu} and Voirol \cite{Vo, Vo2}. The case $M>0$ was not considered.\smallskip

The equation $(\ref{I-0})$ is the stationary part of the associated parabolic equation
\bel{I-0-2}
\prt_tu-\Gd u-|u|^{p-1}u-M\abs{\nabla u}^{q}=0.
\ee
which is studied in \cite{ChWe} and \cite{Fi}, where one of the aims was to find conditions for the blow-up of positive solutions. A general survey with several open problems can be found in \cite{So}.\smallskip

In the non radial case an important contribution dealing with a priori estimates of local positive solutions of $(\ref{I-0-1})$ and existence or non-existence of entire positive solution in $\BBR^N$ is due to the authors \cite{BiGaVe2}. In this paper we complete the results of \cite{BiGaVe2} in presenting a quite exhaustive study of the radial solutions of $(\ref{I-0})$ for any real number $M$.\smallskip

The radial solutions of $(\ref{I-0})$ are functions $r\mapsto u(r)$ defined in $(0,\infty)$ where they satisfy
\bel{I-1}
-u_{rr}-\myfrac{N-1}{r}u_{r}=|u|^{p-1}u+M\abs{u_{r}}^{\frac{2p}{p+1}}.
\ee
Because of the invariance of $(\ref{I-1})$ under the transformation $T_k$ there exists
 an autonomous variant of $(\ref{I-0})$ obtained by setting
\bel{I-1'} u(r)=r^{-\frac{2}{p-1}}x(t)\quad\text{ with } \;t=\ln r.  \ee
Then $(\ref{I-1})$ becomes
\bel{I-3}
x_{tt}+Lx_{t}-\myfrac{2K}{p-1}x+|x|^{p-1}x+M\abs{\myfrac{2K}{p-1}x-x_{t}}^{\frac{2p}{p+1}}=0
\ee
with
\bel{I-3v}\BA {lll}\displaystyle
K=\myfrac{(N-2)p-N}{p-1}\,\;\text {and }\;L=\myfrac{(N-2)p-(N+2)}{p-1}=K-\myfrac{2}{p-1}.
\EA\ee
Putting $y(t)=-r^{\frac{p+1}{p-1}}u_{r}(r)$, then $(x(t),y(t))$ satisfies the system
\bel{I-3x}\BA {lll}\displaystyle
x_{t}=\myfrac{2}{p-1}x-y\\[2mm]
y_{t}=-Ky+|x|^{p-1}x+M\abs{y}^{\frac{2p}{p+1}}.
\EA\ee
We are mainly interested in the trajectories of the system which remain in the first quarter {\bf Q}=$\{(x,y)\in\BBR^2:x>0,y>0\}$. Indeed, among these trajectories, we find the ones corresponding to {\it ground states}, i.e. positive $C^2$ solutions $u$ of $(\ref{I-1})$  which are defined on $[0,\infty)$. They verify $u_r(0)=0$ and actually they are $C^\infty$ on $ (0,\infty)$. Using the invariance of the equation under $T_k$ all the ground states can be derived by scaling from a unique one which satisfies $u(0)=1$. Since it is easy to prove that such a solution $u$ is decreasing, in the variables $(x,y)$, a ground state is a trajectory of $(\ref{I-3x})$ in {\bf Q}, defined on $\BBR$ and satisfying $\displaystyle\lim_{t\to-\infty}\frac{y(t)}{x(t)}=0$.  The corresponding trajectory is denoted by ${\bf T}_{reg}$.\smallskip

Contrarily to the case of the Lane-Emden equation $(\ref{I-0-1c})$, there exists no natural Lyapunov function when $M\neq 0$. This makes the study much more delicate and it is based upon a phase plane analysis. The solutions of $(\ref{I-1})$  invariant under $T_k$ for any $k>0$ correspond to constant solutions of $(\ref{I-3})$ and have the form
\bel{I-4}
U(r)=Xr^{-\frac{2}{p-1}}\quad\text{for all } r>0,
\ee
where $X$ is a positive root of
\bel{I-5}
X^{p-1}+M\left(\frac{2}{p-1}\right)^{\frac{2p}{p+1}}X^{\frac{p-1}{p+1}}-\frac{2K}{p-1}=0.
\ee
This equation plays a fundamental role in the description of the set of solutions of $(\ref{I-1})$.  The following constant, defined for $N=1,2$ and $p>1$ or $N\geq 3$ and $1<p\leq  \frac{N}{N-2}$ has an important role in the description of the set roots of $(\ref{I-5})$,
\bel{I-6}
\gm^*(N)=(p+1)\left(\frac{N-(N-2)p}{2p}\right)^{\frac{p}{p+1}}.
\ee
When the is no ambiguity we write $\gm^*:=\gm^*(N)$. This set is described in the following proposition. \medskip

\nind{\bf Proposition 1} {\it 1- If $M\geq 0$, equation $(\ref{I-5})$ admits a positive root, necessarily unique, if and only if $N\geq 3$ and
$p>\frac{N}{N-2}$.\smallskip

\nind 2- If $M<0$ and $p\geq \frac{N}{N-2}$, equation $(\ref{I-5})$ admits a unique positive root $X_M$.\smallskip

\nind 3- If $M<0$ and either $N=1,2$ and $p>1$ or  $N\geq 3$ and
$1<p\leq \frac{N}{N-2}$, there exists no positive root of $(\ref{I-5})$ if $-\gm^*<M<0$, a unique positive root if $M=-\gm^*<0$ and two positive roots $X_{1,M}<X_{2,M}$ if $M<-\gm^*$.
}
\medskip

We also set $Y_M=\frac{2}{p-1}X_M$ and $P_M=(X_M,Y_M)$ (resp. $Y_{j,M}=\frac{2}{p-1}X_{j,M}$ and $P_{j,M}=(X_{j,M},Y_{j,M})$, for j=1,2) and define the corresponding singular solutions $U_M(r)=X_Mr^{-\frac{2}{p-1}}$ (resp. $U_{j,M}(r)=X_{j,M}r^{-\frac{2}{p-1}}$).\smallskip

Recall briefly the description of the positive solutions of the Lane-Emden equation $(\ref{I-0-1c})$, i.e. $M=0$: there exists radial ground states if and only if $N\geq 3$ and $p\geq\frac{N+2}{N-2}$. If $p=\frac{N+2}{N-2}$ these ground states are explicit and they satisfy $\lim_{r\to\infty}r^{N-2}u(r)=c>0$. There exist infinitely many singular solutions $u$ ondulating around $U_M$. Note that a ground state corresponds to a homoclinic orbit at $0$ for system $(\ref{I-3x})$ and these singular solutions are cycles surrounding $P_M$. We recall that an orbit of $(\ref{I-3x})$  which connects two different equilibria (resp. the same equilibrium) when $t\in\BBR$ is called {\it heteroclinic} (resp. {\it  homoclinic}).\smallskip


Without the radiality assumptions, and using a delicate combination of refined Bernstein techniques and Keller-Osserman estimate we have obtained in \cite[Theorems C, D]{BiGaVe2} a series of general a priori estimates for any positive solution of $(\ref{I-0})$, in an arbitrary domain of $\BBR^N$ in the case $N\geq 1$, $p>1$  and $q=\frac{2p}{p+1}$ and $M>0$. In particular we proved there that if $p>1$ and
$M>M_\dag:=\left(\frac{p-1}{p+1}\right)^{\frac{p-1}{p+1}}\left(\frac{N(p+1)^2}{4p}\right)^{\frac{p}{p+1}}$,
or if $N\geq 2$, $p<\frac{N+3}{N-1}$ and $M>0$
equation $(\ref{I-0})$ admits no ground state. \smallskip

In the sequel we describe the ground states and the singular global solutions of $(\ref{I-1})$ in $\BBR^N\setminus\{0\}$. Concerning the ground states, we discuss according to the sign of $M$ and the value of $p$. The next value of
$M$ appears when we linearize the system $(\ref{I-3x})$ at the equilibrium $P_M$,
\bel{I-8}\overline M=\overline M(N,p)=\myfrac{(p+1)\left((N-2)p-N-2\right)}{(4p)^{\frac{p}{p+1}}\left((N-2)(p-1)^2+4\right)^\frac{1}{p+1}}.
\ee
Then $\overline M$ is positive (resp. negative) if $p>\frac{N+2}{N-2}$ (resp. $p<\frac{N+2}{N-2}$  and we set $\overline\gm=\abs{\overline M}$). It is easy to see that if $M=\overline M$ then the characteristic values of the linearized operator at $P_M$ are purely imaginary.  Notice that $\overline M$ is positive (resp. negative) according
$p>\frac{N+2}{N-2}$ (resp. $p<\frac{N+2}{N-2}$).
\medskip


\nind{\bf Theorem A} {\it Let $N\geq 1$, $p>1$ and $M> 0$.\smallskip

\nind 1- For any $1<p\leq\frac{N+2}{N-2}$ if $N\geq 3$, and any $p>1$ if $N=1,2$, then equation $(\ref{I-1})$ admits no ground state.

\smallskip

\nind 2- If $N\geq 3$ and $p>\frac{N+2}{N-2}$, there exist constants $\tilde M_{min},\tilde M_{max}$ verifying
$$0<\overline M<\tilde M_{min}\leq \tilde M_{max},$$
 such that \\
- if $0<M<\tilde M_{min}$ there exist ground states $u$ satisfying $u(r)\sim U_M(r)$ when $r\to\infty$.\\
- if $M=\tilde M_{min}$ or $M=\tilde M_{max}$ there exists a  ground state $u$ minimal at infinity, that is satisfying $\displaystyle\lim_{r\to\infty}r^{N-2}u(r)= c>0$.\\
- for $M>\tilde M_{max}$ there exists no radial ground state.
}\medskip

The values of $\tilde M_{min}$ and $\tilde M_{max}$ appear as transition values for which the ground state still exists but it is
smaller than the others at infinity; it is of order $r^{2-N}$ instead of $r^{-\frac{2}{p-1}}$. They are not explicit but they can be estimated in function of $N$ and $p$. It is a numerical evidence that $\tilde M_{min}=\tilde M_{max}$ in the phase plane analysis of system $(\ref{I-3x})$ and we conjecture that this is true. When $M=\tilde M_{min}$ or $\tilde M_{max}$, the system $(\ref{I-3x})$ admits homoclinic trajectories.  We prove that the system $(\ref{I-3x})$ admits a Hopf bifurcation when $M=\overline M$. When $p>\frac{N+2}{N-2}$ we also prove the existence of different types of positive singular solutions\medskip


\nind{\bf Theorem A'} {\it Let $N\geq 3$. \smallskip

\nind 1- If $\frac{N}{N-2}<p\leq \frac{N+2}{N-2}$ for any $M>0$ there exist a unique (always up to a scaling transformation) positive singular solutions $u$ of $(\ref{I-1})$ satisfying $\displaystyle\lim_{r\to 0}r^{\frac{2}{p-1}}u(r)=X_M$ and  $\displaystyle\lim_{r\to \infty}r^{N-2}u(r)=c>0$
\smallskip

\nind 2- If $p>\frac{N+2}{N-2}$, then\\
\nind (i) If $M>\tilde M_{max}$, there exists a unique singular solution $u$ of $(\ref{I-1})$ with the same behaviour as in 1. \\
\nind (ii)  If $\overline M<M<\tilde M_{min}$ there exist positive singular solutions $u$ ondulating around $U_M$ on $\BBR$.}\medskip

In terms of the system $(\ref{I-3x})$ the 1) and 2-(i) correspond to the existence of a heteroclinic orbit in {\bf Q} connecting $P_M$ to $(0,0)$ and (ii)  to the existence of a cycle in {\bf Q} surrounding $P_M$.\medskip

 When $M$ is negative, the precise description of the trajectories of $(\ref{I-3x})$ depends also on the value of $p$ with respect to $\frac N{N-2}$. It is proved in \cite[Th. B, E]{BiGaVe2} that for $N\geq 3$ and $1<p<\frac{N+2}{N-2}$ there exists $\ge_0>0$ such that if $\abs{M}\leq\ge_0$ equation $(\ref{I-0})$ admits no positive solution in $\BBR^N$. The same conclusion holds if $N\geq 3$,  $1<p\leq\frac N{N-2}$ (or $N=2$ and $p>1$) and $M>-\gm^*$. We first consider the case $p\geq \frac{N}{N-2}$ for which there exists a unique explicit singular solution $U_M$, and the results present some similarity with the ones of Theorem A.

\medskip

\nind{\bf Theorem B} {\it Let $N\geq 3$, $p\geq \frac N{N-2}$ and $M<0$.  Then \smallskip

\nind 1- If $p\geq\frac{N+2}{N-2}$, then equation $(\ref{I-1})$ admits ground states $u$. Moreover they satisfy $u(r)\sim U_M(r)$ as $r\to\infty$.

\smallskip

\nind 2- If  $\frac{N}{N-2}\leq p<\frac{N+2}{N-2}$, there exist numbers $ \tilde\gm_{min}$ and $\tilde\gm_{max}$ verifying
$$0<\overline\gm<\tilde\gm_{min}\leq \tilde\gm_{max}<\gm^*(1),
$$
such that \\
(i) for $M<-\tilde\gm_{max}$ there exist ground states $u$ such that $u(r)\sim U_{M}(r)$ when $r\to\infty$. \\
(ii) for $M=- \tilde\gm_{min}$ or for $M=- \tilde\gm_{max}$ there exist ground states minimal at infinity in the sense that $u(r)\sim cr^{2-N}$ when $r\to\infty$, $c>0$.
\\
(iii) for $- \tilde\gm_{min}<M<0$ there exists no radial ground state.
}\medskip

Here also the value of $\tilde \gm_{min}$, $\tilde \gm_{max}$ are not explicit and we conjecture that they coincide.
The next result presents some similarity with Theorem A'.
\medskip

\nind{\bf Theorem B'} {\it Let $N\geq 3$ and $\frac N{N-2}<p<\frac{N+2}{N-2}$. \smallskip

\nind (i) If $\overline M<M<0$ there exists a unique (up to scaling) positive singular solution $u$ of $(\ref{I-1})$, such that $u(r)\sim U_M(r)$ when $r\to 0$ and
$u(r)\sim cr^{2-N}$ when $r\to\infty$  for some $c>0$.  \smallskip

\nind (ii) If $- \tilde\gm_{min}<M<0$ there exist positive singular solutions $u$ ondulating around $U_M$ on $[0,\infty)$ and singular solution ondulating around $U_M$ in a neighbourhood of  $0$ and satisfying $u(r)\sim cr^{2-N}$  for some $c>0$ when $r\to\infty$.}\medskip

In terms of the system $(\ref{I-3x})$, (i) corresponds to a heteroclinic orbit connecting $P_M$ and $(0,0)$, while
(ii) to the existence of a periodic solution in {\bf Q} around $P_M$, and the existence of a solution in {\bf Q} converging to $(0,0)$ at $\infty$ and having a limit cycle at $t=-\infty$ which is a periodic orbit around $P_M$.\smallskip

The situation is more complicated when $1<p<\frac{N}{N-2}$ and $M<-\gm^*$ because {\it there exist two explicit singular solutions
$U_{1,M}$ and $U_{2,M}$} which coincide when $M=-\gm^*$.\medskip

\nind{\bf Theorem C} {\it Let $M<0$, $N\geq 3$ and $1<p<\frac{N}{N-2}$, or $N=2$ and $p>1$. Then there exist
two constants $\tilde \gm_{min}$ and $\tilde \gm_{max}$ verifying
$$\gm^*\leq\overline\gm<\tilde \gm_{min}\leq \tilde \gm_{max}<\gm^*(1),
$$
such that \smallskip

\nind 1- If $M<-\tilde  \gm_{max}$ then equation $(\ref{I-1})$ admits ground states $u$ either ondulating around $U_{2,M}$ or such that $u(r)\sim U_{2,M}(r)$ as $r\to\infty$.\smallskip

\nind 2- If $M=-\tilde  \gm_{min}$ or $M=-\tilde  \gm_{max}$ there exists a ground state $u$ such that $u(r)\sim U_{1,M}(r)$ as $r\to\infty$.\smallskip

\nind 3- If $-\tilde  \gm_{min}<M<0$ there exists no radial ground state.
}\medskip

Here again $\tilde \gm_{min}$ and $\tilde \gm_{max}$ appear as transition values for which the ground state still exists but it is smaller than the others at infinity: it behaves like $U_1$ instead of $U_2$. The proof of this theorem is very elaborate in particular in the case $N=2$. In the case $N=1$ the result is already proved in \cite{ChWe}. The nonexistence of a ground state, not necessarily radial for $M>-\gm^*$ is proved in \cite{AGMQ} and independently in \cite{BiGaVe2} with a different method. In the radial case it was obtained much before in the case $N=1$ in \cite{ChWe} and then by Fila and Quittner  \cite{FiQu}
who raised the question whether the condition $-\tilde  \gm_{min}<M<0$ is optimal for the non-existence of radial ground state.  This question received a negative answer in the work of Voirol \cite{Vo} who extended the domain of non-existence to $-\gm^*-\ge<M\leq -\gm^*$. The next result is the counterpart of Theorem C when dealing with singular solutions.
\medskip

\nind{\bf Theorem C'} {\it Let $M<0$, $N\geq 3$ and $1<p<\frac{N}{N-2}$. \smallskip
\nind (i) If $M<-\gm^*$ there exist positive singular solutions $u$ such that $u(r)\sim U_{1,M}(r)$ as $r\to \infty$ and
$u(r)\sim cr^{2-N}$ with $c>0$ when $r\to0$.\smallskip

\nind (ii) If $\overline M\leq M<-\gm^*$ there exists a unique up to scaling positive singular solution $u$, such that $u(r)=U_{2,M}(r)$ as $r\to 0$ and $u(r)=U_{1,M}(r)$ as $r\to \infty$. Furthermore $u(r)>U_{1,M}(r)$ for all $r>0$.
\smallskip

\nind (iii) If $-\hat\gm_{min}<M<-\overline\gm$  there exist positive singular solutions $u$ ondulating around $U_{2, M}$ at $0$ and such that
$u(r)\sim U_{1,M}(r)$ as $r\to \infty$, and positive singular solutions $u$ ondulating around $U_{2, M}$ on $\BBR$.\smallskip

\nind (iv) If $M=-\tilde\gm_{min}$ or $M=-\hat\gm_{max}$ there exists a positive singular solutions $u$ different from $U_{1, M}$ such that
$u(r)\sim U_{1,M}(r)$ when $r\to 0$ and $r\to \infty$.\smallskip

\nind (v) If $-\tilde \gm_{min}<M<-\hat\gm_{max}$ there exists a positive singular solution $u$ such that $\displaystyle \lim_{r\to 0}r^{N-2}u(r)=c>0$ and either ondulating around $U_{2, M}$  or such that $u(r)\sim U_{2,M}(r)$ when $r\to \infty$.

\nind (vi) If $N\geq 3$ and $M=-\gm^*$, there exist positive singular solutions $u$ satisfying $\displaystyle \lim_{r\to 0}r^{N-2}u(r)=c>0$ and $u(r)\sim U_{-\gm^*}(r)$ as $r\to \infty$.
   }\smallskip

In terms of the system $(\ref{I-3x})$ (i) corresponds to a heteroclinic orbit connecting $P_{1,M}$ to $(0,0)$; (ii) to a heteroclinic orbit connecting $P_{2,M}$ to $P_{1,M}$; (iii) to a trajectory having a periodic orbit around $P_{2, M}$ for limit cycle at $-\infty$ and converging to
$P_{2, M}$ at $\infty$ and to a periodic orbit around $P_{2, M}$; (iv) corresponds to homoclinic orbit at $P_{1,M}$; (v) corresponds to a trajectory connecting $(0,0)$ at $-\infty$ and either converging to  $P_{2,M}$ at $\infty$ or having a periodic orbit around $P_{2,M}$ for limit set at $\infty$; (vi) corresponds to a heteroclinic orbit connecting from $(0,0)$ to $P_{-\gm^*}$.



\medskip

\nind{\bf Acknowledgements} This article has been prepared with the support of the collaboration
programs ECOS C14E08 and FONDECYT grant 1160540 for the three authors. \smallskip

The authors thank the anonymous referee for the careful reading of the manuscript which allowed to eliminate some ambiguities in the presentation and the proof of some of our results.  
\mysection{General properties of the system}

\subsection{Reduction to autonomous equation and system}

\subsubsection{The standard reduction}
We recall that if $u$ is a $C^3$ function defined on some interval $I\subset[0,\infty)$ verifying $(\ref{I-1})$
and if
$$u(r)=r^{-\frac{2}{p-1}}x(t)\quad\text{with }\;t=\ln r,
$$
then $x$ satisfies the autonomous equation
       \bel{I-3-2}\BA {lll}\displaystyle
x_{tt}+Lx_{t}-\myfrac{2K}{p-1}x+\abs x^{p-1}x+M\abs{\frac{2x}{p-1}-x_{t}}^{\frac{2p}{p+1}}=0,
             \EA\ee
on $\ln (I)$  where $K$ and $L$ are defined in $(\ref{I-3v})$.   Setting $u_{r}=-r^{-\frac{p+1}{p-1}}y(t)$, then $(x(t),y(t))$ satisfies
            \bel{I-3-4}\BA {lll}\displaystyle
x_{t}=H_1(x,y)\\[2mm]
y_{t}=H_2(x,y),
             \EA\ee
where
             \bel{I-3-4'}\BA {lll}\displaystyle
&H_1(x,y)=\myfrac{2x}{p-1}-y\\[2mm]
&H_2(x,y)=-Ky+\abs x^{p-1}x+M\abs y^{\frac{2p}{p+1}}.
             \EA\ee
             and we denote by {\bf H} the vector field of $\BBR^2$ with components $H_1$ and $H_2$.

\subsubsection{The geometry of the vector field {\bf H}}

Let us denote by ${\bf Q}:=\{(x,y):x>0,y>0\}$ the first quadrant. The vector field is inward in (resp. outward of) ${\bf Q}$  on the axis $\{(x,y):x>0,y=0\}$ (resp. $\{(x,y):x=0, y>0\}$). We set
\bel{I-3-20}\BA{lll}\displaystyle
\CL:=\left\{(x,y)\in{\bf Q}:y=\frac{2x}{p-1}\right\}\quad\text{and } \;\CC:=\left\{(x,y)\in{\bf Q}:x=\left(Ky-My^{\frac{2p}{p+1}}\right)^{\frac{1}{p}}\right\}
\EA\ee
and $\psi(y)=\left(Ky-My^{\frac{2p}{p+1}}\right)^{\frac{1}{p}}$. Then $x_{t}=0$ on $\CL$ and $y_{t}=0$ on $\CC$. The curves $\CL$ and $\CC$ have zero, one or two intersections in {\bf Q} according the value of $K$ and $M$.
If $M,K>0$, then $\CC\subset \left[0,\left(\frac{p-1}{2p}\right)^{\frac 1p}K^{\frac 2{p+1}}\left(\frac{(p+1)}{2pM}\right)^{\frac{p+1}{p(p-1)}}\right]\ti \left[0,\left(\frac{K}{M}\right)^{\frac{p+1}{p-1}}\right]$. The points $(0,0)$, $P_M$ and $\Bigl(0,\left(\frac{K}{M}\right)^{\frac{p+1}{p-1}}\Bigr)$ belong to $\CC$. The function $\psi$ is increasing on $\Bigl(0,\left(\frac{K}{qM}\right)^{\frac{p+1}{p-1}}\Bigr)$ and decreasing on $\Bigl(\left(\frac{K}{qM}\right)^{\frac{p+1}{p-1}},\left(\frac{K}{M}\right)^{\frac{p+1}{p-1}}\Bigr)$.\\
If $M<0$ and $K>0$, $\psi$ is concave and increasing on $(0,\infty)$ with $\psi(y)=(-M)^{\frac{1}{p}}y^{\frac{2}{p+1}}(1+o(1))$ as $y\to\infty$.
If $M<0$ and $K<0$, $\psi$ is still concave and increasing on $\Bigl(\left(\frac{K}{M}\right)^{\frac{p+1}{p-1}},\infty\Bigr)$ with the same asymptotic as above.  We quote below the possible connected components of ${\bf Q}\setminus\left(\CL\cup\CC\right)$.\\

\nind${\bf A}=\left\{(x,y):\frac{2x}{p-1}-y<0\right\}\cap\left\{(x,y):-Ky+x^p+My^{\frac{2p}{p+1}}<0\right\}=\{(x,y):x_t<0,y_t<0\}$.\\
${\bf B}=\left\{(x,y):\frac{2x}{p-1}-y>0\right\}\cap\left\{(x,y):-Ky+x^p+My^{\frac{2p}{p+1}}<0\right\}=\{(x,y):x_t>0,y_t<0\}$.\\
${\bf C}=\left\{(x,y):\frac{2x}{p-1}-y>0\right\}\cap\left\{(x,y):-Ky+x^p+My^{\frac{2p}{p+1}}>0\right\}=\{(x,y):x_t>0,y_t>0\}$.\\
${\bf D}=\left\{(x,y):\frac{2x}{p-1}-y<0\right\}\cap\left\{(x,y):-Ky+x^p+My^{\frac{2p}{p+1}}>0\right\}\cap\left\{(x,y):x>X_{2,M}^{\phantom{r}}\right\}\\\phantom{{\bf D}}=\left\{(x,y):x_t<0,y_t>0\right\}\cap\left\{(x,y):x>X_{2,M}^{\phantom{r}}\right\}$.\\
${\bf E}=\left\{(x,y):\frac{2x}{p-1}-y<0\right\}\cap\left\{(x,y):-Ky+x^p+My^{\frac{2p}{p+1}}>0\right\}\cap\left\{(x,y):x<X_{1,M}^{\phantom{r}}\right\}\\\phantom{{\bf E}}=\{(x,y):x_t<0,y_t>0\}\cap\left\{(x,y):x<X_{1,M}^{\phantom{r}}\right\}$.\smallskip

\nind These connected components are \\
\nind {\bf A}, {\bf B}, {\bf C}, {\bf D} if $K\geq 0$, $M<0$ or $K,M> 0$.\\
\nind  {\bf A}, {\bf C}, {\bf D} if $K<0$ and $-\gm^*<M<0$.\\
\nind {\bf A}, {\bf C}, {\bf D}, {\bf E} if $K<0$ and $M=-\gm^*$.\\
\nind {\bf A}, {\bf B}, {\bf C}, {\bf D}, {\bf E} if $K<0$ and $M<-\gm^*$.
\subsubsection{Graphic representation of the vector field \bf H}
We present below some graphics of the vector field $H$ associated to system $(\ref{I-3x})$.

\begin{figure}[!h]
\begin{center}
 \includegraphics[keepaspectratio, width=10cm]{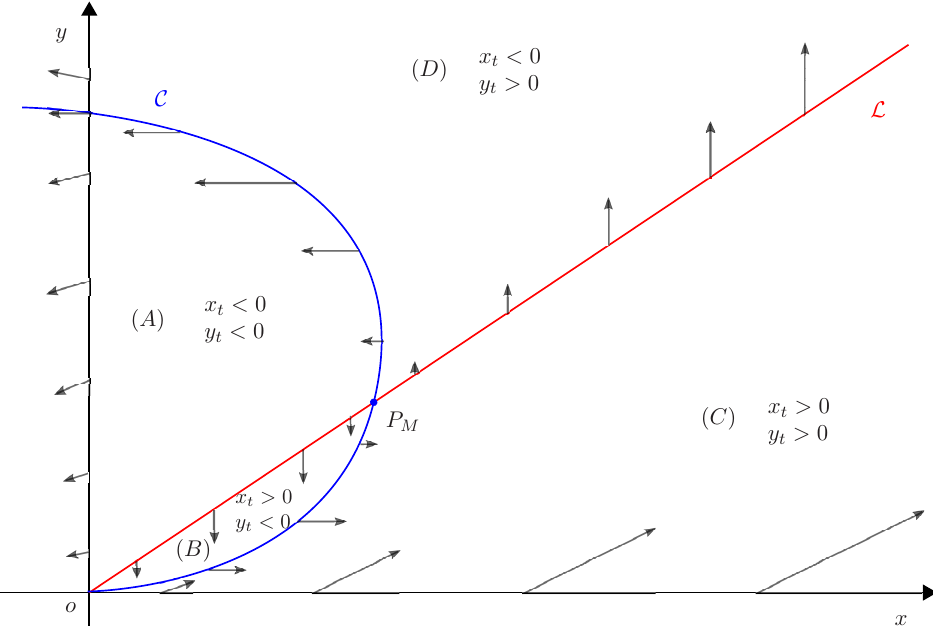}
 \end{center}
 \caption{$M>0$, $K>0$.}
 \end{figure}

\mbox{   }

\begin{figure}[!h]

 \begin{center}
 \includegraphics[keepaspectratio, width=10cm]{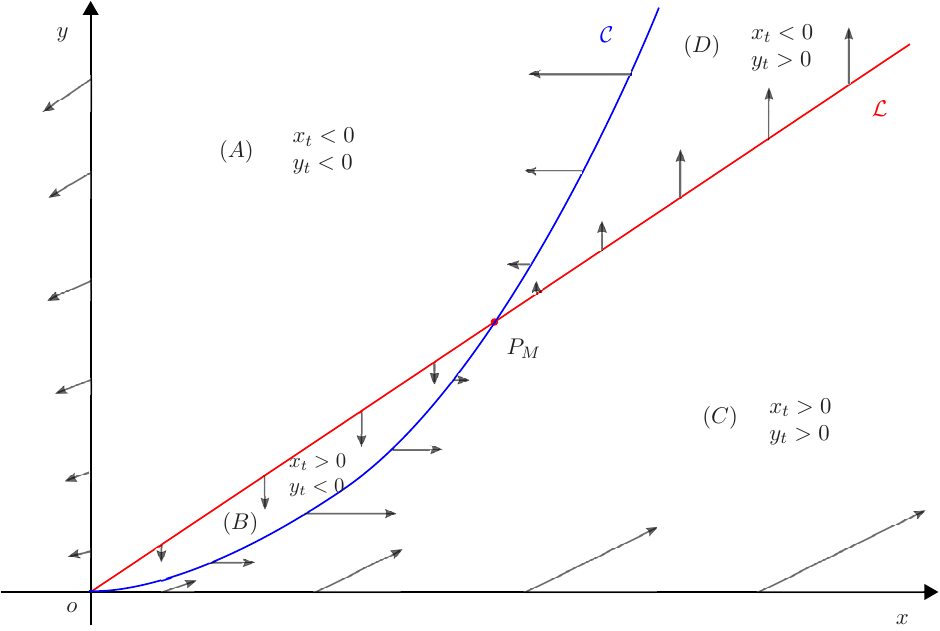}
 \end{center}
 \caption{$M<0$, $K\geq 0$.}
 \begin{center}
\includegraphics[keepaspectratio, width=10cm]{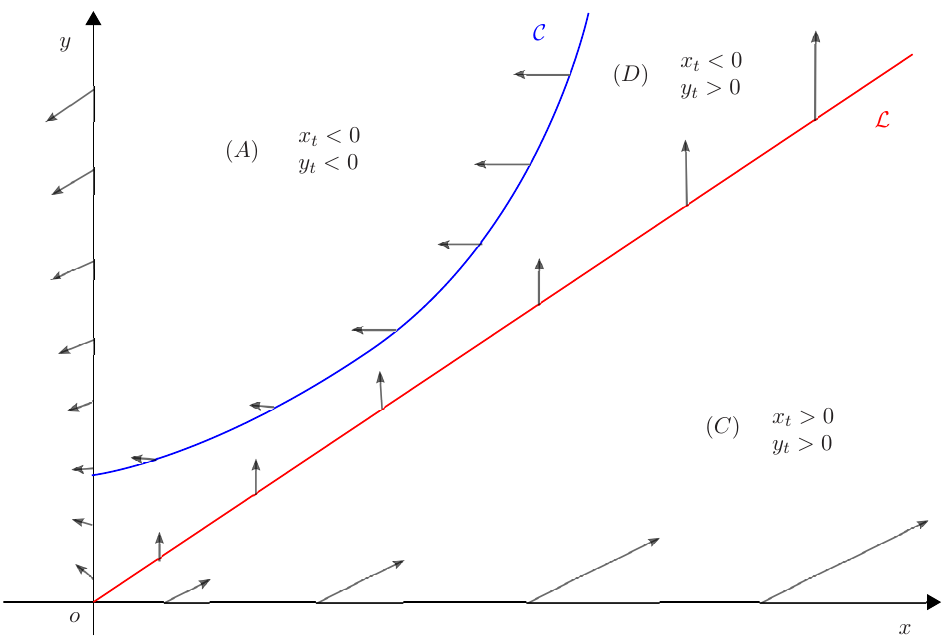}
\end{center}
\caption{$-\gm^*<M<0$, $K<0$.}

 \end{figure}
\medskip
\newpage
\mbox{   }

\begin{figure}[!h]
 \begin{center}
 \includegraphics[keepaspectratio, width=10cm]{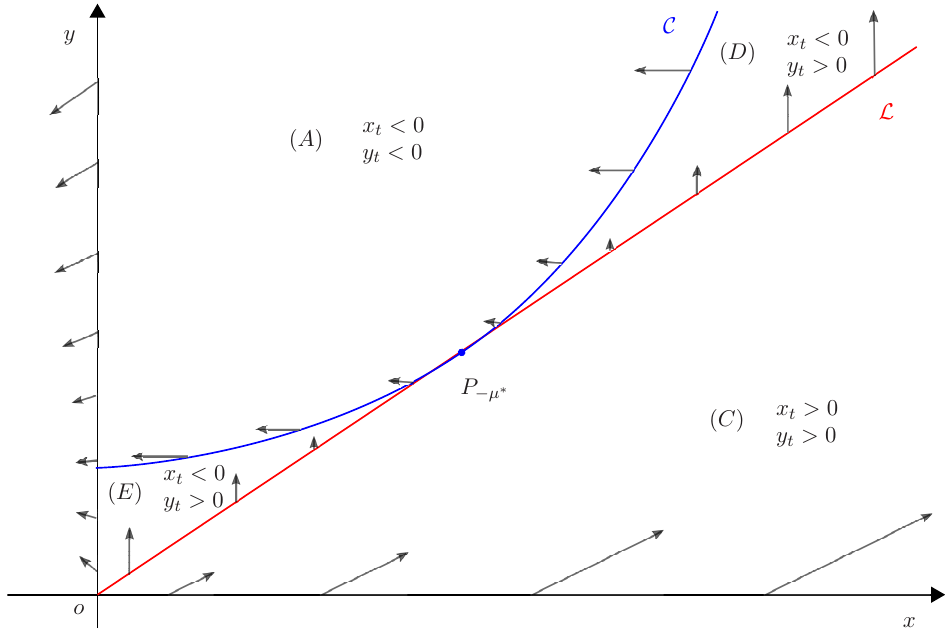}
 \end{center}
 \caption{$M=-\gm^*$, $K<0$.}
 \begin{center}
 \includegraphics[keepaspectratio, width=10cm]{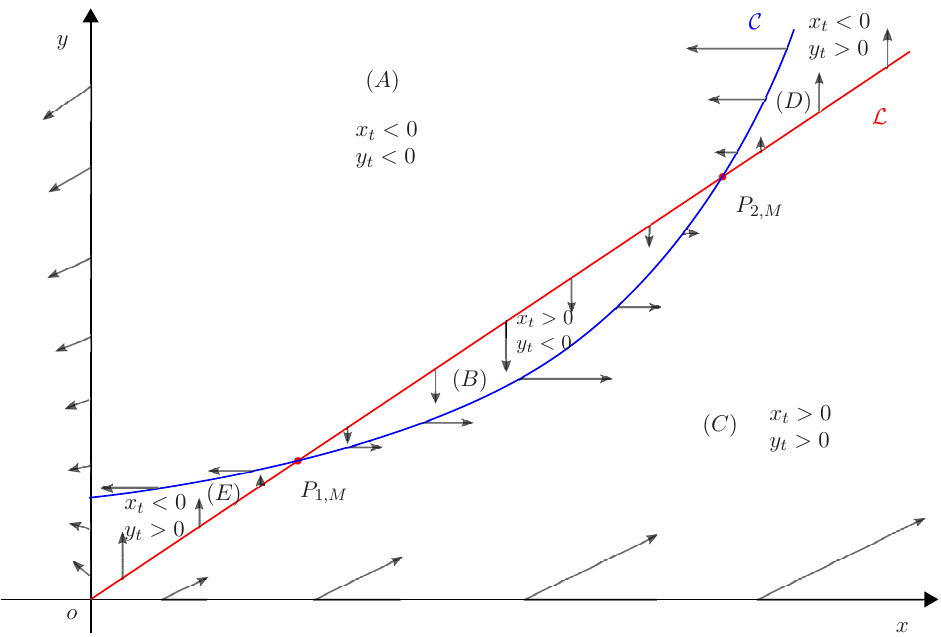}
 \end{center}
 \caption{$M<-\gm^*$, $K<0$.}
 \end{figure}
 \vfill


             \subsubsection{Other reduction}
             The following change of unknowns, already used in \cite{BiGi} when $M=0$,
              \bel{I-3-5}\BA {lll}\displaystyle
\gs(t)=\frac{y(t)}{x(t)}=-r\frac{u_{r}(r)}{u(r)}\quad\text{and }\; z(t)=\frac{\abs x^{p-1}x(t)}{y(t)}=-r\frac{\abs u^{p-1}u(r)}{u_{r}(r)},
             \EA\ee
and valid if $u_{r}\neq 0$, transforms  $(\ref{I-3x})$ into a {\it Kolmogorov system}   with  vector field ${\bf V}=(V_1,V_2)$
          \bel{I-3-6}\BA {lll}\displaystyle
&\gs_{t}=\gs\left(\gs+2-N+z+M\abs{\abs\gs^{p-1}\gs z}^{\frac{1}{p+1}}\right)=V_1(\gs;z)\\[4mm]
&z_{t}=z\left(N-p\gs-z-M\abs{\abs\gs^{p-1}\gs z}^{\frac{1}{p+1}}\right)=V_2(\gs;z).
             \EA\ee
Since $\gs$ and $z$ are in factor the two axis $\{\gs=0\}$ and $\{z=0\}$ are trajectories, actually not admissible for $(\ref{I-3-4})$ in view of $(\ref{I-3-5})$. The system is singular on these two axis however it can be desingularized  by setting $\gs=\tilde \gs^{2k+1}$ and $z=\tilde z^{2k+1}$ for some integer $k>p+1$, which transforms $(\ref{I-3-6})$ into a new nonsingular Kolmogorov system,
          \bel{I-3-6-k}\BA {lll}\displaystyle
&\tilde \gs_{t}=\myfrac{1}{2k+1}\tilde \gs\left(\tilde \gs^{2k+1}+2-N+\tilde z^{2k+1}+M\abs{(\abs{\tilde \gs}^{p-1}\tilde \gs)^{2k+1} \tilde z^{2k+1}}^{\frac{1}{p+1}}\right)\\[4mm]
&\tilde z_{t}=\myfrac{1}{2k+1}\tilde z\left(N-p\tilde\gs^{2k+1}-\tilde z^{2k+1}-M\abs{(\abs{\tilde \gs}^{p-1}\tilde \gs)^{2k+1} \tilde z^{2k+1}}^{\frac{1}{p+1}}\right).
             \EA\ee
Therefore no other trajectory can intersect them in finite time and the quadrant ${\CQ}:=\{(\gs,z):\gs>0,z>0\}$ is invariant. Furthermore $\gs z=r^2\abs u^{p-1}$. It is noticeable that if $M=0$ the initial system is quadratic and regular. \smallskip

   The system   $(\ref{I-3-6})$  will be used in the most delicate cases. It corresponds to the differentiation of the initial equation  $(\ref{I-3-2})$.
\subsection{Regular solutions and ground states}
\bdef{reg-gs} A regular solution of $(\ref{I-1})$ is a $C^{2}$ solution defined on some maximal interval $[0,r_0)$ satisfying $u(0)=u_0>0$ and $u_{r}(0)=0$. A ground state is a nonnegative $C^{2}$ solution defined on $[0,\infty)$.
\es
The existence and uniqueness of a regular solution is standard by the Cauchy-Lipschitz integral method. If $u$ is a $C^2$ solution it satisfies $u_{r}<0$ on $(0,r_0)$. Indeed $r^{N-1}u_{r}(r)$ is decreasing near $0$, hence $u_{r}<0$ on some maximal interval
$(0,r_1)\subset (0,r_0)$ and $u_{r}(r_1)=0$ if $r_1<r_0$. If $u(r_1)=0$ then $u\equiv 0$ by uniqueness. If $u(r_1)>0$ then $u_{rr}<0$ near $r_1$ which would imply that $u(r)<u(r_1)$ for $r_1-\ge\leq r<r_1$ which contradict the negativity of  $u_{r}$ on $(0,r_1)$. Hence $u(r_1)<0$ which implies that $u_{r}(r)<0$ on the maximal interval $(0,r_2)$ where $u>0$. Thus, if $u$ is a ground state $u_{r}<0$ on $(0,\infty)$. Hence the trajectory of a ground state expressed in the system $(\ref{I-3x})$ lies in ${\bf Q}$ and expressed in the system $(\ref{I-3-6})$ it lies in the quadrant ${\CQ}$. It is easy to check that the regular solution $u:=u_{u_0}$, such that $u(0)=u_0$ satisfies
\bel{I-3-8x}
u(r)=u_0\left(1-\myfrac{u_0^{p-1}r^2}{2N}-\myfrac{M(p+1)^2(u_0^{p-1}r^2)^{\frac{2p+1}{p+1}}}{(4p+2)((N+2)p+N)N^q}(1+o(1))\right)\quad\text{as }r\to 0.
\ee
Under the scaling transformation $T_k$, $u$ can be transformed into the regular solution $(\ref{I-1})$ $u:=u_1$ satisfying $u(0)=1$. If one considers the system $(\ref{I-3x})$ the transformation $T_k$ becomes the time shift which transforms $t\mapsto (x(t),y(t))$ into $t\mapsto (x(t+\ln k),y(t+\ln k)) $, and the trajectory $(x(t), y(t))$ of $(\ref{I-3x})$ corresponding to a ground state is therefore uniquely determined and denoted by ${\bf T}_{reg}$ and satisfies
\bel{I-3-8y}\lim_{t\to\infty}(x(t), y(t))=(0,0)\,\text{,  }\lim_{t\to\infty}\myfrac{y(t)}{x(t)}=0\quad\text{and }\lim_{t\to\infty}\myfrac{y(t)}{x^p(t)}
=\frac1N.\ee
Hence in the system $(\ref{I-3-5})$ there holds on the corresponding trajectory
\bel{I-3-8z}\lim_{t\to\infty}\gs(t)=0\,\text{,  }\lim_{t\to\infty}z(t)=N.\ee

\medskip

\subsection{Explicit singular solutions}

Explicit self-similar solutions of $(\ref{I-1})$, necessarily under the form $u=Ar^{-\frac{2}{p-1}}$, play a fundamental role in the study, whenever they exist. The following result covers Proposition 1.

\bprop{fixed-p} 1- Let $M\geq 0$. Then there exists a unique self-similar solution of $(\ref{I-1})$ if and only if $N\geq 3$ and $p>\frac{N}{N-2}$. We denote it by $U_M(r)=X_Mr^{-\frac{2}{p-1}}$, where $X_M>0$ depends also on $M$, $N$ and $p$. To this solution is associated the equilibrium
$(X_M,\frac{2}{p-1}X_M)$ of the system $(\ref{I-3x})$. Furthermore the mappings $M\mapsto X_M$ is continuous and decreasing  on $[0,\infty)$, $M\mapsto MX_M^{\frac{p-1}{p+1}}$
is increasing and there holds
\bel{I-3-8a}\BA {lll}
(i)&\qquad\qquad X_0=\left(\myfrac{2K}{p-1}\right)^{\frac {1}{p-1}},\\[4mm]
(ii)& \myfrac{p-1}{2}\left(\myfrac{K}{M}\right)^{\frac {p+1}{p-1}}\left(1-\myfrac{1}{M}\left(\myfrac{p-1}{2}\right)^p\left(\myfrac{K}{M}\right)^p\right)_+^{\frac{p+1}{p-1}}\leq X_M\leq \myfrac{p-1}{2}\left(\myfrac{K}{M}\right)^{\frac {p+1}{p-1}}.
\EA\ee

 \smallskip

\nind 2- Let $M<0$. If $N\geq 3$ and $p\geq \frac{N}{N-2}$ there exists a unique self-similar solution of $(\ref{I-1})$ $U_M(r)$. The mapping $M\mapsto X_M$ is continuous and decreasing,  $M\mapsto MX_M^{\frac{p-1}{p+1}}$ is decreasing and there holds
\bel{I-3-8ab}\BA {lll}
\max\left\{\left(\myfrac{2K}{p-1}\right)^\frac{1}{p-1},\left(\myfrac{2}{p-1}\right)^\frac{2}{p-1}\abs M^{\frac{p+1}{p(p-1)}} \right\} \leq X_M\\[4mm]
\phantom{---------------}
\leq 2^{\frac{2}{p-1}}\left(\left(\myfrac{2K}{p-1}\right)^\frac{1}{p-1}+\left(\myfrac{2}{p-1}\right)^\frac{2}{p-1}\abs M^{\frac{p+1}{p(p-1)}}\right).
\EA\ee
 \smallskip

\nind 3-  Let $M<0$. If $N=1,2$ and $p>1$, or $N\geq 3$ and $1<p< \frac{N}{N-2}$, there exists no self-similar solution of $(\ref{I-1})$ if $-\gm^*<M<0$ where $\gm^*=\gm^*(N)>0$ is defined in  $(\ref{I-6})$. If
$M=-\gm^*$ there exists a unique self-similar solution $U_{-\gm^*}(r)$. If $M<-\gm^*<0$ there exist two-self-similar solutions
$U_{1,M}(r)=X_{1,M}r^{-\frac{2}{p-1}}$ and $U_{2,M}(r)=X_{2,M}r^{-\frac{2}{p-1}}$ with $X_{1,M}<X_{2,M}$. Furthermore the mappings $M\mapsto X_{1,M}$ and $M\mapsto X_{2,M}$ are continuous, respectively increasing and decreasing on $(-\infty,-\gm^*)$, while $M\mapsto MX^{\frac{p-1}{p+1}}_{1,M}$ and $M\mapsto MX^{\frac{p-1}{p+1}}_{2,M}$ are  respectively decreasing and increasing. Furthermore there holds  for $\abs{M}$ large enough;
\bel{I-3-8b}\BA {lll}
(i)&\qquad\qquad X_{j,-\gm^*}=\left(\myfrac{2}{p-1}\right)^{\frac{p}{p-1}}\left(\myfrac{-K}{p}\right)^{\frac{1}{p-1}},\\[4mm]
(ii)&\myfrac{p-1}{2} \left(\myfrac{K}{M}\right)^{\frac{p+1}{p-1}}<X_{1,M}<\myfrac{p-1}{2} \left(\myfrac{K}{M}\right)^{\frac{p+1}{p-1}}
\left(1-\myfrac{2}{K}\left(\myfrac{(p-1)K}{2M}\right)^{p+1}\right)^{\frac{p+1}{p-1}},\\[4mm]
(iii)&\left(\myfrac{2}{p-1}\right)^{\frac{2}{p-1}}(-M)^{\frac{p+1}{p(p-1)}}\left(1-\myfrac{K}{M\abs M^{\frac 1p}}\right)^{\frac{p+1}{p-1}}<X_{2,M}<\left(\myfrac{2}{p-1}\right)^{\frac{2}{p-1}}(-M)^{\frac{p+1}{p(p-1)}}.
\EA\ee
\es

\nind\Proof  The function $U_M=X_Mr^{-\frac 2{p-1}}$ is a self-similar solution of $(\ref{I-1})$ if and only if $X_M$ is a positive root of
  \bel{I-3-9-}\tilde f_M(x):=x^{p-1}+\left(\frac{2}{p-1}\right)^{\frac{2p}{p+1}}Mx^{\frac{p-1}{p+1}}-\frac{2K}{p-1}=0.
\ee
Equivalently $P_M=(X_M,Y_M)=(\frac{p-1}{2}Y_M,Y_M)$ is a fixed point of system $(\ref{I-3x})$, where $Y_M$ is the positive root of
            \bel{I-3-9}\BA {lll}\displaystyle
f_M(y)=\left(\frac{p-1}{2}\right)^py^{p-1}+My^{\frac{p-1}{p+1}}-K=0.
             \EA\ee
The use of the variable $y$ is a little easer than $x$. Since $X_0$ is explicit if $M=0$, we shall study the cases $M\neq 0$. \smallskip

 \nind 1- {\it Case $M>0$}. If $K>0$, equivalently $p>\frac{N}{N-2}$, and $M\geq 0$, $f_M$ is an increasing function tending to $\infty$ at $\infty$ and negative at $y=0$. Hence $Y_M$ is the unique positive root of   $(\ref{I-3-9})$. If $K<0$, $f_M$ is positive on $[0,\infty)$, hence no such solution exists. Since
 $$\left(\frac{p-1}{2}\right)^pY_M^{p-1}+MY_M^{\frac{p-1}{p+1}}-K=0,
$$
by the implicit function theorem, $M\mapsto Y_M$ is $C^1$. For $M>M'>0$, $f_M(y)>f_{M'}(y)$ for all $y>0$.
Hence $M\mapsto Y_M$ is decreasing on $[0,\infty)$. Actually the expression of $f_M$ shows more, namely that $M\mapsto MY_M^{\frac{p-1}{p+1}}$
is increasing on $[0,\infty)$. Furthermore
$$MY_M^{\frac{p-1}{p+1}}<K\Longrightarrow Y_M<\left(\frac{K}{M}\right)^{\frac{p+1}{p-1}},
$$
and
$$\BA {llll}\left(\myfrac{p-1}{2}\right)^{p}\left(\myfrac{K}{M}\right)^{p+1}+MY_M^{\frac{p-1}{p+1}}>K\\[4mm]
\phantom{---------------}
\Longrightarrow Y_M>\left(\myfrac{K}{M}\right)^{\frac{p+1}{p-1}}\left(1-\left(\myfrac{p-1}{2}\right)^p\myfrac{1}{M}\left(\myfrac{K}{M}\right)^{p}\right)_+^{\frac{p+1}{p-1}},
\EA$$
and we get $(\ref{I-3-8a})$.
 \smallskip

 \nind 2- {\it Case $M<0$}. Clearly $f_M$ has a minimum at $y=y_{0,M}$ where
\bel{1-3-10a}y_{0,M}=\left(\frac{2}{p-1}\left(\frac{-M}{p+1}\right)^{\frac{1}{p}}\right)^{\frac{p+1}{p-1}}\quad\text{and }\;f_M(y_{0,M})=
- \frac{2p}{p-1}\left(\frac{-M}{p+1}\right)^{\frac{p+1}{p}}-K.
\ee
We encounter two possibilities:\smallskip

\nind 2-a.  If $N\geq 3$ and $p\geq\frac{N}{N-2}$, then $f_M(0)<0$, $f_M$ is decreasing on $(0,y_{0,M})$ and increasing on $(y_{0,M},\infty)$, hence $f_M(y)=0$ has a unique positive root $Y_M>y_{0,M}$. Since for $M<M'<0$ and $y>0$,
$f_{M'}(y)>f_M(y)$, the mapping $M\mapsto Y_M$ is continuous and decreasing and $M\mapsto MY^{\frac{p-1}{p+1}}_M$ is increasing.\\
Then the left-hand side of $(\ref{I-3-8ab})$ is clear. Next we put
$$A^p=\left(\myfrac{2}{p-1}\right)^{\frac{2p}{p+1}}\abs M\,,\; B^{p+1}=\myfrac{2K}{p-1}\,,\; \xi=X^{\frac{p-1}{p+1}}\,,\;
a=\myfrac{B}{A}\,\text{ and }\;\eta=\myfrac{\xi}{A}.
$$
Then $\phi(\eta)=\eta^{p+1}-\eta-a^{p+1}=0$. Since
$$\BA {lll}\phi(1+a)=(1+a)^{p+1}-1-a-a^{p+1}\geq (1+a)^{p+1}-1-a-a^p-a^{p+1}\\
\phantom{\phi(1+a)}
\geq (1+a)\left((1+a)^{p}-(1+a^p)\right)\geq 0
\EA$$
we derive $\eta\leq 1+a$, which implies
the right-hand side of $(\ref{I-3-8ab})$.\smallskip

\nind  2-b. If $N=1,2$ or $N\geq 3$ and $1<p<\frac{N}{N-2}$, then $K<0$. Hence, if $f_M(y_{0,M})>0$, or equivalently $-\gm^*<M\leq 0$, the equation $f_M(y)=0$ has no positive root, if $M=-\gm^*$, it has a double root $Y_{-\gm^*}$ where
            \bel{I-3-10}\BA {lll}\displaystyle
Y_{-\gm^*}=\left(\frac 2{p-1}\right)^{\frac{p}{p-1}}\left(\frac {-K}{p}\right)^{\frac{1}{p-1}}\Longrightarrow \gm^*Y_{-\gm^*}^{\frac{p-1}{p+1}}=-\frac{p+1}{p}K,
 \EA\ee
and if $M<-\gm^*$ the equation $(\ref{I-3-9})$ has two positive roots $0<Y_{1,M}<Y_{2,M}$ and since $f'_M$ does not vanish at $Y_{j,M}$, they are $C^1$ functions of $M\in (-\infty,M^*)$, respectively increasing and decreasing. Since $M,K<0$, we obtain from $f_M(Y_M)=0$,
$$Y_{1,M}>\left(\myfrac{K}{M}\right)^{\frac{p+1}{p-1}},
$$
and
$$\left(\frac{p-1}{2}\right)^pY_{2,M}^{p-1}< -MY_{2,M}^{\frac{p-1}{p+1}}\Longrightarrow Y_{2,M}<\left(\frac{2}{p-1}\right)^{\frac{p+1}{p-1}}(-M)^{\frac{p+1}{p(p-1)}}
$$
For a sharper estimate, we have for $M$ large enough,
$$\tilde Y=\left(\myfrac{K}{M}\right)^{\frac{p+1}{p-1}}\left(1-\frac{2}{K}\left(\myfrac{(p-1)K}{2M}\right)^{p+1}\right)^{\frac{p+1}{p-1}}\Longrightarrow f_M(\tilde Y)<0.
$$
Hence
$$\left(\myfrac{K}{M}\right)^{\frac{p+1}{p-1}}<Y_{1,M}<\left(\myfrac{K}{M}\right)^{\frac{p+1}{p-1}}\left(1-\frac{2}{K}\left(\myfrac{(p-1)K}{2M}\right)^{p+1}\right)^{\frac{p+1}{p-1}}.
$$
Similarly
$$\left(\frac{2}{p-1}\right)^{\frac{p+1}{p-1}}(-M)^{\frac{p+1}{p(p-1)}}\left(1-\myfrac{K}{M\abs M^{\frac 1p}}\right)^{\frac{p+1}{p-1}}<Y_{2,M}<\left(\frac{2}{p-1}\right)^{\frac{p+1}{p-1}}(-M)^{\frac{p+1}{p(p-1)}}.
$$
The estimates $(\ref{I-3-8b})$ follow.
\qeda

\subsubsection{Upper estimate of the regular solutions}
We first recall the following estimate in the case $M\geq 0$,  consequence of the fact that the positive solutions of $(\ref{I-0})$ are superharmonic and proved in a more general setting in \cite[Prop. 2.1 ]{BiGaVe2}.

\bprop{zero} 1- There exists no positive solution of $(\ref{I-1})$ in $(R,\infty)$, $R\geq 0$ if  $M\geq 0$ and either $N=1,2$ and $p>1$ or $N\geq 3$ and $1<p\leq \frac{N}{N-2}$. In particular there exists no ground state.\smallskip

\nind 2- If $N\geq 3$, $p>\frac{N}{N-2}$, $M\geq 0$ and $u$ is a positive solution of $(\ref{I-1})$ in $(R,\infty)$, $R\geq 0$. Then  there exists $\gr\geq R$ such that
\bel{I-1-01}
 u(r)\leq \min \left\{\left(\frac{2N}{p-1}\right)^{\frac{1}{p-1}},\myfrac{p-1}{2}\left(\myfrac{p(N-2)-N}{(p-1)M}\right)^{\frac{p+1}{p-1}}\right\}r^{-\frac{2}{p-1}},\qquad\text{for all }\;r>\gr ,
\ee
and
\bel{I-1-01'}
\abs{u_{r}(r)}\leq \min \left\{(N-2)\left(\frac{2N}{p-1}\right)^{\frac{1}{p-1}},\left(\myfrac{p(N-2)-N}{(p-1)M}\right)^{\frac{p+1}{p-1}}\right\}r^{-\frac{p+1}{p-1}}\qquad\text{for all }\;r>\gr.
\ee
 Furthermore, if $R=0$, inequalities $(\ref{I-1-01})$,  and $(\ref{I-1-01'})$ hold with $\gr=0$.
\es

The next estimate is verified by any ground state, independently of the sign of $M$.

\bprop{bound} Let $p>1$ and $N\geq 1$. Then the ground state  $u$ of $(\ref{I-1})$ with $u(0)=1$  satisfies
\bel{I-3-21}\BA{lll}\displaystyle
u(r)\leq \min\left\{1,c_{N,p,M}r^{-\frac{2}{p-1}}\right\}\,\text{ and }\;\abs{u_{r}(r)}\leq c_{N,p,M}r^{-\frac{p+1}{p-1}}\quad\forall r>0.
\EA
\ee
\es
\Proof
\nind The trajectory ${\bf T}_{reg}$ starts from $(0,0)$ and enters the region ${\bf C}$. If it stays in ${\bf C}$, then $x$ is increasing on
$\BBR$. Since $y<\frac{2x}{p-1}$, if $x(t)$ tends to some finite limit as $t\to\infty$, it implies that the limit set of the trajectory exists. It cannot be a cycle since $x$ is monotone, thus it is one of the equilibrium of the system. Hence
$x(t)\leq X_{j,M}$ for $j=1$ or $2$ if $K<0$ and $M<-\gm^*$  or $x(t)\leq X_M$ if $K\geq 0$ and either $M>0$ or $-\gm^*\leq M\leq 0$. This implies that $(\ref{I-3-21})$ holds.
If $x(t)$ tends to $\infty$, so does $y$ because $y(t)\geq x^p(t)-\frac{2|K|}{p-1}x(t)-|M|\left(\frac{2x(t)}{p-1}\right)^{\frac{2p}{p+1}}\to\infty$ as $t\to\infty$. Therefore $y_{t}\geq Cy^p$ for some $C>0$ which would imply that $t\mapsto y^{1-p}(t)+C(p-1)t$ is increasing, which is impossible. \smallskip

Next we suppose that the trajectory leaves ${\bf C}$ by crossing the line $\CL$. Since it cannot enter ${\bf B}$ through $\CC$ (in the case $K<0$, $M\leq -\gm^*$ ), it leaves  ${\bf C}$ by intersecting $\CL$ and we denote by $t_1$ the first time where ${\bf T}_{reg}$ intersects $\CL$. Then
$x_{t}(t_1)=0$ and $x_{tt}(t_1)=-y_{t}(t_1)<0$. Therefore $t_1$ is a local maximum.  Now the trajectory cannot cross again the half line
$\left\{(x,y):x=x(t_1), y>\frac{2x(t_1)}{p-1}\right\}$ because on it there holds $x_{t}=\frac{2x(t_1)}{p-1}-y<0$. Hence $x(t)\leq x(t_1)$
for any $t\geq t_1$. \smallskip

In the same way, either $y$ is increasing on $\BBR$ and since $x_{t}=\frac{2}{p-1}x-y$ and $x$ is bounded, $y$ cannot tend to infinity when $t\to\infty$, thus $y(t)\to y_0>0$, or $y$ is not monotone and ${\bf T}_{reg}$ crosses $\CC$ at a first value $t_2$, necessarily larger than $t_1$ and where $x_{t}(t_2)<0$.  Then $y_{tt}(t_2)\!=\!p\abs{x(t_2)}^{p-1}\!x_{t}(t_2)<0$ and $t_2$ is a local maximum of $y$. Therefore $\cup_{t\geq t_2}(x(t),y(t))$ remains in the subset of ${\bf Q}$ bordered by the portion of trajectory of ${\bf T}_{reg}$ for $t\leq t_2$ and $\{(x,y):0<x\leq x(t_2), y=y(t_2)\}$. This implies that $y(t)\leq y(t_2)$ for all $t\in\BBR$. Noticing that $u(r)\leq 1$ since $u$ is decreasing, we get the conclusion.

\qeda
\medskip

\nind\Remark The above method does not give an explicit estimate of the upper bounds of $x$ and $y$ and such a bound can be estimated in some cases.
If $M\geq 0$   it follows from \rprop{zero} that for any $p>1$ there holds
$$u(r)\leq \left(\frac{2N}{p-1}\right)^{\frac{1}{p-1}}r^{-\frac{2}{p-1}},
$$
thus this estimate is independent of $M$. Here a new critical value is involved in all dimension $N\geq 2$, namely 
\bel{I-3-222}
\gm^*(2)=(p+1)\left(\myfrac{1}{p}\right)^{\frac{p}{p+1}},
\ee
corresponding to the definition $(\ref{I-6})$.  If $-\gm^*(2)<M<0$, it is easy to check that the function $v=\ln u$ satisfies
\bel{I-3-22}\BA{lll}\displaystyle
-\Gd v\geq ae^{(p-1)v}
\EA
\ee
with $a=1-\left(\frac{\abs M}{\gm^*(2)}\right)^{p+1}>0$. We derive that the function
$w(r)=-(r^{N-1}v_r(r))$ is increasing, with limit $\ell\in (0,\infty]$. Hence
$$e^{(1-p)v(r)}\geq e^{(1-p)v(0)}+\frac{a(p-1)r^2}{4N}
$$
This implies
\bel{I-3-24}u(r)\leq \min\left\{1,c_0a^{-\frac{1}{p-1}}r^{-\frac{2}{p-1}}\right\}.
\ee

\subsection{Linearization of the system $(\ref{I-3x})$ near equilibria}

\subsubsection{Linearization at $(0,0)$}
The linearized system at $(0,0)$ is
            \bel{I-3-11}\BA {lll}\displaystyle
x_{t}=\myfrac{2x}{p-1}-y\\[2mm]
y_{t}=-Ky,
 \EA\ee
 which admits the eigenvalues
  \bel{I-3-11a}\gl_1=-K\qquad \gl_2=\frac{2}{p-1}.
  \ee
  Note that $\gl_2-\gl_1=N-2$.
 \smallskip

 \nind (a) Assume that $N\geq 3$ and  $p>\frac{N}{N-2}$,  or equivalently $K>0$. Then $(0,0)$ is a saddle point. There exists a unique {\it unstable trajectory} ${\bf T}_{unst}$ such that
   \bel{I-3-11b}\displaystyle
    \lim_{t\to-\infty}(x(t), y(t))=(0,0) \quad\text{and }\; \lim_{t\to-\infty}\myfrac{y(t)}{x(t)}=0,
  \ee
 and more precisely, from $(\ref{I-3-8x})$,
 \bel{I-3-11w}\BA {lll}x(t)=e^{\frac{2t}{p-1}}\left(1-\myfrac{e^{2t}}{2N}-\myfrac{M(p+1)^2e^{\frac{4p+2}{p+1}t} }{2\left(N(p+1)+2p\right)(2p+1)} (1+o(1))\right)\\[5mm]
y(t)=e^{\frac{2pt}{p-1}}\left(\myfrac1N+\myfrac{M(p+1)}{N(p+1)+2p}e^{\frac{(p-1)t}{p+1}}(1+o(1))\right)\quad\text{as }\,t\to-\inftyÊ.
\EA\ee
 From Definition 2.1 and the lines which follow, the unstable trajectory ${\bf T}_{unst}$ coincides with the regular trajectory ${\bf T}_{reg}$.   This is included in the region {\bf C} for $t<T_0$ for some $-\infty<T_0\leq\infty$.  There exists also a unique {\it stable trajectory} ${\bf T}_{st}$ such that
    \bel{I-3-11c}\displaystyle
    \lim_{t\to\infty}(x(t), y(t))=(0,0) \quad\text{and }\; \lim_{t\to-\infty}\myfrac{y(t)}{x(t)}=\myfrac{2}{p-1}+K=N-2.
  \ee
  Since $N-2>\frac2{p-1}$, ${\bf T}_{st}$ belongs to the region {\bf A} for $t>T_1$ for some $T_1<\infty$. \medskip

   \nind\Remark If ${\bf T}_{reg}\subset{\rm\bf Q}$ satisfies $\displaystyle\lim_{t\to\infty}(x(t),y(t))=(0,0)$, then ${\bf T}_{reg}={\bf T}_{st}$ and the corresponding solution is a ground state. The same conclusion holds, if ${\bf T}_{st}\subset{\rm\bf Q}$ satisfies $\displaystyle\lim_{t\to-\infty}(x(t),y(t))=(0,0)$. Such a solution is called a {\it homoclinic orbit at $(0,0)$}. Because of the uniqueness of the stable and unstable trajectories of a saddle point, it is unique in the class of solutions satisfying $(\ref{I-3-11c})$. Equivalently the class of ground states $u$ of $(\ref{I-1})$ satisfying
 $u(r)\sim cr^{2-N}$ for some $c>0$ is then a one parameter family characterized by $u(0)=u_0$.

 \smallskip

 \nind (b) Assume that $N\geq 3$ and $1<p<\frac N{N-2}$. Then $K <0$ and $0<\gl_1<\gl_2$. Hence $(0,0)$ is a source and all the trajectories of $(\ref{I-3x})$ in some neighbourhood of $(0,0)$  converge to $(0,0)$ when $t\to-\infty$. Among those trajectories there exists {\it one fast trajectory} which satisfies
  $$\lim_{t\to-\infty}\frac{y(t)}{x(t)}=0.
  $$
 It is actually the regular trajectory ${\bf T}_{reg}$. There exist also infinitely many {\it slow trajectories} which satisfy
    $$\lim_{t\to-\infty}\frac{y(t)}{x(t)}=N-2.
  $$
   \smallskip

 \nind (c) If $p=\frac N{N-2}$, then $K =0$ and $\gl_1=0<\gl_2=N-2$. We still find the regular trajectory   ${\bf T}_{reg}$ associated to $\gl_2$ and the corresponding eigenvector $(1,0)$.
By the central manifold theorem corresponding to $\gl_1$ there exists an invariant curve passing through $(0,0)$ with slope $N-2$. Using the matched asymptotic expansion method, one finds that if $M<0$ there exists a solution $x$ of $(\ref{I-3-2})$ such that
$x(t)\sim C_N\abs M^{-\frac{1}{N-1}}t^{1-N}$ when $t\to \infty$, i.e. $u(r)\sim C_N\abs M^{-\frac{1}{N-1}}r^{2-N}(\ln r))^{1-N}$ when $r\to \infty$, and if $M>0$
 there exists a solution $x$ of $(\ref{I-3-2})$ such that
$x(t)\sim C_NM^{-\frac{1}{N-1}}(-t)^{1-N}$ when $t\to -\infty$, equivalently $u(r)\sim C_N\abs M^{-\frac{1}{N-1}}r^{2-N}(\ln (\frac 1r))^{1-N}$ when $r\to 0$.\smallskip

 \nind (d) If $N=2$, $(0,0)$ is a source with $\gl_1=\gl_2=\frac{2}{p-1}$, with corresponding eigenspace $(1,0)$. The linearized problem is equivalent to equation
 $$x_{tt}-\myfrac{4}{q-1}x_{t}+\myfrac{4}{(q-1)^2}x=0
 $$
 with general solutions $x(t)=ae^{\frac{2t}{q-1}}+bte^{\frac{2t}{q-1}}$ for  some real parameters $a,b$.
 Hence there exists infinitely trajectories of $(\ref{I-3-2})$ tending to $0$ when $t\to-\infty$ and they are tangent to $(1,0)$ at $(0,0)$. The regular trajectory  ${\bf T}_{reg}$ corresponds to $b=0$ and the other trajectories correspond to singular solutions $u$ of $(\ref{I-0})$. They satisfy
 $u(r)\sim b\ln r$ as $r\to 0$ and there holds
 $$-\Gd u=u^p+M\abs{\nabla u}^{\frac{2p}{p+1}}-2\gp b\gd_{0},
 $$
 in the sense of distributions in $B_\ge$ for some $\ge>0$.\medskip

 Next we give a general result in case the system admits only one equilibrium in {\bf Q}.

 \blemma{tetra} Let $N\geq 3$, $p>\frac{N}{N-2}$ and $M\in\BBR$. If $u$ is a regular solution the following tetrachotomy occurs:\\
 (i) either $\lim_{r\to\infty}r^{N-2}u(r)=c$ for some $c>0$,\\
 (ii) or $u(r)\sim U_M(r)$ as $r\to\infty$,\\
 (iii) or $u(r)$ has an $\omega$-limit cycle surrounding $P_M$,\\
 (iv) or $u(r)$ changes sign for some $r>0$.
 \es
 \Proof By assumption $P_M$ is the unique equilibrium. The trajectory {\bf T}$_{reg}$ starts from $(0,0)$ and remains in the region {\bf C} where $x_t,y_t>0$ for $t\leq t_0\leq\infty$. If $t_0=\infty$, $u$ is a ground state, hence it is bounded from \rprop{bound}. Its $\omega$-limit set is non-empty. Because $x$ and $y$ are monotone, it converges when $t\to\infty$ to some point  which is necessarily $P_M$.  If $t_0<\infty$, then at $t=t_0$ the trajectory leaves {\bf C}  through $\CL$ since it cannot enter in {\bf B}, and it enters the region {\bf D} where $x_t<0$, $y_t>0$. Moreover $x_t(t_0)=0$ and $x(t_0)>X_M$. Then three possibilities occur:\\
 ($\ga$) either $x(t)\to X_M$ monotonically when $t\to\infty$; thus the trajectory converges to $P_M$. \\
 ($\gb$) either $x(t)\to 0$ monotonicaly  when $t\to\infty$. Since $(0,0)$ is a saddle point, then {\bf T}$_{reg}$={\bf T}$_{st}$. This implies that {\bf T}$_{reg}$ is a homoclinic trajectory at $(0,0)$.\\
($\gg$) or there exists $t_1>t_0$ such that $x_t(t_1)=0$. Then $x(t_1)<X_M$. Hence {\bf T}$_{reg}$ enters the region {\bf B} and by continuity there exists $t'<t_0$ such that $x(t')=x(t_1)$ and $y(t')<y(t_1)$. Therefore the bounded region of $\BBR^2$ bordered by the segment $I=\{(x,y):x=x(t_1), y(t')< y< y(t_1)\}$ and the portion of {\bf T}$_{reg}$ defined by $\{(x(t),y(t))\in {\bf T}_{reg}:t'\leq t\leq t_1\}$ is positively invariant (notice that $x_t>0$ on $I$) and it contains $P_M$ and no other equilibrium. Therefore either the trajectory converges to $P_M$ or it admits an $\gw$-limit cycle which is a closed orbit surrounding $P_M$.
\qeda

\subsubsection{Linearization at a fixed point $P_M:=(X_M,Y_M)$}
Suppose that $P_M$ (or  $P_{j,M}$) exists. Then setting $(x,y)=(X_M,Y_M)+(\overline x,\overline y)$, the linearized system at this point is

            \bel{I-3-12}\BA {lll}\displaystyle
\overline x_{t}=\frac{2\overline x}{p-1}-\overline y\\[2mm]
\overline y_{t}=pX_M^{p-1}\overline x+\left(\frac{2p}{p+1}MY_M^{\frac{p-1}{p+1}}-K\right)\overline y,
 \EA\ee
Using equation $(\ref{I-3-9})$, the eigenvalues of its matrix are the roots of the trinomial
             \bel{I-3-13}\BA {lll}\displaystyle
T(\gl)=\gl^2-\left(\frac{2p}{p+1}MY_M^{\frac{p-1}{p+1}}-L\right)\gl+2K-\frac{2p}{p+1}MY_M^{\frac{p-1}{p+1}}.
 \EA\ee
  If $M$ is such that $\frac{2p}{p+1}MY_M^{\frac{p-1}{p+1}}=L$, then $2K-\frac{2p}{p+1}MY_M^{\frac{p-1}{p+1}}=2K-L=K+\frac{2}{p-1}=N-2$, and we denote by $\overline M$ such a value of $M$ which  is characterized by
              \bel{I-3-13'}\BA {lll}\displaystyle
\myfrac{2p}{p+1} {\overline M}Y_M^{\frac{p-1}{p+1}}=L.
 \EA\ee
 Since $Y_M$ is a positive root of $(\ref{I-3-9})$, we get
 $$\left(\frac{2}{p-1}\right)^pY_M^{p-1}=A:=K-\frac{p+1}{2p}L=\frac{(N-2)(p-1)^2+4}{2p(p-1)},
 $$
 hence  $\overline M$, well defined for $N\geq 2$, is given by
\bel{I-3-14}\overline M=\frac{(p+1)L}{2p}\left(\frac 1A\left(\frac{p-1}{2}\right)^p\right)^{\frac{1}{p+1}}
 =\frac{(p+1)\left((N-2)p-N-2\right)}{(4p)^\frac{p}{p+1}\left((N-2)(p-1)^2+4\right)^{\frac{1}{p+1}}}.
\ee
That is $(\ref{I-8})$.
\medskip

 \nind\Remark We see that $\overline M>0$ (resp. $\overline M<0$)  if and only if $N\geq 3$ and $p>\frac{N+2}{N-2}$ (resp. $1<p<\frac{N+2}{N-2}$ if $N\geq 3$ and any $p>1$ if $N=2$). If $M=\overline M$ and $N\geq 3$, the eigenvalues of the linearized system are purely imaginary. If $p=\frac{N+2}{N-2}$, then $\overline M=0$. It is known that in that case the point $P_0=\left(\left(\frac{N-2}{2}\right)^{\frac{N-2}{2}},\left(\frac{N-2}{2}\right)^{\frac{N}{2}}\right)$
is a center for the system $(\ref{I-3x})$ associated to  $-\Gd u=u^{\frac{N+2}{N-2}}$  and that there exist infinitely many cycles turning around $P_0$ with equation
\bel{sob}\frac 12y^2+\frac {N-2}{2N}x^{\frac{2N}{N-2}}-\frac{N-2}{2}xy=E.
\ee
 It can be verified that $E<0$, in particular using the function $F$ defined in $(\ref{I-3-15})$.
 For $N\geq 2$ and $1<p<\frac{N}{N-2}$, there always holds $\overline M\leq -\gm^*$ and more precisely

\blemma{comp} If $N\geq 3$ and $1<p<\frac{N}{N-2}$, then $\overline M<-\gm^*$. If $N=2$ and $p>1$, then $\overline M=-\gm^*$.
\es
\Proof From \rprop{fixed-p} and identities $(\ref{I-3-10})$ and $(\ref{I-3-13})$, there holds,
$$\gm^*<\overline\gm\Longleftrightarrow \gm^*Y_{-\gm^*}^{\frac{p-1}{p+1}}<\overline\gm Y_{\overline M}^{\frac{p-1}{p+1}}\Longleftrightarrow2\abs K<\abs L.
$$
But $\abs L-2\abs K=\frac{N-2}{p-1}$ and the conclusion follows. If $N=2$ we just replace $<$ by $=$ in the above series of equivalences.
\qeda

 \blemma{barM} Assume $N\geq 3$. If $p>\frac{N+2}{N-2}$ and $M=\overline M>0$ then  $P_M$ is a weak sink and a Andronov-Hopf bifurcation point. If $p<\frac{N+2}{N-2}$ and $M=\overline M<0$ then $P_M$ is a weak source.
 \es
 \Proof We recall that a weak sink is an asymptotically stable equilibrium which attracts the nearby points as $t\to\infty$ at a rate slower than the usual exponential rate. A weak source is a weak sink of the system obtained by changing $t$ into $-t$ (see \cite[Chap. 9]{HuWe}). We write
 $\bar x=x-X_M$, $\bar y=y-Y_M$ and obtain the new nonlinear system
\bel{I-3-14a}\BA {lll}
 \bar x_{t}=\frac{2\bar x}{p-1}-\bar y\\[2mm]
 \bar y_{t}=pX_{\overline M}^{p-1}\bar x+\frac{2p}{p+1}\overline MY_{\overline M}^{\frac{p-1}{p+1}}\bar y-K\bar y+h(\bar x,\bar y),
 \EA\ee
 where
 $$h(\bar x,\bar y)=c_1\bar x^2+c_2\bar y^2+c_3\bar x^3+c_4\bar y^3+...,
 $$
 with
 $$\BA {lll}\displaystyle c_1=\frac{p(p-1)}{2}X_{\overline M}^{p-2}\,,\; c_2=\frac{p(p-1)}{(p+1)^2}\overline MY_{\overline M}^{-\frac{2}{p+1}}\\[4mm] \displaystyle  c_3=\frac{p(p-1)(p-2)}{6}X_{\overline M}^{p-3}\,,\;
c_4=-\frac{2p(p-1)}{3(p+1)^3}\overline MY_{\overline M}^{-\frac{p+3}{p+1}}.
\EA $$
 Setting $\ga^2=\frac{4}{(p-1)^2}+N-2$, we have $pX_{\overline M}^{p-1}=\ga^2$,  since $X_{\overline M}=\frac{p-1}{2}Y_{\overline M}$ and $\overline M$ satisfies  $(\ref{I-3-13'})$, thus
\bel{I-3-14b}c_1=\myfrac{\ga^2}{Y_{\overline M}}\;,\;\; c_2=\myfrac{(p-1)L}{2(p+1)Y_{\overline M}}
\;,\;\; c_3=-\myfrac{(p-1)L}{3(p+1)^2Y^2_{\overline M}}.
\ee
In order to compute the Lyapunov coefficients we transform the system by setting
 $$\gg=\sqrt{N-2}\,,\; s=\gg t\,,\;\gg w=\frac{2\bar x}{p-1}-\bar y.  $$
 The new system is
  $$\BA {lll}
 \bar x_s=w\\[2mm]
w_s=-\bar x-\frac{1}{\gg^2}h\left(\bar x,\frac{2\bar x}{p-1}-\gg w\right).
 \EA$$
By integrating the first line and using the expansion of $h$, we obtain
$$w_s=-\bar x-\frac{1}{\gg^2}\left[c_1\bar x^2+c_2\left(\frac{2\bar x}{p-1}-\gg w\right)^2+c_3\bar x^3+c_4\left(\frac{2\bar x}{p-1}-\gg w\right)^3+...\right.
$$
This can be written in the following way
$$w_s=-\bar x+\gn_{2,0}\bar x^2+\gn_{1,1}\bar x w+\gn_{0,2}w^2+\gn_{3,0}\bar x^3+\gn_{2,1}\bar x^2 w+\gn_{1,2}\bar x w^2+\gn_{0,3}w^3+...
$$
By \cite[Th. 9.2.3]{HuWe}, the Lyapunov coefficient is given by $\Gl=\gn_{2,0}+3\gn_{0,3}+\gn_{1,1}\left(\gn_{2,0}+\gn_{0,2}\right)$ which yields by computation
$$\gg\Gl=\ga^2\left(3c_4-\frac{4c_2}{p-1}(c_2\ga+c_1)\right)=-\ga^2\frac{(p-1)(N+1)}{(p+1)^2Y_{\overline M}^2}L.
$$
If $L>0$ (resp. $L<0$) $P_{\overline M}$ is a weak sink (resp. a weak source).\qeda
\subsubsection{Energy and Lyapunov functionals for system $(\ref{I-3x})$}

If $x(t)$ is a solution of $(\ref{I-3-2})$ we set
\bel{I-3-15}
F(t):=\frac{x_{t}^2}{2}+\frac{|x|^{p+1}}{p+1}-\frac{Kx^2}{p-1},
\ee
Then, if $(x(t),y(t))\in {\bf Q}$, we have
$$\BA{lll}\displaystyle F_t(t)=-Lx_{t}^2-M\left(\frac{2x}{p-1}-x_{t}\right)^{\frac{2p}{p+1}}x_{t}=-\left(Lx_{t}+M y^q\right)x_{t}\\[4mm]
\displaystyle\phantom{F'(t)}=-\left(L\left(\frac{2x}{p-1}-y\right)+M y^q\right)\left(\frac{2x}{p-1}-y\right).
\EA$$
Hence, if $LM>0$,  $F$ is monotone in the region $\left\{(x,y)\in {\bf Q}:\frac{2}{p-1}x-y>0\right\}$, located under $\CL$.\medskip

\nind\Remark a) This function was introduced classically in the case $M=0$, leading in particular to $(\ref{sob})$ when $p=\frac{N+2}{N-2}$.\\
b) Using this function we can deduce an upper estimate for regular solutions, completing, \rprop{zero} and \rprop{bound}, namely if $M>0$ and $p\geq \frac{N+2}{N-2}$ any ground state satisfies
\bel{I-3-15a}
u(r)\leq \left(\myfrac{p+1}{p-1}K\right)^{\frac{1}{p-1}}r^{-\frac{2}{p-1}}\quad\text{for all }r>0.
\ee

Next we construct a Lyapunov functional adapting the method initiated by \cite{AnLe} and already used in \cite{Bi} and \cite{Bi2}.
\blemma{Lyap} We define on $\BBR^2$
\bel{I-3-16}
\CJ(x,y)=\frac{Kx^2}{p-1}-\frac{|x|^{p+1}}{p+1}-M\left(\frac{2}{p-1}\right)^{\frac{2p}{p+1}}\frac{(p+1)|x|^{\frac{3p+1}{p+1}}}{3p+1}-\frac 12\left(\frac{2x}{p-1}-y\right)^2,
\ee
and, if $(x(t),y(t))$ is a solution of $(\ref{I-3-17})$, $\CV(t)=\CJ(x(t),y(t))$.
If $(x(t),y(t))=(x,y)\in {\bf Q}$ we have
\bel{I-3-17}\BA{lll}\displaystyle
\CV_t(t)=\left(\frac{2x}{p-1}-y\right)\left[L\left(\frac{2x}{p-1}-y\right)-M\left(\left(\frac{2x}{p-1}\right)^{\frac{2p}{p+1}}-y^{\frac{2p}{p+1}}\right)\right]\\[4mm]
\phantom{\CV'(t)}\displaystyle
=L\left(\frac{2x}{p-1}-y\right)^2-M\left(\left(\frac{2x}{p-1}\right)^{\frac{2p}{p+1}}-y^{\frac{2p}{p+1}}\right)\left(\frac{2x}{p-1}-y\right).
\EA\ee
Consequently, the function $t\mapsto \CV(t)$ is decreasing if $M>0$ and $1<p\leq \frac{N+2}{N-2}$,  and increasing if $M<0$ and $p\geq \frac{N+2}{N-2}$.
\es
\Proof We recall the ansatz  introduced in \cite{AnLe} for finding a Lyapunov function for a system of the form
\bel{I-3-18}\BA{lll}\displaystyle
x_{t}=f(x,y)\\
y_{t}=g(x,y).
\EA\ee
 If $f(x,y)=0\Longleftrightarrow y=h(x)$, then consider the function
\bel{I-3-19}\BA{lll}\displaystyle
{\bf L}(x,y)=\myint{h(x)}{y}f(x,t)dt-\myint{0}{x}g(t,h(t))dt.
\EA\ee
In the case of system $(\ref{I-3x})$, $h(x)=\frac{2x}{p-1}$, and we find ${\bf L}(x,y)=\CJ(x,y)$. Then $(\ref{I-3-17})$ and the conclusion follow.\qeda
\medskip

\nind\Remark If $LM>0$ we set
$$\CH=\left\{(x,y)\in{\bf Q}:\frac{\left(\frac{2x}{p-1}\right)^\frac{2p}{p+1}-y^\frac{2p}{p+1}}{\frac{2x}{p-1}-y}=\myfrac LM\,,\;\text{or }\;\frac{2x}{p-1}=y=\left(\myfrac{L}{Mq}\right)^{\frac{1}{q-1}}\right\}.
$$
Then 
$$\CV_t(t)=0\Longleftrightarrow (x(t),y(t))\in\CH\cup\left\{(x,y)\in{\bf Q}:\frac{2x}{p-1}=y\right\}.$$
Moreover $\CH$ is starshaped with respect to $0$ and we set
$$\CR:=\displaystyle\bigcup_{0\leq\gth\leq 1}\gth\CH=\left\{(x,y)\in{\bf Q}:\frac{\left(\frac{2x}{p-1}\right)^\frac{2p}{p+1}-y^\frac{2p}{p+1}}{\frac{2x}{p-1}-y}\leq \myfrac LM\,,\;\text{or }\;\frac{2x}{p-1}=y\leq \left(\myfrac{L}{Mq}\right)^{\frac{1}{q-1}}\right\}.
$$
If $M>0$ we have $\CV_t(t)\geq 0$ if $(x(t),y(t))\in\CR$ and $\CV_t(t)\leq 0$ if $(x(t),y(t))\in\CR^c\cap {\bf Q}$. If $M<0$, the signs of  $\CV_t(t)$ in the same regions are opposite. \medskip

In \cite{BiGaVe2} we also used a function introduced in \cite{SeZo} for equation $(\ref{I-0-1})$. When $q=\frac{2p}{p+1}$ it reduces to
\bel{I-21.1}\BA{lll}\displaystyle
\CZ(r)=r^a\left(\myfrac{p+1}{2}u_r^2+u^{p+1}+a\myfrac{uu_r}{r}+Mu\abs{u_r}^{\frac{2p}{p+1}}\right),
\EA\ee
where $a=\frac{2(p+1)(N-1)}{p+3}$. Since $r=e^t$, we find that in {\bf Q},
\bel{I-21.2}\BA{lll}\displaystyle
\CZ(r)=e^{\frac{2(p+1)L}{p+3}t}\left(\myfrac{p+1}{2}y^2+x^{p+1}-a xy+Mx{y}^{\frac{2p}{p+1}}\right).
\EA\ee
The function satisfies the relation
\bel{I-21.3}\BA{lll}\displaystyle
\CZ_r-\frac{2p}{p+1}M{u_r}^{\frac{p-1}{p+1}}\CZ=r^{a-1}\CU,
\EA\ee
where
\bel{I-21.4}\BA{lll}\displaystyle
\CU=\myfrac{2(N-1)(p^2-1)}{(p+3)^2}\myfrac{uu_r}{r}\left(-L+\myfrac{p(p+3)}{(p+1)^2}Mr{u_r}^{\frac{p-1}{p+1}}\right)\\[4mm]\phantom{\CU}
=\myfrac{2(N-1)(p^2-1)}{(p+3)^2}x y\left(-L+\myfrac{p(p+3)}{(p+1)^2}M{y}^{\frac{p-1}{p+1}}\right)r^{-\frac{p+3}{p-1}}
\EA\ee
Note that $\CU$ has a constant sign in $\bf Q$ if $LM<0$.
\subsubsection{Comparison results}
\blemma{nointer} Let $N\geq 1$, $p>1$ and $M,M'\in\BBR$ such that $M<M'$. Then, as long as they lie below the line $\CL$, i.e. $x_t>0$, the regular trajectories {\bf T}$_{reg}^{_{^{M}}}$ and {\bf T}$_{reg}^{_{^{M'}}}$ associated to $M$ and $M'$ respectively do not intersect. Furthermore {\bf T}$_{reg}^{_{^{M}}}$ is below {\bf T}$_{reg}^{_{^{M'}}}$.
\es
\Proof We use the expansion $(\ref{I-3-11w})$ and deduce, for $M_i=M$ or $M'$, that
\bel{ni-1}
\myfrac{y^{_{^{M_i}}}(t)}{x^{_{^{M_i}}}(t)}=e^{2t}\left(\myfrac{1}{N}+M_i\myfrac{e^{\frac{p-1}{p+1}t}}{N(N+\frac{2p}{p+1})}(1+o(1))\right)
\ee
as $t\to-\infty$. Hence {\bf T}$_{reg}^{_{^{M}}}$ is below {\bf T}$_{reg}^{_{^{M'}}}$ for $t\leq t^*$ for some $t^*\in\BBR$. Suppose that the trajectories intersect for a first time at some point $(x_0,y_0)$ below $\CL$. Since the system is autonomous,  there will exist two solutions of the systems relative to $M$ and $M'$ satisfying at the same time $t_0$, $x^{_{^{M}}}(t_0)=x^{_{^{M'}}}(t_0)=x_0$ and $y^{_{^{M}}}(t_0)=y^{_{^{M'}}}(t_0)=y_0$. From $(\ref{I-3-4})$, $x_t^{_{^{M}}}\!(t_0)=x_t^{_{^{M'}}}\!(t_0)>0$ and
$y_t^{_{^{M}}}\!(t_0)=y_t^{_{^{M'}}}\!(t_0)+(M-M')y_0^{\frac{2p}{p+1}}<y_t^{_{^{M'}}}(t_0)$. Hence the intersection of the two trajectories is transverse and the slope of {\bf T}$_{reg}^{_{^{M}}}$ at this intersection point is smaller that the one of {\bf T}$_{reg}^{_{^{M'}}}$ which is impossible. The second assertion follows immediately.\qeda
\mysection{Study of ground states of $(\ref{I-3x})$ when $M>0$}
When $M\geq 0$ and either $N\geq 3$ and $1<p\leq\frac{N}{N-2}$, or $N=1,2 $ and $p>1$, there exists no ground state by \rprop{zero}; in this section we assume $N\geq 3$ and $p>\frac{N}{N-2}$.
\subsection{Behaviour near the equilibrium}
Since $M\geq 0$ and $p>\frac{N}{N-2}$ there exists a unique equilibrium $P_M\in{\bf Q}$. The sign of the real part of the characteristic roots of the
linearization of the system $(\ref{I-3-4})$ at $P_M$ depends on the position of $M$ with respect to $\overline M$ defined in $(\ref{I-3-14})$.
\blemma{pointsM>0} Assume $M>0$ and $N\geq 3$. \smallskip

\nind 1) If  $\frac{N}{N-2}<p\leq\frac{N+2}{N-2}$, then $P_M$ is a source.\smallskip

\nind 2) If  $p>\frac{N+2}{N-2}$; then $P_M$ is a sink when $M<\overline M$  and it is a source when $M>\overline M$.

\nind 3) If $p>\frac{N+2}{N-2}$, then $P_{\overline M}$ a weak sink. Furthermore, if
$0<M-\overline M<\ge$, for $\ge$ small enough, there exists a periodic trajectory in ${\bf Q}$ surrounding $P_M$.
\es
\Proof {\it Step 1:} Assume $\frac{N}{N-2}<p\leq\frac{N+2}{N-2}$. The linearized system is given in $(\ref{I-3-12})$. Because $M>0$, the product of the characteristic roots given by equation $(\ref{I-3-13})$ is positive since it is given by
$$\BA {lll} \displaystyle
2K-\myfrac{2p}{p+1}MY_M^{\frac{p-1}{p+1}}=\myfrac{2}{p+1}K+\myfrac{2p}{p+1}\left(\myfrac{p-1}{2}\right)^pY_M^{p-1}.
\EA$$
The sum (or the real part) of the characteristic roots is equal to $\myfrac{2p}{p+1}MY_M^{\frac{p-1}{p+1}}-L$ which is positive, as $p\leq\frac{N+2}{N-2}$ implies $L<0$. Hence $P_M$ is a source. \smallskip

\nind {\it Step 2:} Assume $p>\frac{N+2}{N-2}$, hence $L>0$. As in Step 1, the product of the characteristic roots is positive. By \rprop{fixed-p}-(1), $MY_M^{\frac{p-1}{p+1}}$ is an increasing function of $M$, then the sum (or the real part of) of the characteristic roots, given by $\myfrac{2p}{p+1}MY_M^{\frac{p-1}{p+1}}-L$ is also increasing and vanishes if $M=\overline M$. It is negative if $0<M<\overline M$ and positive if $M>\overline M$. It implies assertion 2.\smallskip

\nind {\it Step 3:} If $M=\overline M$, and $p>\frac{N+2}{N-2}$, then $P_{\overline M}$ is a weak sink by
 \rlemma{barM}. The appearance of the limit cycle, which is the called the Andronov-Hopf bifurcation, occurs for $M>\overline M$ when $M-\overline M$ is small enough (see \cite[Chap. 9]{HuWe}). This implies assertion 2.\qeda\medskip

 \nind\Remark The product of the characteristic roots is also expressed by
$$\BA {lll} \displaystyle
2K-\myfrac{2p}{p+1}MY_M^{\frac{p-1}{p+1}}=\myfrac{2M}{p+1}Y_M^{\frac{p-1}{p+1}}+2\left(\myfrac{p-1}{2}\right)^2Y_M^{p-1}.
\EA$$
 Hence it is positive for any $M\geq 0$, $p>1$. \medskip

 Next we give some sufficient conditions for nonexistence of a periodic solution or a homoclinic orbit at $(0,0)$.

 \blemma {cycleM>0} Assume $N\geq 3$. If $M>0$ and $1<p\leq \frac{N+2}{N-2}$ the system $(\ref{I-3x})$ admits no closed orbit in $\overline{\bf Q}$.  If $0<M\leq \overline M$ and
 $p>\frac{N+2}{N-2}$, the system $(\ref{I-3x})$ admits no cycle in ${\bf Q}$ surrounding $P_M$.
 \es
 \Proof If $\gg$ is a non-trivial closed orbit  it corresponds either to a $T$-periodic solution or a solution such that
 $$\lim_{t\to -\infty}(x(t),y(t))=\lim_{t\to \infty}(x(t),y(t)).
 $$
 The function $\CV$ defined in $(\ref{I-3-16})$ is monotone and it satisfies either $\CV(0)=\CV(T)$ or
 $\displaystyle\lim_{t\to -\infty}\CV(t)=\lim_{t\to \infty}\CV(t).$ Hence it is constant and by $(\ref{I-3-17})$ it implies $y(t)=\frac{2}{p-1}x(t)$ for all $t$, a contradiction.

 \smallskip

 Next we suppose $p> \frac{N+2}{N-2}$, $0<M\leq \overline M$ and that there exists a $T$-periodic solution $(x(t),y(t)$ with trajectory $\gg\subset{\bf Q}$ surrounding $P_M$, hence $P_M$ belongs to the bounded connected component
 $\Gg$ of $\BBR^2\setminus\{\gg\}$ bordered by $\gg$. Since $P_M$ is a sink or a weak sink by \rlemma{pointsM>0}, there exists a neighbourhood $\CO$ of $P_M$ such that all the trajectories issued from $\CO$ converge to  $P_M$ as $t\to\infty$. Hence any trajectory issued from $\CO$, necessarily contained in $\Gg$, has an $\ga$-limit set in $\overline \Gg$ which is either a stationary point different from $P_M$, which is excluded, or a limit cycle $\{(x(t),y(t)\}_{t\in [0,\gt)}:=\gg'\subset \overline \Gg$ ($\gt$ is its period). This limit cycle is not stable, hence, by Floquet's theory
 \bel{Flo}
 \myint{0}{\gt}\left(\myfrac{\partial H_1(x(t),y(t))}{\partial x}+ \myfrac{\partial H_2(x(t),y(t))}{\partial y}\right)dt= \myint{0}{\gt}\left(\myfrac{2pM}{p+1}y^{\frac{p-1}{p+1}}(t)-L\right)dt\geq 0
 \ee

 We perform the change of unknowns $\bar x=x-X_M$, $\bar y=y-Y_M$ used in \rlemma{barM} which leads to the system $(\ref{I-3-14a})$. The explicit value of the remaining term is
 \bel{Flo1}\BA {lll}h(\bar x,\bar y)=\left(\bar x+X_M\right)^p-X_M^p-pX_M^{p-1}\bar x+
 M\left(\bar y+Y_M\right)^{\frac{2p}{p+1}}-MY_M^{\frac{2p}{p+1}}-\frac{2p}{p+1}MY_M^{\frac{p-1}{p+1}}\bar y\\
 \phantom{h(\bar x,\bar y)}=\Phi(\bar x)+M\Psi(\bar y),
 \EA \ee
 where $\Phi$ and $\Psi$ are defined accordingly. It is positive by convexity because $M>0$. Since from $(\ref{I-3-14a})$,
 $$\myfrac{2}{p-1}\myint{0}{\gt}\bar x(t) dt-\myint{0}{\gt}\bar y(t) dt=0,
 $$
and
 $$\BA {lll}
0= pX_M^{p-1}\myint{0}{\gt}\bar x(t) dt+\myfrac{2p}{p+1}MY_M^{\frac{p-1}{p+1}}\myint{0}{\gt}\bar y(t) dt
 -K\myint{0}{\gt}\bar y(t) dt+\myint{0}{\gt}h(\bar x,\bar y) dt
 \EA$$
 we derive
$$
 \left(\myfrac{p(p-1)^p}{2^p}Y_M^{p-1}+\myfrac{2p}{p+1}MY_M^{\frac{p-1}{p+1}}-K\right)\myint{0}{\gt}\bar y(t) dt<0.
$$
Using the equation $(\ref{I-3-9})$ satisfied by $Y_M$ it yields
$$(p-1)\left(K-\myfrac{M}{p+1}Y_M^{\frac{p-1}{p+1}}\right)\myint{0}{\gt}\bar y(t) dt<0.
$$
By $(\ref{I-3-8a})$-(ii), $K-\myfrac{M}{p+1}Y_M^{\frac{p-1}{p+1}}>0$, hence $\myint{0}{\gt}\bar y(t) dt<0$, therefore
$\gt^{-1}\myint{0}{\gt}y(t) dt< Y_M$ and by concavity,
 \bel{Flo2}
 \myfrac 1\gt\myint{0}{\gt}(y(t))^{\frac{p-1}{p+1}}dt\leq \left( \myfrac 1\gt\myint{0}{\gt}y(t)dt\right)^{\frac{p-1}{p+1}}
<Y_M^{\frac{p-1}{p+1}}.
\ee
Combining $(\ref{Flo})$ and $(\ref{Flo2})$, we obtain
 \bel{Flo3}
 L<\myfrac {2pM}{p+1}Y_M^{\frac{p-1}{p+1}}.
\ee
Since $M\mapsto MY_M^{\frac{p-1}{p+1}}$ is decreasing by \rprop{fixed-p}-(1-ii) we have for $0\leq M\leq\overline M$
$$0=\myfrac {2p\overline M}{p+1}Y_{\overline M}^{\frac{p-1}{p+1}}-L\geq \myfrac {2pM}{p+1}Y_M^{\frac{p-1}{p+1}}-L,
$$
which contradicts $(\ref{Flo3})$.Ê\qeda\medskip

\nind\Remark Up to changing the sense of variation of $\CV(t)$, the proof of the first assertion shows that there exists no closed orbit
in $\overline{\bf Q}$ if $M<0$ and $p\geq \frac{N+2}{N-2}$. However the proof of the second assertion is not valid when $M<0$.\medskip

The nonexistence of any periodic solution can also be proved when the equilibrium is a node (i.e. the two characteristic values are real with the same sign).

\blemma{node} 1- Let $N\geq 3$ and $p>\frac{N+2}{N-2}$. There exists a unique and explicit $M_0> \overline M$ such that for any $M\geq M_0$, $P_M$ is a repelling node, degenerate if $M= M_0$. If $\overline M<M<M_0$, $P_M$ is a repelling spiraling point. If $3\leq N\leq 10$ and
$0<M<\overline M$, $P_M$ is an attracting spiraling point. \\
2- If $ N\geq  11$ and $\frac{N+2}{N-2}<p<\frac{N-2\sqrt{N-1}}{N-4-2\sqrt{N-1}}$ there exists a unique and explicit $0<M_1< \overline M$ such that if $P_M$ is an attracting node if $0<M<M_1$ and an attracting spiraling point if $M_1<M<\overline M$.\\
3- If $M>M_0$ or if $0<M<M_1$, there exist no periodic trajectory in ${\bf Q}$ around $P_M$, neither no homoclinic trajectory in ${\bf Q}$ at $(0,0)$ surrounding $P_M$.

\es
\Proof The characteristic values of the trinomial $T(\gl)$ defined in $(\ref{I-3-13})$ are real if and only if its discriminant $D$ is nonnegative. By computation we find
 \bel{node1}\BA {lll}
D=\left(\myfrac{2p}{p+1}MY_M^{\frac{p-1}{p+1}}-L\right)^2-8\left(K-\myfrac{p}{p+1}MY_M^{\frac{p-1}{p+1}}\right)\\[4mm]
\phantom{D}
=\left(\myfrac{2p}{p+1}MY_M^{\frac{p-1}{p+1}}\right)^2+(2-L)\myfrac{4p}{p+1}MY_M^{\frac{p-1}{p+1}}+L^2-8K\\[4mm]
\phantom{D}
=\left(\myfrac{2p}{p+1}MY_M^{\frac{p-1}{p+1}}-L+2\right)^2-4(1-L+2K).
\EA\ee
Observing that $1-L+2K=N-1$ by $(\ref{I-3v})$, we deduce
$$D=\left(\myfrac{2p}{p+1}MY_M^{\frac{p-1}{p+1}}-L+2+2\sqrt{N-1}\right)\left(\myfrac{2p}{p+1}MY_M^{\frac{p-1}{p+1}}-L+2-2\sqrt{N-1}\right).
$$
Hence  $(\ref{node1})$ is satisfied if one of the following conditions holds:
\bel{node1a}\BA {lll}
(i)\qquad\qquad \qquad  &MY_M^{\frac{p-1}{p+1}}\leq \myfrac{p+1}{2p}\left(L-2-2\sqrt{N-1}\right),\qquad \qquad\qquad\qquad\quad\quad \\
(ii)\qquad\qquad \qquad  &MY_M^{\frac{p-1}{p+1}}\geq \myfrac{p+1}{2p}\left(L-2+2\sqrt{N-1}\right).\qquad \qquad\qquad\qquad\quad\quad
\EA
\ee
It is easy to check that for $p>\frac{N+2}{N-2}$, one has $0<L-2+2\sqrt{N-1}<qK$.
Since $M>0$, by \rprop{fixed-p}-(1), the mapping $M\mapsto MY_M^{\frac{p-1}{p+1}}$ is continuous, increasing and from $[0,\infty)$ into $[0,\infty)$. Therefore there exists a unique $M_0$ such that
 \bel{node2}\BA {lll}
M_0Y_{M_0}^{\frac{p-1}{p+1}}= \myfrac{p+1}{2p}\left(L-2+2\sqrt{N-1}\right).
\EA\ee
Using $(\ref{I-3-9})$ we find
 \bel{node3}\BA {lll}
M_0=\myfrac{p+1}{2p}\left(\myfrac{p-1}{2}\right)^{\frac{p}{p+1}}
\myfrac{K-\frac{2p}{p-1}+2\sqrt{N-1}}{\left(\frac{p-1}{2p}K+\frac{p+1}{p-1}-\frac{(p+1)\sqrt{N-1}}{p}\right)^{\frac{1}{p+1}}}.
\EA\ee
Concerning the upper bound in (i), there holds
 \bel{node4}L-2-2\sqrt{N-1}>0\Longleftrightarrow p>\myfrac{N-2\sqrt{N-1}}{N-2\sqrt{N-1}-4}\quad\text{and }\;N\geq 11.
\ee
Then we can define $M_1$ by
\bel{node5}\BA {lll}
M_1Y_{M_1}^{\frac{p-1}{p+1}}= \myfrac{p+1}{2p}\left(L-2-2\sqrt{N-1}\right).
\EA\ee
which leads to
 \bel{node6}\BA {lll}
M_1=\myfrac{p+1}{2p}\left(\myfrac{p-1}{2}\right)^{\frac{p}{p-1}}\myfrac{K-\frac{2p}{p-1}-2\sqrt{N-1}}{\left(\frac{p-1}{2p}K+\frac{p+1}{p-1}+\frac{(p+1)\sqrt{N-1}}{p}\right)^{\frac{1}{p+1}}}.
\EA\ee
Note that
$$\BA {lll}\frac{2p}{p+1}M_1Y_{M_1}^{\frac{p-1}{p+1}}=L-2-2\sqrt{N-1}<L=\frac{2p}{p+1}\overline MY_{\overline M}^{\frac{p-1}{p+1}}\leq L-2+2\sqrt{N-1}= \frac{2p}{p+1}M_0Y_{M_0}^{\frac{p-1}{p+1}}.
\EA$$
Hence
\bel{M01}
M_1<\overline M<M_0.
\ee
Next we prove 3) by adapting an argument introduced in \cite{ChTi} for quadratic systems. We return to system $(\ref{I-3-14a})$ that we write under the form
 \bel{node7}\BA {lll}
\bar x_{t}=a\bar x-\bar y\\
\bar y_i=c\bar x+d\bar y+\Phi(\bar x)+M\Psi(\bar y)
\EA\ee
which defines $a$, $c$, $d$ with $\Phi$ and $\Psi$ given by $(\ref{node7})$, and the trinomial $(\ref{I-3-13})$ for characteristic values endows the form
$$T(\gl)=\gl^2-(d+a)\gl+ad+c.
$$
In the range of values of $M$, the discriminant $D=(d+a)^2-4(ad+c)=(d-a)^2-4c$ is positive. We consider the intersection of a straight line $\ell$  passing through $P_M$ with equation $\bar y=A\bar x$ with a trajectory $(\bar x(t),\bar y(t))$. Then
$$\BA {lll}U=\bar y_{t}-A\bar x_{t}=(A^2+(d-a)A+c)\bar x+\Phi(\bar x)+M\Psi(\bar y)\\[2mm]
\phantom{U=\bar y_{t}-A\bar x_{t}}=(A^2+(d-a)A+c)\bar x+\Phi(\bar x)+M\Psi(A\bar x).
\EA$$
We can choose $A\neq 0$ such that $A^2+(d-a)A+c=0$ since $D>0$. Since $\Phi$ and $\Psi$ achieve positive values, we derive from the expression of $h$ that  $U>0$ for $\bar x\neq 0$. This proves that any closed orbit around $P_M$ or passing by $P_M$ can intersect $\ell$ only one time which is a contradiction. $\phantom{----------}$\qeda

\subsection{Existence or nonexistence of ground states}

\bprop{subcritM>0} Let $N\geq 3$, $\frac{N}{N-2}< p\leq \frac{N+2}{N-2}$ and $M>0$. Then there exists no radial ground state, and there exists a singular solutions $u(r)$ such that
$\displaystyle\lim_{r\to 0}r^{\frac{2}{p-1}}u(r)= X_M$ and  $\displaystyle\lim_{r\to \infty}r^{N-2}u(r)= c$ for some $c>0$.
\es
\Proof Assume that ${\bf T}_{reg}$ remains in {\bf Q}, by \rlemma{tetra} we have three possibilities:\\
($\ga$) either ${\bf T}_{reg}$ converges to $P_M$ when $t\to\infty$, which is impossible since $P_M$ is a source,\\
($\gb$) or ${\bf T}_{reg}$ has a limit cycle at $\infty$, and this is impossible by  \rlemma{cycleM>0},\\
($\gg$) or ${\bf T}_{reg}$ converges to $(0,0)$ when $t\to\infty$, hence it is a homoclinic orbit. The function
$\CV$ defined in $(\ref{I-3-16})$ is decreasing by \rlemma{Lyap}. Since $\CV(-\infty)=\CV(\infty)=0$, it is identically $0$ and so is $\CV_t$. This implies that $\frac{2x(t)}{p-1}-y(t)=0$ for all $t$ which is a contradiction. \\
Hence ${\bf T}_{reg}$ does not remain in {\bf Q}. \smallskip

We denote by $\Gx$ the connected region of {\bf Q} bordered by
the semi-axis $\{(x,y):x=0,y>0\}$ and  ${\bf T}_{reg}$. Since {\bf H} is outward on the semi-axis, $\Gx$ is negatively invariant. By Section 3.4.1, $(0,0)$ is a saddle point, hence the stable trajectory ${\bf T}_{st}=\{(x(t),y(t)\}$ satisfies
$\displaystyle\lim_{t\to\infty}\frac{y(t)}{x(t)}=N-2$ which implies $u(r)\sim cr^{2-N}$ for some $c>0$. Its $\ga$-limit set $\ga({\bf T}_{st})$ cannot be a limit cycle as we have seen it above. If it contains $(0,0)$ it implies again that
$\CV$, which is monotone, is equal to $0$, hence $\CV_t\equiv 0$ and  $\frac{2x(t)}{p-1}-y(t)\equiv 0$, which is impossible. Hence $\ga({\bf T}_{st})$ contains $P_M$. Since $P_M$ is a source it implies that ${\bf T}_{st}$ converges to $P_M$ when $t\to-\infty$. \qeda

\bprop{0<M<Mbar} Let $N\geq 3$ and $p>\frac{N+2}{N-2}$. Then for any $0<M\leq \overline M$ there exists a ground state $u$ which satisfies
$u\sim U_M$ at $\infty$.
\es
\Proof If $0<M<\overline M$ (resp. $M=\overline M$), $P_M$ is a sink (resp. a weak sink). Suppose first that the trajectory ${\bf T}_{reg}$ does not stay in {\bf Q}, then it leaves {\bf Q} at some point $(0,y_s)$ with $y_s>0$. As a consequence, the stable trajectory ${\bf T}_{st}$ at $(0,0)$ remains in the negatively invariant region $\Gx$ defined in the proof of
\rprop{subcritM>0}. Since it cannot converge to $P_M$ when $t\to-\infty$ it admits a limit cycle surrounding $P_M$ which contradicts \rlemma {cycleM>0}. Therefore ${\bf T}_{reg}\subset$ {\bf Q} and, again using \rlemma {cycleM>0}, either it converges to $P_M$ when $t\to\infty$ and the proof is complete, or to $(0,0)$ and ${\bf T}_{reg}={\bf T}_{st}$ is a homoclinic trajectory. The trace of the linearized system $(\ref{I-3-11})$ at $(0,0)$ is equal to $\frac{2}{p-1}-K=-L<0$. Therefore, from \cite[Th. 9.3.3]{HuWe} the connection is attracting and the trajectories inside the bounded region $\CT$ bordered by the homoclinic trajectory ${\bf T}_{reg}$ spiral towards it when $t\to\infty$. Hence any such trajectory inside $\CT$ either has a limit cycle when $t\to-\infty$ which is impossible by \rlemma {cycleM>0} or converges to $P_M$ which is also impossible. Consequently there exists no homoclinic trajectory at $(0,0)$ which ends the proof.
\qeda\medskip

Next we study the case where $M$ is large enough. We have already proved in \cite{BiGaVe2} that for any $p>1$,  there exists $M_\dag=M_\dag(N,p)>0$ (see introduction) such that if $M>M_\dag$ there exists no ground state, radial or non-radial. In the radial case we have a more precise result.
\bprop{largeM} Let $N\geq 3$ and $p>\frac{N+2}{N-2}$. Then for any $M\geq M_0$ there exists no ground state, but there exist
singular solutions $u$ which satisfy $\displaystyle\lim_{r\to 0}r^{\frac{2}{p-1}}u(r)=X_M$ and $\displaystyle\lim_{r\to \infty}r^{N-2}u(r)= c>0$.
\es
\Proof Since $M\geq M_0>\overline M$, $P_M$ is a source by \rlemma {pointsM>0}, and there exist no periodic orbit neither a homoclinic trajectory at $(0,0)$ by \rlemma {node}. Thus ${\bf T}_{reg}$ leaves ${\bf Q}$ through the semi-axis $\{(0,y):y>0\}$ by \rlemma {tetra}. As in the proof of \rprop{0<M<Mbar} the stable trajectory ${\bf T}_{st}$ at $(0,0)$ remains in the negatively invariant region $\Gx$ already defined. Then it converges necessarily to $P_M$ when $t\to-\infty$. \qeda\medskip

Next we study the case $\overline M<M<M_0$.
\bth{Mbar<M} Let $N\geq 3$, $p>\frac{N+2}{N-2}$. There exist two positive real numbers $\tilde M_{min}$ and $\tilde M_{max}$ such that
$\overline M<\tilde M_{min}\leq \tilde M_{max}<M_0$ such that,\smallskip

\nind 1- For $\overline M<M<\tilde M_{min}$ there exist ground states ondulating around $U_M$ when $r\to\infty$ and positive singular solutions ondulating around $U_M$ on $[0,\infty)$. \smallskip

\nind 2- For $M=\tilde M_{min}$ and for $M=\tilde M_{max}$ there exist ground states $u$ such that $\displaystyle\lim_{r\to \infty}r^{N-2}u(r)= c>0$.\smallskip

\nind 3- For $\tilde M_{max}<M<M_0$ there exists no ground state and there exist singular solutions such that $\displaystyle\lim_{r\to 0}r^{\frac{2}{p-1}}u(r)=X_M$ or turning around $P_M$ when $r\to 0$, and $\displaystyle\lim_{r\to \infty}r^{N-2}u(r)= c>0$.
\es
\Proof Since $M$ is subject to vary, we put it in exponent in the different specific trajectories of the system. For $\overline M<M<M_0$, $P_M$ is a source and the trajectories converging toward this point when $t\to-\infty$ are spiralling. We have three possibilities:\\
(i) either ${\bf T}^{^{M}}_{reg}$ leaves {\bf Q} at some point $(0,y_s)$, $y_s>0$, \\
(ii) or ${\bf T}^{^{M}}_{reg}$ has a $\gw$-limit set which is a periodic orbit surrounding $P_M$,\\
(iii) or ${\bf T}^{^{M}}_{reg}$ converges to $(0,0)$ when $t\to\infty$. \smallskip

\nind If (i) holds, then ${\bf T}^{^{M}}_{st}$ remains in the region $\Gx:=\Gx^M$ of ${\bf Q}$ bordered by ${\bf T}^{^{M}}_{reg}$ and the semi-axis $\{(0,y):y>0\}$. Then, either it converges to $P_M$ when $t\to-\infty$, or it admits an unstable (from outside) $\ga$-limit cycle. We denote by $({x}^{_{M}}_{st}, {y}^{_{M}}_{st})$ the first backward intersection of ${\bf T}^{^{M}}_{st}$ with the straight line $\CL$. Furthermore
${\bf T}^{^{M}}_{reg}$ intersects $\CL$ at some point (necessarily unique) $({x}^{_{M}}_{reg}, {y}^{_{M}}_{reg})$. Because of the relative position of ${\bf T}^{^{M}}_{st}$ and ${\bf T}^{^{M}}_{reg}$, there holds
$$g(M)={x}^{_{M}}_{reg}-{x}^{_{M}}_{st}>0.
$$

\nind If (ii) holds, then we claim that ${\bf T}^{^{M}}_{st}$ leaves ${\bf Q}$ at some point $(x_s,0)$ with $x_s>0$. Indeed, if ${\bf T}^{^{M}}_{st}\subset {\bf Q}$, it is bounded by \rprop{zero}. Hence the $\ga$-limit set is non-empty. It cannot be $(0,0)$ since ${\bf T}^{^{M}}_{st}\neq {\bf T}^{^{M}}_{reg}$ and ${\bf T}^{^{M}}_{st}$ cannot converge to $P_M$ or have a limit cycle around $P_M$ since it would imply that ${\bf T}^{^{M}}_{st}\cap {\bf T}^{^{M}}_{reg}\neq\{\emptyset\}$. Hence ${\bf T}^{^{M}}_{st}$ intersects $\CL$ at some point $({x}^{_{M}}_{st}, {y}^{_{M}}_{st})$  and
${\bf T}^{^{M}}_{reg}$ intersect the first time $\CL$ and $({x}^{_{M}}_{reg}, {y}^{_{M}}_{reg})$. Because
${\bf T}^{^{M}}_{reg}$ lies in the region of {\bf Q} bordered by ${\bf T}^{^{M}}_{st}$ and the semi-axis $\{(x,0):x>0\}$, there holds
$$g(M)={x}^{_{M}}_{reg}-{x}^{_{M}}_{st}<0.$$

\nind If (iii) holds then $g(M)=0$. \smallskip

\nind The function $g$ defined for any $M\in (\overline M,M_0)$ is continuous. We know that there exists $\ge>0$ such that for any $M\in (\overline M,\overline M+\ge)$, ${\bf T}^{^{M}}_{reg}$ has a $\omega$-limit cycle surrounding $P_M$, hence $g(M)<0$. \\
If $M=M_0$, ${\bf T}^{^{M}}_{reg}$ leaves {\bf Q} at some point $(0,y_s)$ with $y_s>0$ from \rprop{largeM} and the intersection of ${\bf T}^{^{M}}_{reg}$ with the semi-axis $\{(0,y):y>0\}$ is transverse at $(0,y_s)$. By continuity with respect to the parameter $M$, for any $M\in (M_0-\ge',M_0]$, for $\ge'>0$ small enough, ${\bf T}^{^{M}}_{reg}$ leaves {\bf Q} at some point of the semi-axis. So we are in situation (i) and $g(M)>0$. \\
Then, by continuity of $g$ there exist $\tilde M_{min}$ and $\tilde M_{max}$ such that $\overline M<\tilde M_{min}\leq \tilde M_{max}<M_0$
 such that $g(\tilde M_{min})=g(\tilde M_{max})=0$. If $M=\tilde M_{min}$ or $M=\tilde M_{max}$, ${\bf T}^{^{M}}_{st}={\bf T}^{^{M}}_{reg}$
and the trajectory ${\bf T}^{^{M}}_{reg}$ is  homoclinic at $(0,0)$. For $\overline M<M<M_{min}$ we are in situation (ii), which proves 1. For $M>\tilde M_{max}$ we are in situation (i) and
either ${\bf T}^{^{M}}_{st}$ converges to $P_M$ or has a $\alpha$-limit cycle surrounding $P_M$ when $t\to-\infty$.\qeda

\medskip

Theorems A and A' follow from the previous results.\medskip

\nind\Remark It is a challenging question to prove that there is a unique $M$ such that there is a homoclinic trajectory at $(0,0)$. Up to now all we can prove is that if there exist two parameters $0<M_1<M_2$ such that for each of them there exists a homoclinic trajectory of $(0,0)$ in {\bf Q}
${\bf T}^{^{M_i}}:={\bf T}^{^{M_i}}_{st}={\bf T}^{^{M_i}}_{unst}$ (i=1, 2), then ${\bf T}^{^{M_2}}$ is a subset of the domain of {\bf Q} limited by ${\bf T}^{^{M_1}}$.
\mysection{ Study of ground states of ($\ref{I-3x})$ when $M<0$}

We will distinguish the cases $p\geq \frac{N}{N-2}$ where system ($\ref{I-3x})$ admits a unique non-trivial equilibrium and $1<p<\frac{N}{N-2}$ where the existence of zero, one or two equilibria depends on the value of $M$ with respect to $-\gm^*$ defined in $(\ref{I-6})$. In order to avoid confusion we set
\bel{M1}
\overline \gm=\abs{\overline M}=\frac{(p+1)\abs{p(N-2)-(N+2)}}{(4p)^{\frac{p}{p+1}}\left((N-2)(p-1)^2+4\right)^\frac{1}{p+1}}.
\ee
Then $\overline\gm=\overline M$ if  $p>\frac{N+2}{N-2}$ and $\overline\gm=-\overline M$ if  $p<\frac{N+2}{N-2}$.
\bprop{1point} let $N\geq 3$, $p\geq\frac{N}{N-2}$ and $M<0$. \\
1- If $p\geq \frac{N+2}{N-2}$, $P_M$ is a sink.\\
2- If $\frac{N}{N-2}\leq p<\frac{N+2}{N-2}$, then $P_M$ is a sink if $M<\overline M$, it is a source if
$\overline M<M<0$, and $P_{\bar M}$ is a weak source. Moreover there exists $M_1<\overline M$ such that $P_M$ is a node if $M\leq M_1$. If
$\frac{N}{N-2}\leq p<\frac{N+2\sqrt{N-1}}{N-4+2\sqrt{N-1}}$, there exists $M_0\in (\overline M,0)$ such that $P_M$ is a node for
$M\geq M_0$ and a spiraling equilibrium if $M_1<M<M_0$.
\es
\Proof The equation $(\ref{I-3-13})$ satisfied by the characteristic roots can be  written under the form
\bel{M2}
T(\gl)= \gl^2+\left(\frac{2p}{p+1}\abs MY_M^{\frac{p-1}{p+1}}+L\right)\gl+2K+\frac{2p}{p+1}\abs MY_M^{\frac{p-1}{p+1}}=0.
\ee
Since $K\geq 0$ the product of the roots is positive; the real part of the roots is positive if and only if $\frac{2p}{p+1}\abs MY_M^\frac{p-1}{p+1}+L<0$. If $p\geq\frac{N+2}{N-2}$, then $L\geq 0$ and $P_M$ is a sink. If $\frac{N}{N-2}\leq p<\frac{N+2}{N-2}$, then $\overline M$ is characterized by
\bel{M3}
\frac{2p}{p+1}\overline  M Y^{\frac{p-1}{p+1}}_{\overline  M}= L<0.
\ee
By \rprop{fixed-p} the function $M\mapsto MY^{\frac{p-1}{p+1}}_M$ is increasing and onto on $(-\infty,0)$. Therefore, $P_M$ is a sink if $M<\overline M$ and a source if $\overline M<M<0$. Finally, if $M=\overline M$, the two roots are imaginary and by \rlemma{barM} $P_M$ is a weak source. Finally, from $(\ref{node1a})$, the characteristic roots are real if and only
 \bel{M5}
\myfrac{2p\abs M}{p+1}Y_M^{\frac{p-1}{p+1}}\leq \abs L+2-2\sqrt{N-1}\;\text{ or }\; \myfrac{2p\abs M}{p+1}Y_M^{\frac{p-1}{p+1}}\geq \abs L+2+2\sqrt{N-1}.
\ee
Notice that there always hold $\abs L+2-2\sqrt{N-1}<\abs L$ since $N\geq 3$.  The first condition in $(\ref{M5})$ requires  $ \abs L+2-2\sqrt{N-1}>0$, equivalently $p<\frac{N+2\sqrt{N-1}}{N-4+2\sqrt{N-1}}$. Since there holds $\frac{N}{N-2}<\frac{N+2\sqrt{N-1}}{N-4+2\sqrt{N-1}}<\frac{N+2}{N-2}$, the conclusion follows and $M_0$ and
$M_1$ are given by $(\ref{node3})$ and $(\ref{node6})$.\qeda\medskip

The most intricate case corresponds to $1<p<\frac{N}{N-2}$ or $N=1,2$ where there may exist 0, 1 or 2 equilibria.
\bprop {2points} Assume $N=1, 2$ and $p>1$ or $N\geq 3$ and $1<p<\frac{N}{N-2}$, $M<-\gm^*$ and let $P_{j,M}$, j=1 or 2, be the two equilibria  of $(\ref{I-3x})$.\\
1- Then $P_{1,M}$ is a saddle point.  \\
2- Let $N\geq 3$. If  $\overline M<M<-\gm^*$, then $P_{2,M}$ is a source; if $M<\overline M$, then $P_{2,M}$ is a sink; $P_{2,\overline M}$ is a weak source. Moreover there exist $M_1<\overline M$
such that $P_{2,M}$ is a node for $M\leq M_1$; there exists also $M_0\in (\overline M,-\gm^*)$ such that $P_{2,M}$ is a node for $M\leq M_1$ or for
$M_0\leq M<-\gm^*$, and it is a spiraling equilibrium if $M_1<M<M_0$.\\
3- Let $N=2$. Then $P_{2,M}$ is a sink. There exists $M_1<-\gm^*$ such that $P_{2,M}$ is a node if and only if $M\leq M _1$.\\
4- If $N=1$, then $P_{2,M}$ is a sink and a node.
\es
\Proof Recall that, from \rprop{fixed-p}, $M\mapsto MY_{1,M}^{\frac{p-1}{p+1}}$ is decreasing and $M\mapsto MY_{2,M}^{\frac{p-1}{p+1}}$ is increasing.
We first consider  the linearized operator at $P_{1,M}$. The product of the roots of $(\ref{M2})$ is equal to
$2K-\frac{2p}{p+1}MY_{1,M}^{\frac{p-1}{p+1}}$. Since $M<-\gm^*$,
$$
2K-\frac{2p}{p+1}MY_{1,M}^{\frac{p-1}{p+1}}<2K+\frac{2p}{p+1}\gm^* Y_{-\gm^*}^{\frac{p-1}{p+1}}=0.
$$
Hence $P_{1,M}$ is a saddle point.\smallskip

\nind Next we consider $P_{2,M}$. If $\overline M<M<-\gm^*$, since  $M\mapsto MY_{2,M}^{\frac{p-1}{p+1}}$ is increasing,
\bel{M7}MY_{2,M}^{\frac{p-1}{p+1}}>\overline M Y_{2,\overline M}^{\frac{p-1}{p+1}}=\frac{p+1}{2p} L,\ee
and
\bel{M8}
2K-\frac{2p}{p+1}MY_{2,M}^{\frac{p-1}{p+1}}>2K+\frac{2p}{p+1}\gm^* Y_{-\gm^*}^{\frac{p-1}{p+1}}=0.
\ee
Hence $P_{2,M}$ is a source. If $M<\overline M$, then the sign in $(\ref{M7})$ is reversed and
\bel{M9}
2K-\frac{2p}{p+1}MY_{2,M}^{\frac{p-1}{p+1}}>2K-\frac{2p}{p+1}\overline\gm Y_{\overline\gm}^{\frac{p-1}{p+1}}>2K+\frac{2p}{p+1}\gm^* Y_{-\gm^*}^{\frac{p-1}{p+1}}=0.
\ee
Thus $P_{2,M}$ is a sink. By \rlemma{barM} $P_{2,\overline M}$ is a weak source, and assertion 2  is proved. Finally, if $N=2$ and $M<-\gm^*=\overline M$ from \rlemma{comp}, therefore $P_{2,M}$ is a sink.\smallskip

\nind Next we look for conditions which insure that $P_{2,M}$ is a node. The characteristic roots are real if and only if one of the two conditions $(\ref{node1a})$ where $Y_M$ is replaced by $Y_{2,M}$ holds:
\bel{M10}\BA {lll}
(i)\qquad\qquad \qquad  &\myfrac{2p}{p+1}MY_{2,M}^{\frac{p-1}{p+1}}\geq L-2+2\sqrt{N-1},\qquad \qquad\qquad\qquad\quad\quad\\[4mm]
 (ii)\qquad\qquad \qquad  &\myfrac{2p}{p+1}MY_{2,M}^{\frac{p-1}{p+1}}\leq L-2-2\sqrt{N-1}.\qquad \qquad\qquad\qquad\quad\quad
\EA\ee
For $1<p<\frac N{N-2}$, $M\mapsto \Phi(M):=\frac{2p}{p+1}MY_{2,M}^{\frac{p-1}{p+1}}$ is an increasing diffeomorphism from $(-\infty,-\gm^*]$ to
$(-\infty,2K]$. Since  $L-2\sqrt{N-1}-2<2K$ there exists a unique $ M_1<-\gm^*$ such that $\Phi( M_1)=L-2-2\sqrt{N-1}$ and $(\ref{M10})$-(ii) holds when $ M\leq M_1$.Since  $L-2+2\sqrt{N-1}>2K$ is equivalent to $(N-2)^2<0$, there exists no $M\leq -\gm^*$ such that $(\ref{M10})$-(i) holds. When $N=2$ only $M=-\gm^*$ satisfies $(\ref{M10})$-(i) with equality. If $N\geq 3$ there exists a unique $ M_0<-\gm^*$ such that $(\ref{M10})$-(i) holds. Hence for any $ M_0\leq M<-\gm^*$, inequality $(\ref{M10})$-(i) holds with $M_0$ and $M_1$ defined by
$(\ref{node3})$ and $(\ref{node6})$.\smallskip

\nind At end assume $N=1$. For any $M\leq -\gm^*(1)=-(p+1)\left(\myfrac{p+1}{2p}\right)^\frac{p+1}{p}$, there holds
\bel{M11}\BA {lll}
\myfrac{2p}{p+1} MY_{2,M}^{\frac{p-1}{p+1}}- L\leq -\myfrac{2p}{p+1}\gm_*Y_{2,-\gm^*}^{\frac{p-1}{p+1}}=2 K- L=-1.
\EA\ee
Therefore $P_{2,M}$ is always a sink and the discriminant of $(\ref{I-3-13})$ is equal to
$\left(\frac{2p}{p+1} MY_{2,M}^{\frac{p-1}{p+1}}- L+2\right)^2$, hence the characteristic roots are negative.\qeda\medskip

\nind{\it Notations}. Under the assumptions of \rprop{2points} there exist two stable trajectories ${\bf T}^{^{1,j}}_{st}$, j=1,2, converging to
$P_{1,M}$ when $t\to\infty$, associated to a negative characteristic value $\gl$;  the common slope at  $P_{1,M}$ is $\frac{2}{p-1}+\abs\gl$. We assume that ${\bf T}^{^{1,1}}_{st}$ is {\it locally below} $\CL$ and ${\bf T}^{^{1,2}}_{st}$ {\it locally above} $\CL$, thus these trajectories are coming from the regions {\bf C} and {\bf A} defined in the proof of \rprop{bound}. There exist also  two unstable trajectories ${\bf T}^{^{1,j}}_{unst}$, j=3,4, converging to
$P_{1,M}$ when $t\to-\infty$, associated to a positive characteristic value $\gl'$ and the common slope at  $P_{1,M}$ is $\frac{2}{p-1}-\gl'$. We assume that ${\bf T}^{^{1,3}}_{st}$ is {\it locally below} $\CL$ and ${\bf T}^{^{1,4}}_{st}$ {\it locally above} $\CL$, thus, in a neighbourhood of $P_{1,M}$, these trajectories belong also to the regions {\bf C} and {\bf A}. In particular the trajectory
${\bf T}^{^{1,1}}_{st}$  cannot cross $\CL$, neither the segment $\{(x,y):x=X_{1M}, 0\leq y\leq Y_{1M}\}$. Hence either it converges to $(0,0)$ when $t\to-\infty$, or it crosses the axis $\{y=0\}$ at some point with positive $x$-coordinate less than $X_{1,M}$.\medskip


\nind\Remark In the critical case $M=-\gm^*$ there holds $-\frac{2p}{p+1}\gm^*Y^{\frac{p-1}{p+1}}_{-\gm^*}=2K=2-N-L$, hence the characteristic polynomial is $T(\gl)=\gl(\gl+2-N)$. If $N\geq 3$ the characteristic values are $0$ and $N-2$, with respective corresponding eigenvectors $(1,\frac{2}{p-1})$ and $(1, \frac{2}{p-1}+2-N)$. There exists an invariant curve  $\Gg$ passing through $P_{1,M}=P_{2,M}$, tangent
to $(1, \frac{2}{p-1})$ by the center manifold theorem. If $N=2$ the two characteristic values are $0$ with the  eigenspace generated by $(1,\frac{2}{p-1})$ which is tangent to the central manifold (a curve) at $(0,0)$. \medskip

Next we look for the existence of limit cycles. Since $M<0$, we cannot argue using the convexity argument used in \rlemma{cycleM>0}. We use system $(\ref{I-3-6})$ which also has a convexity property,  as we will see it in the proof below.
\blemma{cycleM<0} 1- If $N\geq 3$ and $p\geq\frac{N+2}{N-2}$ and $M<0$, there is no cycle surrounding $P_M$. \\
2- If $N\geq 3$, $\frac{N}{N-2}<p<\frac{N+2}{N-2}$ and $\overline M \leq M<0$, there is no cycle surrounding $P_M$.\\
3- If $N\geq 3$ and $1<p<\frac{N}{N-2}$ or $N=2$ and $\overline M \leq M\leq -\gm^*$, there is no cycle surrounding $P_{2,M}$.
\es
\Proof For 1, assume that there exists a periodic trajectory $(x(t), y(t)$ surrounding $P_M$. By Green's formula,
 \bel{I-3-23}\BA {lll}\displaystyle
0=\myint{\gg}{} \!\!\!\!\!\!\!\phantom{'}{\rm O}\left(-H_2(x,y)dx+H_1(x,y)dy\right)=\int\!\!\int_\Gg\left(\myfrac{\partial H_1(x,y)}{\partial x}+ \myfrac{\partial H_2(x,y)}{\partial y}\right)dxdy\\[4mm]\displaystyle
\phantom{0=\myint{\gg}{} \!\!\!\!\!\!\!\phantom{'}{\rm O}\left(-H_2(x,y)dx+H_1(x,y)dy\right)}
=\int\!\!\int_\Gg\left(\myfrac{2pM}{p+1}|y|^{\frac{p-1}{p+1}}-L\right)dxdy.
\EA \ee
Since $L\geq 0$  and $M<0$ we obtain a contradiction. This proves the first assertion. \\
For 2 and 3 we assume that $\gg_0\subset${\bf Q} is a cycle surrounding $P_M$ or $P_{2,M}$ and we denote by $\Gg_0$ the bounded domain bordered by $\gg_0$. By \rprop{1point} and \rprop{2points}, $P_M$ and $P_{2,M}$ are sources. We use the system $(\ref{I-3-6})$. Setting $J(\gs,z)=(\gs^pz)^\frac{1}{p+1}$, it becomes
          \bel{M12}\BA {lll}\displaystyle
\gs_{t}=\gs\left(\gs+2-N+z+MJ(\gs,z)\right):=F(\gs,z)\\[1mm]
z_{t}=z\left(N-p\gs-z-MJ(\gs,z)\right):=G(\gs,z).
             \EA\ee
             In the phase plane $(\gs,z)$ the equilibrium $P_M$ becomes
             $${\CP}_{\!M}=(\gs_M,z_M)=\left(\tfrac{Y_M}{X_M},\tfrac{X^p_M}{Y_M}\right)=\left(\tfrac{2}{p-1},\left(\tfrac{p-1}{2}\right)^pY_M^{p-1}\right).$$
Similarly
   $P_{2,M}$ becomes ${\CP}_{\!2,M}=(\gs_{2,M},z_{2,M})=\Bigl(\frac{Y_{2,M}}{X_{2,M}},\frac{X^p_{2,M}}{Y_{2,M}}\Bigr)=
   \Bigl(\tfrac{2}{p-1},\left(\tfrac{p-1}{2}\right)^pY_{2,M}^{p-1}\Bigr)$ and ${\CP}_{\!M}$ and ${\CP}_{\!2,M}$ are sources. Hence, any trajectory converging to this point ${\CP}_{\!M}$ (or ${\CP}_{\!2,M}$) admits an omega-limit cycle $\gg\subset \overline\Gg_0$. Consequently this cycle is not unstable which implies that the Floquet integral relative to this T-periodic solution is nonpositive:
              \bel{M13}\BA {lll}\displaystyle
\CI:=\myint{0}{T}\left(\myfrac{\partial F(\gs,z)}{\partial \gs}+\myfrac{\partial G(\gs,z)}{\partial z}\right)(t)dt\leq 0.
             \EA\ee
 The computation gives
 $$\BA {lll}
 \myfrac{\partial F}{\partial \gs}+\myfrac{\partial G}{\partial z}=(2-p)\gs+2-z+M\myfrac{2p+1}{p+1}(\gs^pz)^\frac{1}{p+1}-M\myfrac{p+2}{p+1}
(\gs^pz)^\frac{1}{p+1} \\[2mm]
\phantom{ \myfrac{\partial F}{\partial \gs}+\myfrac{\partial G}{\partial z}}
=(2-p)\gs+2-z+M\myfrac{p-1}{p+1}J(\gs,z).
\EA$$
 We set $\gs=\gs_M+\overline\gs$ and $z=z_M+\overline z$ (the computation would be the same with $(\gs_M,z_M)$ replaced by $(\gs_{2,M},z_{2,M})$), then
 $$\BA {lll}
 \myfrac{\CI}{T}= \myfrac{1}{T}\myint{0}{T}\left( \myfrac{\partial F}{\partial \gs}+\myfrac{\partial G}{\partial z}\right)(\gs(t),z(t))dt\\[4mm]
 \phantom{ \myfrac{\CI}{T}}
 =(2-p)\left(\gs_M+\myfrac{1}{T}\myint{0}{T}\overline\gs(t)dt\right)+2-z_M-\myfrac{1}{T}\myint{0}{T}\overline z(t)dt.
 +\myfrac{p-1}{p+1}\myfrac{M}{T}\myint{0}{T}J(\gs(t),z(t))dt.
 \EA$$
 By addition
 $$\BA {lll}
 \myfrac{\gs_{t}}{\gs}+ \myfrac{z_{t}}{z}=\gs-(N-2)+z+N-p\gs-z=2-(p-1)\gs\\[1mm]
 \phantom{ \myfrac{\gs_{t}}{\gs}+ \myfrac{z_{t}}{z}}
 =2-(p-1)\gs_M-(p-1)\overline\gs=-(p-1)\overline\gs.
 \EA$$
Integrating on a period, we get  $\int_{0}^{T}\overline\gs(t)dt=0$. We also derive from $(\ref{M12})$
$$\BA{lll}
0=\gs_M+\myfrac{1}{T}\myint{0}{T}\overline\gs(t)dt+2-N+z_M+\myfrac{1}{T}\myint{0}{T}\overline z(t)dt+\myfrac{M}{T}\myint{0}{T}J(\gs,z)dt\\[4mm]\phantom{0}
=\myfrac{2}{p-1}+2-N+z_M+\myfrac{1}{T}\myint{0}{T}\overline z(t)dt+\myfrac{M}{T}\myint{0}{T}J(\gs,z)dt
\\[4mm]\phantom{0}
=-K+z_M+\myfrac{1}{T}\myint{0}{T}\overline z(t)dt+\myfrac{M}{T}\myint{0}{T}J(\gs,z)dt.
\EA$$
Furthermore, $\gs_M+2-N+z_M+MJ(\gs_M,z_M)=0$. Indeed
              \bel{M14}\BA {lll}\displaystyle
\myfrac{1}{T}\myint{0}{T}\overline z(t)dt=-\myfrac{M}{T}\myint{0}{T}\left(J(\gs,z)-J(\gs_M,z_M)\right)dt.
             \EA\ee
             Next we show that the function $J$ is {\it concave} on ${\CQ}:=\{(\gs,z):\gs>0, z>0\}$: indeed
             $$J(\gs,z)-J(\gs_M,z_M)=\myfrac{\partial J}{\partial\gs}(\gs_M,z_M)\overline\gs+\myfrac{\partial J}{\partial z}(\gs_M,z_M)\overline z
             +\myfrac{1}{2}\left(a\overline\gs^2+2b\overline\gs\overline z+c\overline z^2\right),
             $$
  where $a$, $b$ and $c$ depend on $\gs$, $\gs_M$, $z$ and $z_M$ and:
  $$\myfrac{\partial J}{\partial\gs}(\gs,z)=\myfrac{p}{p+1}\gs^{-\frac{1}{p+1}}z^{\frac{1}{p+1}}\;\text{and }\; \myfrac{\partial J}{\partial z}(\gs,z)
  =\myfrac{1}{p+1}\gs^{\frac{p}{p+1}}z^{-\frac{p}{p+1}}.
  $$
  $$a=-\frac{p}{(p+1)^2}\gs_\gth^{-\frac{p+2}{p+1}}z_\gth^{\frac{1}{p+1}}\,,\; b=\frac{p}{(p+1)^2}\gs_\gth^{-\frac{1}{p+1}}z_\gth^{\frac{p}{p+1}}
  \;\text{and }\; c=-\frac{p}{(p+1)^2}\gs_\gth^{\frac{p}{p+1}}z_\gth^{-\frac{2p+1}{p+1}},
  $$
  where $(\gs_\gth,z_\gth)=(\gth\gs+(1-\gth)\gs_M,\gth z+(1-\gth)z_M)$ for some $\gth\in (0,1)$. Since
  $$b^2-ac=\myfrac{p^2}{(p+1)^4}\left(\gs_\gth^{-\frac{2}{p+1}}z_\gth^{-\frac{2}{p+1}}-\gs_\gth^{-\frac{2}{p+1}}z_\gth^{-\frac{2p}{p+1}}\right)=0,
  $$
  then $a\overline\gs^2+2b\overline\gs\overline z+c\overline z^2=a\left(\overline\gs+\frac{b}{a}\overline z\right)^2=-2R^2(t)\leq 0$ and from $(\ref{M14})$,
  $$\BA {lll}\displaystyle
\myfrac{1}{T}\myint{0}{T}\overline z(t)dt=-\myfrac{M}{T}\myfrac{\partial J}{\partial\gs}(\gs_M,z_M)\myint{0}{T}\overline \gs (t)dt-\myfrac{M}{T}\myfrac{\partial J}{\partial z}(\gs_M,z_M)\myint{0}{T}\overline z(t)dt+\myfrac{M}{T}\myint{0}{T}R^2(t)dt\\[4mm]
\phantom{\myfrac{1}{T}\myint{0}{T}\overline z(t)dt}
=-\myfrac{M}{T(p+1)}\gs^{\frac{p}{p+1}}_Mz_M^{-\frac{p}{p+1}}\myint{0}{T}\overline z(t)dt+\myfrac{M}{T}\myint{0}{T}R^2(t)dt,
  \EA$$
  thus
                \bel{M15}\BA {lll}\displaystyle
\left(1+\myfrac{M\gs^{\frac{p}{p+1}}_Mz_M^{-\frac{p}{p+1}}}{p+1}\right)\myint{0}{T}\overline z(t)dt=M\myint{0}{T}R^2(t)dt.
             \EA\ee
When $M<0$ we have already seen that $Y_M>y_{0,M}:=\left(\frac{2}{p-1}\left(\frac{-M}{p+1}\right)^{\frac{1}{p}}\right)^{\frac{p+1}{p-1}}$, defined in $(\ref{1-3-10a})$ (resp. $Y_{2,M}>y_{0,M}$), which yields
                 \bel{M16}\BA {lll}\displaystyle
1+\myfrac{M\gs^{\frac{p}{p+1}}_Mz_M^{-\frac{p}{p+1}}}{p+1}=\left(Y_M^{\frac{p(p-1)}{p+1}}-y_{0,M}^{\frac{p(p-1)}{p+1}}\right)Y_M^{-\frac{p(p-1)}{p+1}}>0\Longrightarrow \myint{0}{T}\overline z(t)dt<0.
             \EA\ee
             Therefore
                $$\BA {lll}\displaystyle
\myfrac{\CI}{T}=(2-p)\gs_M+2-z_M-\myfrac{1}{T}\myint{0}{T}\overline z(t)dt+\myfrac{M(p-1)}{T(p+1)}\myint{0}{T}J(\gs,z)dt\\[4mm]
\phantom{\myfrac{\CI}{T}}
=(2-p)\gs_M+2-z_M-\myfrac{1}{T}\myint{0}{T}\overline z(t)dt+\myfrac{M(p-1)}{p+1}J(\gs_M,z_M)\\[4mm]
\phantom{\myfrac{\CI}{T}---------------}
+\myfrac{M(p-1)}{T(p+1)}\myint{0}{T}\left(J(\gs,z)-J(\gs_M,z_M)\right)dt.
       \EA$$
       Hence from $(\ref{M14})$
       \bel{M17}\BA {lll}
\myfrac{\CI}{T}
=(2-p)\gs_M+2-z_M-\myfrac{2p}{T(p+1)}\myint{0}{T}\overline z(t)dt+\myfrac{M(p-1)}{p+1}J(\gs_M,z_M)\\[4mm]
\phantom{\myfrac{\CI}{T}}
=\myfrac{2}{p-1}-\left(\myfrac{p-1}{2}\right)^{p}Y_M^{p-1}-\myfrac{2p}{T(p+1)}\myint{0}{T}\overline z(t)dt+\myfrac{M(p-1)}{p+1}Y_M^{\frac{p-1}{p+1}}\\[4mm]
\phantom{\myfrac{\CI}{T}}
=\myfrac{2p}{p+1}MY_M^{\frac{p-1}{p+1}}-L-\myfrac{2p}{T(p+1)}\myint{0}{T}\overline z(t)dt.
             \EA\ee
But from \rprop{fixed-p},
\bel{M18}\overline  M<M<0\Longrightarrow \myfrac{2p}{p+1}MY_M^{\frac{p-1}{p+1}}-L=\myfrac{2p}{p+1}\left(MY_M^{\frac{p-1}{p+1}}-\overline M Y_{\overline M}^{\frac{p-1}{p+1}}\right)>0.
\ee
Combining $(\ref{M16})$ and $(\ref{M18})$ we obtain $\CI>0$ which contradicts the nonpositivity of the Floquet integral in the case $M>\overline  M$. Finally,
if $\overline  M=M$, then $\int_{0}^{T}\overline z(t)dt=0$, which implies by $(\ref{M15})$ that $R\equiv 0$ on the omega-limit cycle $\gg$.
Using the expression of $R(t)$ we infer
$$\overline\gs+\myfrac{b}{a}\overline z\equiv 0\;\text{ where }\,\; \frac{a}{b}=\myfrac{1}{z_\gth\gs^{\frac{p-1}{p+1}}_\gth}\;\text{ and }\,\; (\gs_\gth,z_\gth)=(\gth\gs+(1-\gth)\gs_M,\gth z+(1-\gth)z_M),
$$
for $\gth=\gth(t)\in (0,1)$. This is clearly impossible if one considers points at the intersection of $\gg$ and the straight line $\CL$. The same argument holds when $N=2$ where $\overline M=-\gm^*$.  \qeda\medskip

\nind\Remark The above proof can easily be adapted to recover the second statement of \rlemma{cycleM>0}. Indeed, if $p>\frac{N+2}{N-2}$ and
$0<M<\overline M$, $P_M$ is a sink. Hence if there is a cycle surrounding $P_M$, we can assume that it is an $\ga$-limit cycle, say $\gg,$ and the integral $\CI$ given by $(\ref{M13})$ is nonnegative by Floquet's theory. The inequality $(\ref{M15})$ yields $\int_0^T\overline z(t) dt> 0$ since
$1+\myfrac{M\gs^{\frac{p}{p+1}}_Mz_M^{-\frac{p}{p+1}}}{p+1}>1$. For $M<\overline M$  there holds $\frac{2p}{p+1}MY_M^{\frac{p-1}{p+1}}<L$ from the monotonicity of $M\mapsto MY_M^{\frac{p-1}{p+1}}$, which contradicts the sign of $\CI$.\medskip

In the next statement we extend \cite[Prop. 5.6]{ChWe}, which was proved in the case $1<p<\frac{N}{N-2}$, but valid actually for any $p>1$, and give precision on the behavior of the solutions. The constant
$$\gm^*(1)=\left(\frac{(p+1)^{2p+1}}{(2p)^p}\right)^{\frac 1{p+1}}$$
plays an important role. Recall also that $\gm^*(1)>\gm^*:=\gm^*(N)$ for $N\geq 2$.
\blemma {m1} Let $p>1$. If $N\geq 2$ and $M\leq -\gm^*(1)$ or $N=1$ and $M< -\gm^*(1)$, then there exists a ground state $u$. Furthermore there holds
\bel{M19}
u_{r}^2(r)<\myfrac{2p}{p+1}u^{p+1}(r)\qquad\text{for all }\, r>0.
\ee
As a consequence, the corresponding trajectory ${\bf T}_{reg}=\{(x_{reg}(t),y_{reg}(t))\}_{t\in\BBR}$ does not converge to $(0,0)$ when $t\to\infty$. If $1<p<\frac N{N-2}$ and $N\geq 3$ or $p>1$ and $N=1,2$, ${\bf T}_{reg}$ does not converge to $P_{1,M}$ when $t\to\infty$.
\es
\Proof Let $u$ be a regular solution of $(\ref{I-1})$ with $u(0)=1$ and $u_{r}(0)=0$. As in \cite{ChWe} we set
$$G(r)=\myfrac{1}{2}\left(u_{r}^2(r)-\myfrac{2p}{p+1}|u|^{p+1}(r)\right).
$$
Then $G(0)<0$ and $u>0$ on some maximal interval $[0,r_1)$ with $r_1\leq\infty$. If there exists a minimal $r_0\leq r_1$ such that $G(r_0)=0$, then $G_r(r_0)\geq 0$. From $(\ref{I-1})$ we have
$$\BA {lll}
G_r(r_0)=u_{r}(r_0)\left(u_{rr}(r_0)-pu^{p}(r_0)\right)\\[2mm]
\phantom{G'(r_0)}
=-\myfrac{N-1}{r_0}u_{r}^2(r_0)-u_{r}(r_0)\left((p+1)u^{p}(r_0)+M\abs{u_{r}(r_0)}^{\frac{2p}{p+1}}\right).
\EA$$
Since $G(r_0)=0$, $u^p(r_0)=\left(\frac{p+1}{2p}\right)^{\frac{p}{p+1}}\abs{u_{r}(r_0)}^{\frac{2p}{p+1}}$, hence
$$G_r(r_0)=-\myfrac{N-1}{r_0}u_{r}^2(r_0)-u_{r}(r_0)\abs{u_{r}(r_0)}^{\frac{2p}{p+1}}\left(\gm^*(1)+M\right).
$$
Since $u_{r}(r_0)<0$ we obtain a contradiction. Therefore $G(r)<0$ for any $r\in (0,r_0)$. By continuity this implies  $(\ref{M19})$ and in particular $u(r)>0$, hence $u$ is a ground state. Inequality $(\ref{M19})$ implies
\bel{M20}
u(r)>\myfrac{1}{\left(1+\frac {p-1}{2}\sqrt{\frac{2p}{p+1}}\;r\right)^{\frac{2}{p-1}}}\qquad\text{for all }\, r>0,
\ee
from which follows
\bel{M21}\displaystyle
\liminf_{r\to\infty}r^{\frac{2}{p-1}}u(r)\geq \left(\myfrac{2(p+1)}{p(p-1)^2}\right)^{\frac{1}{p-1}},
\ee
This implies that the trajectory ${\bf T}_{reg}$ does not converge to $(0,0)$ at infinity (a result which was clear when $1<p<\frac{N}{N-2}$ in which case $(0,0)$ is a source).  For the last statement we have with the equation   $(\ref{I-3-9-})$ satisfied by the equilibrium $X_M$ and the fact that $M\leq -\gm^*(1)$,
$$\BA {lll}\tilde f_M\left(\left(\myfrac{2(p+1)}{p(p-1)^2}\right)^{\frac{1}{p-1}}\right)\leq
\myfrac{2(p+1)}{p(p-1)^2}-\gm^*(1)\left(\myfrac{2}{p-1}\right)^{\frac{2p}{p+1}}\left(\myfrac{2(p+1)}{p(p-1)^2}\right)^{\frac{1}{p+1}}-\myfrac{2K}{p-1}
\\[4mm]
\phantom{\tilde f_M\left(\left(\myfrac{2(p+1)}{p(p-1)^2}\right)^{\frac{1}{p-1}}\right)}
\leq -\myfrac{2(N-1)}{p-1}.
\EA$$
Since $\tilde f_M$ has two roots $0<X_{1,M}<X_{2,M}$  by \rprop{fixed-p} and $M\leq-\gm^*(1)<-\gm^*(N)$, and we have $\left(\frac{2(p+1)}{p(p-1)^2}\right)^{\frac{1}{p-1}}>X_{1,M}$, the result follows.
\qeda\medskip

Next we give an alternative proof of a result of \cite{FiQu}.

\blemma{sim} Let $N=1,2$ and $p>1$ or $N\geq 3$ and $1<p<\frac{N}{N-2}$. If $-\gm^*<M\leq 0$ there exists no ground state.
\es
\Proof If $-\gm^*<M<0$, the only equilibrium of $(\ref{I-3x})$ is $(0,0)$ and it is a source by Section 3.4.1. If there exists a ground state then the trajectory ${\bf T}_{reg}$ remains in {\bf Q}. It is bounded by \rprop{bound} hence it admits an $\gw$-limit set which either contains an equilibrium or is a periodic orbit. The two possibilities are excluded. \qeda\medskip

Next we study the case  $M\leq  -\gm^*$. In particular we cover the case $M=  -\gm^*$ studied in \cite{Vo} with a different proof, using the energy function $\CZ$ defined in $(\ref{I-21.1})$.


\blemma{tra} Let $N\geq 3$ and $1<p<\frac{N}{N-2}$, or  $N=2$ and $p>1$. If $M\leq-\gm^*$ then,\smallskip

\nind(i) $\bf{T}_{reg}$ cannot converge to $P_{1,M}$ as $t\to\infty$ with $t\mapsto x(t)$ increasing,\smallskip

\nind(ii) $\bf{T}_{reg}$ cannot intersect $\CL$ at a point between $(0,0)$ and $P_{1,M}$.
\es
\Proof The function $\CZ$ defined by $(\ref{I-21.1})$ satisfies $(\ref{I-21.3})$ and $(\ref{I-21.4})$. Then
$\CU>0$ as soon as
$$\abs{M}y^{\frac{p-1}{p+1}}\leq \abs L\frac{(p+1)^2}{p+3}.
$$
For $M\leq-\gm^*$, there holds from \rprop{fixed-p} and $(\ref{I-3-10})$,
$$\abs{M}Y_{1,M}^{\frac{p-1}{p+1}}\leq\gm^*Y_{-\gm^*}^{\frac{p-1}{p+1}}=\frac{p+1}{p}\abs K<\abs L\frac{(p+1)^2}{p(p+3)}.
$$
Indeed we check that
$$(p+1)\abs L-(p+3)\abs K=2(N-1).
$$
Hence $\CU>0$ for any $y\in [0,Y_{1,M}]$.\smallskip

\nind (i) Suppose that $\bf{T}_{reg}$ converges to $P_{1,M}$ with $x(t)$ increasing, hence $(x(t),y(t))$ is below $\CL$ and $y(t)$ is also increasing. Then from $(\ref{I-21.3})$ $r\mapsto e^{\frac{2p\abs M}{p+1}\int_0^r\abs{u_r}^{\frac{p-1}{p+1}}ds}\CZ(r)$
is increasing. But $\CZ(0)=0$ and $\displaystyle\lim_{r\to\infty}\CZ(r)=0$ from $(\ref{I-21.2})$ since $L<0$ and $x$ and $y$ are bounded. This is a contradiction.\smallskip

\nind (ii) Suppose that $\bf{T}_{reg}$ intersects $\CL$ at a point $(\tilde x_{^M},\tilde y_{^M})$ between $(0,0)$ and $P_{1,M}$, i.e. $\tilde y_{^M}<Y_{1,M}$.\\
(a) If $M<-\gm^*$, then consider the stable trajectory ${\bf T}^{1,1}_{st}$ at $P_{1,M}$ which is below $\CL$: ${\bf T}^{1,1}_{st}$ cannot converge to $(0,0)$ when $t\to-\infty$; indeed it would be a unstable trajectory at the source point $(0,0)$ and  since  ${\bf T}_{reg}$ is the unique fast unstable trajectory at  this point (see Section 2.4.1), it is below  ${\bf T}^{1,1}_{st}$ near zero and the two curve would intersect. Therefore ${\bf T}^{1,1}_{st}$ leaves {\bf Q} through the semi-axis $\{(x,0):x>0\}$ at some
at some $x=x(\gt)$ and $(x(t),y(t))\in {\bf Q}$ for $t>\gt$. Setting $\bar r=e^\gt$, there holds
$\CZ(\bar r)=\bar r^au^{p+1}(\bar r)>0$ and again $r\mapsto e^{\frac{2p\abs M}{p+1}\int_{r}^{\infty}\abs{u_r}^{\frac{p-1}{p+1}}ds}\CZ(r)$ is increasing with limit $0$ as $r\to\infty$, a contradiction.\\
(b) If $M=-\gm^*$ and $\bf{T}_{reg}:=\bf{T}^{{-\gm^*}}_{reg}$ intersects $\CL$ at some point between $(0,0)$ and
$P_{-\gm^*}$, then by continuity and transversality, $\bf{T}^{{_M}}_{reg}$ intersects $\CL$ at some point between $(0,0)$ and  $P_{1,M}$ if $-\gm^*-\ge\leq M<-\gm^*$ provided $\ge>0$ is small enough. This contradicts (a). \qeda
\bprop{cross} Let $N\geq 3$ and $1<p<\frac{N}{N-2}$ or $N=2$ and $p>1$.\smallskip

\nind (i)  If $M<-\gm^*$ and  ${\bf T}_{reg}$ does not converge to $P_{2,M}$, then it intersects the line $\CL$ at some point $(x,y)$ with $x>X_{2,M}$. \smallskip

\nind (ii)  If $M=-\gm^*$ then ${\bf T}_{reg}$ intersects $\CL$ at some point $(x,y)$ with $x>X_{-\gm^*}$ and leaves {\bf Q}; there is no ground state.
\es
\Proof (i) Suppose $M<-\gm^*$. If ${\bf T}_{reg}$ remains below $\CL$, then $x_t>0$, $y_t>0$ and ${\bf T}_{reg}$
converges to $P_{1,M}$ or $P_{2,M}$ The first limit is excluded by \rlemma{tra} and the second by assumption. Hence
 ${\bf T}_{reg}$ intersects $\CL$. This intersection cannot occur between $(0,0)$ and $P_{1,M}$ and between
$P_{1,M}$ and $P_{2,M}$ since the vector field {\bf H} is inward in the region {\bf B} on this segment, so it occurs at some point
$(x,y)$ with $x>X_{2,M}$.\medskip

\nind (ii) Suppose $M=-\gm^*$. By \rlemma{tra}, ${\bf T}_{reg}$ intersects $\CL$ at some point $(x,y)$ with $x>X_{-\gm^*}$ and it enters in  the region {\bf D} (note that the region {\bf B} is empty).
Consider the slope $\gs=\frac{y}{x}$ of the trajectory ${\bf T}_{reg}$. As long as ${\bf T}_{reg}$ stays under $\CL$,
$\gs>\frac{2}{p-1}$. Now we introduce the system $(\ref{I-3-5})$, $(\ref{I-3-6})$ and set $\gk=\frac{z}{\gs}=\frac{x^{p+1}}{y^2}$. Hence
$$\BA {lll}
\myfrac{\gs_t}{\gs}=\gs+z-\gm^*\gs^{\frac{p}{p+1}}z^{\frac{1}{p+1}}+2-N
=\gs\phi(\gk)+2-N,
\EA$$
with
$$\phi(\gk)=1+\gk-\gm^*\gk^{\frac{1}{p+1}}.
$$
The function $\phi$ achieves its minimum at $\gk^*=\left(\frac{\gm^*}{p+1}\right)^{\frac{p+1}{p}}$ and
$\phi(\gk^*)=\frac{(N-2)(p-1)}{2}$. Thus
$$
\myfrac{\gs_t}{\gs}\geq (N-2)\left(\gs-\frac{p-1}{2}\right)\geq 0.
$$
Thus $\gs$ is nondecreasing. Therefore, after crossing $\CL$, $\gs>\frac 2{p-1}$ and ${\bf T}_{reg}$ cannot converge to  $P_{-\gm^*}$.\qeda

\subsection{Ground states and large solutions when $M<0$}
The following result extends \cite[Theorem B']{BiGaVe2} to a larger class of parameters $M$ in the radial case. We recall that $\gm^*(1)$ and $\gm^*(2)$ are defined in $(\ref{I-6})$.
\bprop{largesol} Assume $N\geq 2$, $p>1$ and $M> -\gm^*(1)$. Then there exists no positive solution of $(\ref{I-1})$ in $(a,b)$ for $a<b$ tending to infinity at $r=a$.
\es
\Proof Without loss of generality we can assume $a=1$. If $M\geq -\gm^*(2)$ the result follows from \cite[Theorem B']{BiGaVe2}, hence we can assume
$-\gm^*(1)<M<-\gm^*(2)$. We put $m=\abs M>0$ and
$$\ga=\ga(m)=\left(p^p\left(\myfrac{m}{p+1}\right)^{p+1}-1\right)>0.
$$
If $u$ satisfies $(\ref{I-1})$ and blows-up at $r=1$, then $v=\ln u$ satisfies
\bel{Jo1}\BA {lll}\displaystyle-v_{rr}-\frac{N-1}{r}v_r+\ga e^{(p-1)v}\geq 0\qquad\text{in }\,(1,b)\\[2mm]\phantom{-----;;---}
\!\displaystyle\lim_{r\to 1}v(r)=\infty.
\EA\ee
Up to changing $b$, we can assume that $u(b)\geq 1$, it follows that $v$ is bounded from below on $(1,b)$ by the solution of
\bel{Jo2}\BA {lll}\displaystyle-w_{rr}-\frac{N-1}{r}w_r+\ga e^{(p-1)w}\geq 0\qquad\text{in }\,(1,b)\\[2mm]\qquad\;\;
\!\displaystyle\lim_{r\to 1}w(r)=\infty\qquad v(b)=0.
\EA\ee
It is classical, (see e.g. \cite{VaVe}) that $\abs{w(r)-\frac{2}{p-1}\ln\left(\frac{1}{r-1}\right)}$ remains bounded on  $(1,b)$. Returning to the variable $(x(t),y(t))$ solutions of $(\ref{I-3x})$, then $x(t)\geq ct^{-\frac{2}{p-1}}$ on  $(0,\ln b)$. \smallskip

\nind (i) We first observe that $t^{\frac{2}{p-1}}x(t)$ and
$t^{\frac{p+1}{p-1}}y(t)$ remain bounded on $(0,\ln b)$. Indeed, by the equation satisfied by $y$ we get
$$\left(e^{Kt}y(t)\right)_t\geq -me^{Kt}\abs{y(t)}^{\frac{2p}{p+1}}\geq -c\abs{e^{Kt}y(t)}^{\frac{2p}{p+1}},
$$
and there exists a sequence $\{t_n\}$ converging to $0$ such that $y(t_n)\to\infty$. Furthermore $z(t)=e^{Kt}y(t)$ satisfies
$$-\left(\myfrac{p+1}{p-1}|z|^{-\frac{p-1}{p+1}}sgn(z)\right)'\geq c_1\Longrightarrow |z(t)|^{-\frac{p-1}{p+1}}sgn(z(t))\leq c_2t
\Longrightarrow z(t)\geq c_3t^{-\frac{p+1}{p-1}}.
$$
Therefore $y(t)\geq c_3t^{-\frac{p+1}{p-1}}$. Using the equation $(\ref{I-3-2})$ we first observe that $x_{t}$ remains negative in a right neighbourhood of $0$ otherwise there would exists a sequence $\{t_n\}$ decreasing to $0$ where $x_{t}(t_n)=0$ and $x_{tt}(t_n)\geq 0$ yielding
$$(x(t_n))^p-\frac{K}{p-1}x(t_n)-m\left(\frac{2}{p-1}\right)^{\frac{2p}{p+1}}(x(t_n))^\frac{2}{p-1}\leq 0,
$$
which is impossible since $x(t_n)\to\infty$. Similarly $y_{t}$ remains negative on some interval $(0,\gt_1)$ otherwise there would exists a sequence $\{t'_n\}$decreasing to $0$ such that $y_{t}(t'_n)=0$ and $y_{tt}(t'_n)\geq 0$. Since
$$y_{tt}=-Ky_{t}+px^{p-1}x_{t}-\myfrac{2mp}{p+1}y^{\frac{p-1}{p+1}}y_{t},
$$
and $x_{t}\leq 0$ we derive a contradiction. Therefore $y(t)\to\infty$ as $t\to 0$,  $y_{t}\leq 0$ and
$$-Ky=x^p-my^{\frac{2p}{p+1}}\leq 0\Longrightarrow y(t)\geq c_4x(t)^{\frac{p+1}{2}},
$$
which implies
$$x_{t}\leq \myfrac{2x}{p-1}-c_4x^{\frac{p+1}{2}}\leq c_5x^{\frac{p+1}{2}}\Longrightarrow x(t)\leq c_6ct^{-\frac{2}{p-1}},
$$
near $t=0$. Using again $(\ref{I-3x})$ and the monotonicity of $x(t)$ and $y(t)$,
$$e^{-\tfrac{t}{p-1}}x(\tfrac{t}{2})\geq \myfrac{(p-1)\left(e^{-\frac{t}{p-1}}-e^{-\frac{2t}{p-1}}\right)}{2}y(t)\geq c_6ty(t),
$$
which yields $0\leq y(t)\leq c_7t^{-\frac{p+1}{p-1}t}$ near $t=0$.\smallskip

(ii) Next we set $\gk(\gt)=(r-1)^{\frac{2}{p-1}}u(r)$ with $\gt=\ln(r-1)$. Then $\gk$ satisfies on $(-\infty,\gb]$ for some $\gb\in\BBR$
  \bel{Jo3}\BA {lll}
  -\gk_{\gt\gt}+\myfrac{p+3}{p-1}\gk_{\gt}-\myfrac{2(p+1)}{(p-1)^2}\gk+\myfrac{(N-1)e^\gt}{1+e^\gt}\left(\myfrac{2}{p-1}\gk-\gk_{\gt}\right)
  \\[4mm]\phantom{-------------}
  -\gk^p+m\abs{\myfrac{2}{p-1}\gk-\gk_{\gt}}^{\frac{2p}{p-1}}=0.
  \EA\ee
and  $\gk$ and $\gk_\gt$ remain bounded from above and from below on $(-\infty,a]$. Therefore the limit set $\Gs$ of the corresponding trajectory at $-\infty$ is not empty and it is included in the limit set  $\tilde \Gs$ at $-\infty$, of some trajectory of a nonnegative function satisfying the autonomous system
    \bel{Jo4}\BA {lll}
  -\gf_{\gt\gt}+\myfrac{p+3}{p-1}\gf_{\gt}-\myfrac{2(p+1)}{(p-1)^2}\gf
  -\gf^p+m\abs{\myfrac{2}{p-1}\gf-\gf_{\gt}}^{\frac{2p}{p-1}}=0,
  \EA\ee
  Which is precisely equation $(\ref{I-3-2})$ in dimension 1. By the Poincar\'e-Bendixon theorem, $\tilde \Gs$ either contains an equilibrium or it is a limit cycle. If the limit set contains an equilibrium, say $W$, it is positive and satisfies
  $$-\myfrac{2(p+1)}{(p-1)^2}W-W^p+m\left(\myfrac{2}{p-1}W\right)^{\frac{2p}{p-1}}=0.
  $$
  Since $M>-\gm^*(1)$ the only nonnegative root is zero, which yields a contradiction. If the limit set is a cycle $\gg$, it is a subset of
  ${\bf Q}$. This imply that there  would exist an equilibrium in the region bordered by $\gg$, contradiction. This ends the proof.\qeda

  \bcor{p47} Assume $N\geq 2$, $p>1$ and $-\gm^*(1)<M\leq -\gm^*$. Then any solution $(x(t), y(t))$ of system $(\ref{I-3x})$ issued from a point $(x_0,y_0)\in${\bf Q} and staying in {\bf Q} for $t\in I_{(x_0,y_0)}\cap\BBR_-$ where $I_{(x_0,y_0)}$ is the maximal interval of existence of this solution is defined on $\BBR_-$ and is bounded therein.
  \es
  \Proof Consider any solution such that $(x(0),y(0)=(x_0,y_0)\in {\bf Q}$ and suppose its negative trajectory {\bf T}$_{-}$ is defined on some maximal interval $(\gth,0]$ with $\gth\in (-\infty,0)$ and thus unbounded. We first suppose that
  $t\mapsto x(t)$ is not monotone when $t\to\gth$. Then there exists a sequence $\{t_n\}$ decreasing to $\gth$ such that $x_t(t_n)=0$ and thus $a_n=(x(t_n),y(t_n))\in\CL$, $x_{tt}(t_n)=y_t(t_n)<0$ and $ \lim_{t_n\to\gth}x(t_n)=\infty$.  \\
  Consider now the regular trajectory {\bf T}$_{reg}$.
  If $p\geq\frac{N}{N-2}$ then, from \rlemma{tetra}, either {\bf T}$_{reg}$ converges to $P_M$ when $t\to\infty$, or it crosses $\CL$ at a point $(x^*,\frac{2}{p-1}x^*)$ with $x^*>X_M$. If $n$ is such that $x^*<x(t_n)$ we get a contradiction: indeed for $t<t_{n-1}$, {\bf T}$_{-}$
  stays in the region bordered from above by $\CL$ and {\bf T}$_{reg}$, so it cannot intersect $\CL$ at $a_n$. \\
  If $1<p<\frac{N}{N-2}$, we infer the same contradiction using \rprop{cross}.
  Therefore $x(t)$ decreases monotonically to $\infty$ when $t\to\gth$. Since $M>-\gm(1)$ we derive a contradiction.
  \\ Hence $\inf I_{(x_0,y_0)}=-\infty$.  By the same reasons as above, $t\mapsto x(t)$ is monotone decreasing  and
  $x(t)\to\infty$  when $t\to-\infty$, and $y(t)>\frac{2}{p-1}x(t)\to\infty$. Moreover there exists $\bar t<0$ such that $y_t(t)\leq 0$: indeed , if for some $t_0<\bar t$, $y_t(t_0)> 0$, then $(x(t_0), y(t_0))\in {\bf D}$ and necessarily $(x(t), y(t)$ remains in {\bf D} for
  $t\leq t_0$ because $x_t(t)<0$ implies that $(x(t), y(t)\in{\bf D}\cup{\bf A}$ and the backward trajectory cannot cross the curve $\CC$ where $y_t=0$; now this implies that $y_t(t)>0$ for $t\leq t_0$, a contradiction. Therefore $y_t(t)\leq 0$ for $t\leq \bar t$. Then
   $$\abs My^{\frac{2p}{p+1}}\geq x^p-\abs Ky\geq x^p-\frac{\abs M}{2}y^{\frac{2p}{p+1}}-c_M.
   $$
   Since $x(t)\to\infty$, we deduce for large $|t|$,
   $$\frac{\abs M}{2}y^{\frac{2p}{p+1}}\geq \frac{1}{2}x^p.
   $$
   Returning to the system $(\ref{I-3x})$, we get  for large $|t|$
   $$x_t=\frac{2}{p-1}x-y\leq \frac{2}{p-1}x-\abs M^{-\frac{2p}{p+1}}x^{\frac{p+1}{2}}\leq -c^2x^{\frac{p+1}{2}}.
   $$
  Since $\frac{p+1}{2}>1$, it is straightforward to check by integration  that a positive function $x$ satisfying the above differential inequality cannot be defined on a interval unbounded from below, which ends the proof.$\phantom{-----}$\qeda
\bprop{exepsilon} Assume $p>1$ and $-\gm^*(1)< M< -\gm^*(1)+\ge$ for $\ge>0$ small enough. Then there exists a ground state and
$(\ref{M19})$ still holds.
\es
\Proof Assume $-\gm^*(1)< M< -\gm^*(1)+\ge$ and the regular solution $u$ is not a ground state. Then using the notation of \rlemma{m1} there exists $r_0>0$ such that $G(r_0)=0$ and $G'(r_0)\geq 0$. Hence $u_{r}^2(r_0)=\frac{2p}{p+1}u_{r}^{p+1}(r_0)\neq 0$, thus
$$\myfrac{N-1}{r_0}\leq \abs{u_{r}(r_0)}^{\frac{p-1}{p+1}}\left(\gm^*(1)+M\right).
$$
Put $t_0=\ln r_0$, then $y^2(t_0)=\frac{2p}{p+1}x^{p+1}(t_0)$ and $N-1\leq \abs{y(t_0)}^{\frac{p-1}{p+1}}\left(\gm^*(1)+M\right)$. Equivalently
\bel{Lo1} y(t_0)\geq\left(\myfrac{N-1}{\ge}\right)^{\frac{p+1}{p-1}}
\text{ hence }\; x(t_0)\geq \left(\myfrac{2p}{p+1}\right)^{\frac{1}{p+1}}\left(\myfrac{N-1}{\ge}\right)^{\frac{2}{p-1}}.
\ee
The curve $\{(x,y):y^2=\frac{2p}{p+1}x^{p+1}\}$ cuts $\CL$ at a unique $S_0=(x_0,y_0)\in{\bf Q}$ and $y_0^{p-1}=\frac{p+1}{2p}\left(\frac{2}{p-1}\right)^{p+1}$. If $M=-\gm^*(1)$ and $p\geq \frac{N}{N-2}$ (resp. $1<p< \frac{N}{N-2}$), then $X_{-\gm^*(1)}>x_0$ (resp. $X_{2,-\gm^*(1)}>x_0$). Indeed
this follows from $(\ref{1-3-10a})$. In the same way, if $1<p< \frac{N}{N-2}$, then $X_{1,-\gm^*(1)}<x_0$. These configurations still hold if $\ge$ is small enough, i.e. $X_{-\gm^*(1)+\ge}>x_0$ if $p\geq \frac{N}{N-2}$ and $X_{1,-\gm^*(1)+\ge}<x_0<X_{2,-\gm^*(1)+\ge}$ if $1<p< \frac{N}{N-2}$. However, the regular trajectory associated to $M$, ${\bf T}_{reg}^{^{_{M}}}$ has a unique intersection with $\CL$, at a point $(x_{reg}(t_1),y_{reg}(t_1)$ where $x_{reg}(t_1)$ is maximal, and either $p\geq \frac{N}{N-2}$ and $x_{reg}(t_1)>X_{-\gm^*(1)+\ge}$, or $1<p< \frac{N}{N-2}$ and $x_{reg}(t_0)>X_{1,-\gm^*(1)+\ge}$. In both cases $t_1<t_0$, $y_{reg}(t_1)=\frac 2{p-1}x_{reg}(t_1)$ and
$x_{reg}(t_1)>\left(\frac{p+1}{2p}\right)^{\frac{1}{p+1}}\left(\frac{N-1}{\ge}\right)^{\frac{2}{p-1}}$. Now, for $M>-\gm^*(1)$ the trajectory
${\bf T}_{reg}^{^{_{M}}}$ remains above ${\bf T}_{reg}^{^{_{-\gm^*(1)}}}$ as long as they remain below $\CL$ by \rlemma{nointer}. Then we encounter two possibilities: either ${\bf T}_{reg}^{^{_{-\gm^*(1)}}}$ converges to $P_{-\gm^*(1)}$ (or $P_{2,-\gm^*(1)}$ if $1<p<\frac{N}{N-2}$) with $x(t)$ and $y(t)$ increasing, or ${\bf T}_{reg}^{^{_{-\gm^*(1)}}}$ crosses $\CL$ at some point $(\tilde x,\tilde y)$ depending only on $N$ and $p$. Both possibilities are ruled out if $\ge$ is small enough. Hence ${\bf T}_{reg}^{^{_{M}}}$ is a ground state if
$M\in (-\gm^*(1),-\gm^*(1)+\ge]$ for $\ge>0$ small enough and $G$ remains negative.$\phantom{--------------}$\qeda
\subsection{The case $M<0$, $N\geq 3$ and $p\geq \frac{N+2}{N-2}$}
\bth{fac} Assume $M<0$, $N\geq 3$ and $p\geq \frac{N+2}{N-2}$. Then there exist ground states $u$. Moreover they satisfy $u(r)\sim U_M(r)$ as $r\to\infty$.
\es
\Proof Let $(x_{reg}(t),y_{reg}(t))$ be the regular solution issued from $(0,0)$. By \rlemma{Lyap} the function $\CV$ defined in $(\ref{I-3-16})$ is increasing. Since it vanishes at $t=-\infty$ it is positive. If there exists some $t_0$ such that $x_{reg}(t_0)=0$ then $\CV(t_0)<0$, which is impossible, hence
$x_{reg}(t)>0$ for all $t$. Thus ${\bf T}_{reg}$ is a ground state; it is bounded by \rprop{bound}, and it  cannot converge to $(0,0)$ since
$\CV(t)>\CV(-\infty)=0$. From \rlemma{cycleM<0} there is no cycle , hence ${\bf T}_{reg}$ converges to $P_M$ which is a sink. \qeda\medskip

\nind\Remark The existence of a ground state was already obtained in \cite{SeZo} with the use of the function $\CZ$ defined in
$(\ref{I-21.1})$.
\subsection{The case $M<0$, $N\geq 3$ and $\frac{N}{N-2}< p< \frac{N+2}{N-2}$}
In what follows we give an improvement of \cite[Theorem A]{Vo} in which it is shown that in this range of exponent there exists some $\gw>0$
 such that $-\gw\leq -M<0$ there exists no ground state. The expression of $\gw$ is explicit (and not simple).
 \bth{nbar1point} Assume $N\geq 3$ and $\frac{N}{N-2}< p< \frac{N+2}{N-2}$. If $\overline M\leq M<0$ there exists no ground state. Moreover there exists a unique, up to similarity transformation,  positive solution  $u$ satisfying $u(r)\sim U_M(r)$ as $r \to 0$ such that $u(r)\sim cr^{2-N}$ ($c>0$) as $r\to\infty$.
 \es
 \Proof Suppose that there exists such a ground state, then  ${\bf T}_{reg}$ remains in $\overline {\bf Q}$. By \rprop{1point}, $P_M$ is either a source if $\overline M< M$ or a weak source if $\overline M= M$, and by \rlemma{cycleM<0} there exists no cycle surrounding $P_M$. By \rlemma{tetra}, $(x_{reg}(t),y_{reg}(t))$ converges to  $(0,0)$ when $t\to \infty$, hence ${\bf T}_{reg}$ is a homoclinic orbit equivalently ${\bf T}_{reg}={\bf T}_{st}={\bf T}_{unst}$. Now
 $$Trace \,D{\bf H}(0,0)=\frac{2}{p-1}-K=-L>0.
 $$
 Hence, by \cite[Th 9.3.3]{HuWe} the homoclinic orbit is repelling. Since $P_M$ is also repelling, we derive a contradiction because any trajectory issue from $\CB$ must converge to ${\bf T}_{reg}$. Hence ${\bf T}_{reg}$ intersects the axis $\{x=0\}$ for some positive $y_1>0$, and there exists no ground state. We denote by $\CO$ the region of {\bf Q} delimited by the regular trajectory ${\bf T}_{reg}$ and the segment $(0,y): 0<y<y_1$. It is negatively invariant. The stable trajectory ${\bf T}_{st}=\{(x_{st},y_{st})\}$ of $(0,0)$ satisfies $x_{st}(t)=ce^{-Kt}(1+o(1))$ and
 $y_{st}(t)=c(N-2)e^{-Kt}(1+o(1))$ when $t\to-\infty$, thus it remains in $\CO$. Because there are no cycle in $\CO$, it must converge to $P_M$, hence the corresponding $u_{st}$ is equivalent to $U_M$ near $r=0$, which ends to proof.\qeda \medskip

 \nind\Remark In the previous theorem the positive solution  $u$ satisfying $u(r)\sim U_M(r)$ as $r \to 0$ such that $u(r)\sim cr^{2-N}$ ($c>0$) as $r\to\infty$ is the stable trajectory ${\bf T}_{st}$. It is a heteroclinic orbit connecting $(0,0)$ to $P_M$.
 We conjecture that in the case $p=\frac N{N-2}$ the non existence of a ground state still holds and that there exists a unique solution $u$ such that satisfying $u(r)\sim U_M(r)$ as $r \to 0$ such that $u(r)\sim cr^{2-N}(\ln r)^{\frac{2-N}{2}}$ ($c>0$) as $r\to\infty$.

 \medskip

 The expression of the result presents some similarity with \rth{Mbar<M} in the case $M>0$. However a new type of difficulty appears: in order to define properly an intersection function expressing the distance  between some trajectories as in \cite{Bi2}, we need to find some values of the parameter $M$ for which there exists a ground state, and all the trajectories in {\bf Q} are bounded. It is not the case when $M<-\gm^*(1)$ even if there exists a ground state by \rlemma {m1} but we can easily prove that there exist  large solutions. So we need to prove that for
 $M=-\gm^*(1)+\ge$ there exist a ground state and all the trajectories in {\bf Q} are bounded. This is the object of \rprop{exepsilon} and
 \rprop{largesol}.
\bth{psuperM<0} Let $N\geq 3$ and $\frac{N}{N-2}<p<\frac{N+2}{N-2}$. There exist positive constants $\tilde\gm_{min}\leq\tilde\gm_{max}$, verifying $\overline\gm<\tilde\gm_{min}\leq\tilde\gm_{max}<\gm^*(1)$ such that\smallskip

\nind (i) If $M<-\tilde\gm_{max}$ there exist ground states $u$ such that $u(r)\sim U_M(r)$ or ondulating around $U_M(r)$ when $r\to\infty$.\smallskip

\nind (ii) If $M=-\tilde\gm_{max}$ or If $M=-\tilde\gm_{min}$ there exist ground states $u$ such that $u(r)\sim cr^{2-N}$, $c>0$, when $r\to\infty$.\smallskip

\nind (iii) If $-\tilde\gm_{min}<M<0$ there exists no ground state. Furthermore there exist singular solutions $u$ ondulating around $U_M(r)$ when $r>0$ and singular solutions $u$ ondulating around $U_M(r)$ when $r\to 0$ and such that $u(r)\sim cr^{2-N}$, $c>0$, when $r\to\infty$.
\es
\Proof Recall that $\overline M=-\overline\gm$. First we show that if $-\gm^*(1)<M<\overline M$, the stable trajectory ${\bf T}_{st}:={\bf T}^{_M}_{st}$ either has a limit cycle around $P_M$
or does not stay in {\bf Q}. If we assume that it stays in {\bf Q}, then it is bounded by \rcor{p47}. Since at $-\infty$ it cannot converge to $P_M$ which is a sink by \rprop{1point}, it admits a alpha-limit cycle which is a closed orbit around $P_M$. \smallskip

For $-\gm^*(1)<M< \overline M$ we denote by
 $P^M_{st}=(x^M_{st}, x^M_{st})$ the farthest point of the closure $\overline{{\bf T}_{st}}$ of the trajectory ${\bf T}_{st}$  belonging to the line $\CL$, i.e. the points with the largest $x$ (and $y$)-coordinate. We also denote by $P^M_{reg}=(x^M_{reg},x^M_{reg})$ the farthest point of the intersection of $\CL$ with the closure $\overline{{\bf T}_{reg}}$ of ${\bf T}_{reg}:={\bf T}^{_M}_{reg}$.  More precisely since either ${\bf T}_{st}$ leaves {\bf Q} or has an alpha-limit cycle around
$P_M$, in that case $P^M_{st}$ corresponds to the last intersection of ${\bf T}_{st}$ and $\CL$. If ${\bf T}_{reg}$ converges to $P_M$  monotonically, then $P^M_{reg}=P_M\in \overline{{\bf T}_{reg}}$, while if this convergence is not monotone, or if ${\bf T}_{reg}$ admits a omega-limit cycle around $P_M$, or if ${\bf T}_{reg}$ leaves {\bf Q}, $P^M_{reg}$ is  the first intersection of ${\bf T}_{reg}$ with $\CL$. Both the functions
$M\mapsto x^M_{st}$ and $M\mapsto x^M_{reg}$ are continuous, either by transversality argument or by the continuity of $M\mapsto X_M$.
Hence the function $M\mapsto g(M)=x^M_{reg}-x^M_{st}$ is continuous.  For $M<\overline M$ we encounter three possibilities:\smallskip

\nind  (i) $x^M_{reg}=X_M$ or ${\bf T}_{st}$ converges to $X_M$ non-monotonicaly, or ${\bf T}_{reg}$ has a omega-limit cycle around  $P_M$. In such a case ${\bf T}_{st}$ does not stay in {\bf Q}, thus $g(M)<0$.
\smallskip

\nind  (ii)  ${\bf T}_{reg}$ does not stay in {\bf Q}, then ${\bf T}_{st}$ belongs to the region of {\bf Q} bordered by ${\bf T}_{reg}$ and the axis
$\{x=0\}$. Then thus $g(M)>0$. \smallskip

\nind  (iii)  $g(M)=0$, then ${\bf T}_{st}={\bf T}_{reg}$ is a homoclinic orbit. \smallskip

\nind If $M=\overline M$ there exists no ground state by \rth{nbar1point} hence, by continuity, this still holds for $\overline M-\ge<M<\overline M$ for $\ge>0$ small enough and then $g(\overline M)>0$. By \rprop{exepsilon}, if $-\gm^*(1)<M<-\gm^*(1)+\ge$, there exists a ground state, hence
$g(M)<0$. Since $g$ is continuous there exists $M\in (-\gm^*(1),\overline M)$ such that $g(M)=0$.\smallskip

If  we define
\bel{tilde mu}\BA{lll}
\tilde \gm_{min}=\min\{\abs M\in (\overline \gm,\gm^*(1)):g(M)>0\}\\
\tilde \gm_{max}=\max\{\abs M\in (\overline \gm,\gm^*(1)):g(M)<0\}
\EA\ee
then the trajectories ${\bf T}_{reg}$ corresponding to $M=-\tilde \gm_{min}$ and $M=-\tilde \gm_{max}$ are homoclinic and they satisfy the statements (ii) of \rth{psuperM<0} and the conclusion follows.\qeda\medskip

The proof of Theorem B, B' follows from \rth{nbar1point} and \rth{psuperM<0}.

 \subsection{The case $M<0$ and $1<p<\frac{N}{N-2}$}
 We present first a general existence result of singular solutions.
 \bprop{T11} If $N\geq 3$,  $1<p<\frac{N}{N-2}$ and $M<-\gm^*$ there exists positive singular solutions $u$ such that
 $r^{N-2}u(r)\to c$ for some $c>0$ when $r\to 0$ and $r^{\frac{2}{p-1}}u(r)\to X_{1,M}$ when $r\to \infty$.
 \es
 \Proof By \rprop{2points} $P_{1,M}$ is a saddle point of system $(\ref{I-3x})$. By the remark after this proposition there exist two stable trajectories ${\bf T}_{st}^{1,j}$, $j=1,2$ converging to $P_{1,M}$ as $t\to\infty$; the trajectory ${\bf T}_{st}^{1,1}$ is locally below the line $\CL$ (see \rprop{cross}), hence it belongs to the region {\bf C} for $t>t_0$. By \rprop{cross} either the regular trajectory ${\bf T}_{reg}$ converges to  $P_{2,M}$ or it crosses the line $\CL$ beyond $P_{2,M}$. Hence ${\bf T}_{st}^{1,1}$ cannot {intersect} ${\bf T}_{reg}$, and is trapped when $t$ decreases in the region {\bf C} and the curve ${\bf T}_{reg}$. Thus it converges to $(0,0)$ when $t\to-\infty$. Because $(0,0)$ is a source
 (see Section 3.4.1) with one fast trajectory ${\bf T}_{reg}$, which satisfies
 $\displaystyle\lim_{t\to-\infty}e^{-\frac{2}{p-1}t}x_{reg}(t)=u(0)$, the trajectory ${\bf T}_{st}^{1,1}$ is a slow one and it satisfies
 $\displaystyle\lim_{t\to-\infty}e^{(N-2-\frac{2}{p-1})t}x_{reg}(t)=c$.\qeda\medskip

 \nind\Remark Under the assumptions of \rprop{T11} the trajectory ${\bf T}_{unst}^{1,4}$ leaves {\bf Q} since in this region it stays in the sector
 $\{(x,y):H_j(x,y)<0\}$ for $j=1,2$, and this sector contains no stable equilibrium. The result holds also if $M=-\gm^*$.
 \bth{bar2points} Let $N\geq 3$ and $1<p<\frac{N}{N-2}$. Then\\
 (i) if $\overline  M\leq M\leq 0$ there exists no ground state.\\
 (ii) if $\overline  M\leq M<-\gm^*$ there exists a positive singular solution $u$, unique up to scaling, such that $u(r)\sim U_{2,M}(r)$ as $r\to 0$, $u(r)\sim U_{1,M}(r)$ as $r\to \infty$ and $u(r)> U_{1,M}(r)$ for any $r>0$.
 \es
 \Proof (i) Let $\overline  M\leq M\leq 0$. Suppose that there exists a ground state ${\bf T}_{reg}$, then $\overline  M\leq M< -\gm^*$ by \rlemma{sim}, and by \rprop{cross}. Furthermore this trajectory is bounded by \rprop{bound}. Moreover $P_{2,M}$ is a source or a weak source and there exists no cycle surrounding it by  \rprop{2points} and \rlemma{cycleM<0}. Hence ${\bf T}_{reg}$ converges to an equilibrium which cannot be $(0,0)$ neither $P_{2,M}$. So, from \rprop{cross}, it converges to $P_{1,M}$ from above $\CL$ as $t\to\infty$. Since  $P_{1,M}$ is a saddle point ${\bf T}_{reg}$ must coincide with the stable trajectory  ${\bf T}_{st}^{1,2}$ of this point. Therefore the region bordered by ${\bf T}_{reg}$, ${\bf T}_{st}^{1,1}$ and $(0,0)$ is invariant and it contains only one source equilibrium and no cycle around $P_{2,M}$ by \rlemma{cycleM<0}. Any trajectory starting from this region must converge to $P_{1,M}$ which is impossible. Hence ${\bf T}_{reg}$ is not a ground state. \smallskip
 Since the trajectory ${\bf T}_{reg}$ intersects the axis $\{x=0\}$, the trajectory ${\bf T}_{st}^{1,2}$ which converges to $P_{1,M}$ at infinity from above $\CL$ is trapped in the region bordered by ${\bf T}_{reg}$ and the semi-axis $\{(0,y):y>0\}$ which is negatively invariant. As there is no cycle in this region, it converges to $P_{2,M}$ when $t\to\infty$. To this trajectory corresponds a solution $u$ of $(\ref{I-1})$ which satisfies $u(r)\sim U_{2,M}(r)$ when $r\to\infty$ and $u(r)\sim U_{1,M}(r)$ when $r\to 0$. Furthermore $u(r)> U_{1,M}(r)$ for all $r>0$.\qeda\medskip

 The next result extends \rprop{T11} to the case $M=-\gm^*$.

 \bprop{bar2points*}Let $N\geq 3$ and $1<p<\frac{N}{N-2}$. If $M=-\gm^*$ there exists positive solutions $u$ satisfying $u(r)\sim cr^{2-N}$
when $r\to 0$ for some $c>0$ and $u(r)\sim U_{-\gm^*}(r)$ when $r\to \infty$.
\es
\Proof Since $N\geq 3$, there still exist two unstable trajectories ${\bf T}_{unst}^{1,j}$ (j=3,4) starting from $P_{1,M}$ with slope $N-2$;  ${\bf T}_{unst}^{1,4}$, which is above $\CL$, leaves {\bf Q} through the semi-axis $\{(0,y):y>0\}$ (see the remark above). From \rprop{cross} the trajectory ${\bf T}_{reg}$ crosses $\CL$ at a unique point $A_{reg}=(x_{reg},\frac{2}{p-1}x_{reg})$ with $x_{reg}>X_{-\gm^*}$ and leaves {\bf Q} through the semi-axis $\{(0,y):y>0\}$ and its exit point is above the exit point of ${\bf T}_{unst}^{1,4}$. Hence the other unstable trajectory ${\bf T}_{unst}^{1,3}$ which is trapped in the region of {\bf Q} bordered by  ${\bf T}_{unst}^{1,4}$, ${\bf T}_{reg}$ and the semi-axis $\{(0,y):y>0\}$, either converges to $P_{-\gm^*}$ or leaves {\bf Q} crossing the semi-axis $\{(0,y):y>0\}$. Necessarily it intersects the line $\CL$ at some point $\overline A=(\overline x,\frac{2}{p-1}\overline x)$ with $\overline x>X_{-\gm^*}$ and enters the region $\bf D$. As long as it stays above $\CL$, as in the \rprop{cross}, there holds $\gs_t>0$ (with the notations from this proposition). Thus it cannot converges to $P_{-\gm^*}$. Therefore ${\bf T}_{unst}^{1,3}$ crosses the semi-axis $\{(0,y):y>0\}$ between the exit points of ${\bf T}_{unst}^{1,4}$ and ${\bf T}_{reg}$.\\
Let $\CR$ be the open connected region of {\bf Q} below the line $\CL$ and bordered by ${\bf T}_{reg}$ and ${\bf T}_{unst}^{1,3}$.  If $(\tilde x,\tilde y)\in \CR$ and $\tilde {\bf T}$ is the trajectory through this point,  we have three possibilities:\smallskip

\nind (i) Either $\tilde {\bf T}$ leaves $\CR$ crossing $\CL$ between $(0,0)$ and $P_{-\gm^*}$,\smallskip

\nind (ii) Or $\tilde {\bf T}$ leaves $\CR$ crossing $\CL$ between $P_{-\gm^*}$ and $\overline A$,

\nind (iii) Or $\tilde {\bf T}$ converges to $P_{-\gm^*}$. \smallskip

\nind The set of points satisfying (i) or (ii) are non-empty, disjoint and open. Therefore the set of points satisfying (iii) is non-empty and the corresponding trajectory $\tilde {\bf T}$ converges to $P_{-\gm^*}$ when $t\to\infty$. The backward trajectory remains in  $\CR$ which is negatively invariant. Since there is no fixed point in this region, it converges to $(0,0)$ when $t\to-\infty$ and it is a slow trajectory of this point, which ends the proof. \qeda

\medskip
Next we describe the behaviour of the positive solutions for $M<-\overline\gm$. However, in order to use the method introduced in the proof of \rth{psuperM<0} we are confronted to another difficulty namely that there can exist homoclinic trajectories at $P_{1,M}$.
\bth{psubM<0} Let $N\geq 3$, $M<0$ and $1<p<\frac{N}{N-2}$. Then there exist positive real numbers $\overline \gm<\hat\gm_{min}\leq \hat\gm_{max}<\tilde \gm_{min}\leq \tilde\gm_{max}<\gm^*(1)$  with the following properties \smallskip

\nind 1- for $\overline \gm<|M|<\tilde \gm_{min}$ there is no radial ground state;\smallskip

\nind 2- for $|M|=\tilde \gm_{min}$ or $|M|=\tilde \gm_{max}$ there exist ground states $u$ satisfying $u(r)\sim U_{1,M}(r)$ when $r\to\infty$;\smallskip

\nind 3- for $|M|>\tilde\gm_{max}$ there exist ground states either such that $u(r)\sim U_{2,M}(r)$ when $r\to\infty$ or ondulating around
$U_{2,M}(r)$ when $r\to\infty$.\smallskip

Moreover,\smallskip

\nind 4- For $\overline \gm<|M|<\hat\gm_{min}$ there exist  solutions $u$, necessarily singular, ondulating around  $U_{2,M}$ as $r\to 0$ and $u(r)\sim U_{1,M}(r)$ when $r\to\infty$ and solutions $u$ ondulating around  $U_{2,M}$ on $(0,\infty)$;
\smallskip

\nind 5- for $|M|=\hat\gm_{min}$ or $|M|=\hat\gm_{max}$ there exists a solution $u\neq U_{1,M}$ such that $u(r)\sim U_{1,M}(r)$ both when $r\to 0$ and $r\to\infty$;\smallskip

\nind 6- for $\hat\gm_{max}<|M|<\tilde \gm_{min}$ there exists a solution $u$ such that $u(r)\sim cr^{2-N}$ as $r\to 0$ and either
$u(r)\sim U_{2,M}(r)$ or ondulating around $U_{2,M}$ when $r\to\infty$.
\es
\Proof {\it Step 1}. For $M=\overline M$ we know the behaviour of the solutions from \rth{bar2points}. The trajectories ${\bf T}_{reg}$, ${\bf T}_{unst}^{1, 3}$
and ${\bf T}_{unst}^{1, 4}$ leave {\bf Q} on $\{(0,y):y>0\}$ with transverse intersections, and ${\bf T}_{unst}^{1, 2}$ connects $P_{2,M}$ to $P_{1,M}$. This transversality property is also true for the corresponding trajectories with parameter $\overline M-\ge_0\leq M\leq \overline M$ for $\ge_0>0$ small enough. Let $\overline M_\ge=\overline M-\ge$ where $0<\ge\leq\ge_0$. Then $P_{2,\overline M_\ge}$ is a sink by \rprop{2points}. By \rprop{cross}, the points $P_{1,\overline M_\ge}$ and $P_{2,\overline M_\ge}$
belong to the region $\CR_\ge$ bordered by ${\bf T}_{reg}$, and the semi-axis $\{(0,y):y>0\}$. The trajectory ${\bf T}^{1,2}_{st}$ which converges to
$P_{1,\overline M_\ge}$ at $\infty$, cannot converge to $P_{2,\overline M_\ge}$ at $-\infty$. Hence it has a limit cycle around $P_{2,\overline M_\ge}$. Note that ${\bf T}^{1,2}_{st}$ is included in the region bordered by ${\bf T}^{1,3}_{unst}$ and ${\bf T}^{1,4}_{unst}$.\smallskip

For $\ga>0$ small enough set $M^\ga=-\gm^*(1)+\ga$. By \rprop{exepsilon} the corresponding regular trajectory ${\bf T}_{reg}$ is a ground state. Since $(0,0)$ is a source and $(\ref{M19})$ holds, it cannot converge neither to $(0,0)$ nor to $P_{1,M^\ga}$ as in \rlemma{m1}. Hence either it converges to $P_{2,M^\ga}$
 at $\infty$ or it admits a limit cycle around. Next, ${\bf T}^{1,3}_{unst}$ is trapped in the positively invariant region bordered by ${\bf T}^{1,1}_{st}$ which connects $(0,0)$ to $P_{1,M}$, the portion of the curve $\CC$ below $\CL$ (hence between $P_{1,M}$ and $P_{2,M}$) and ${\bf T}_{reg}$ between $(0,0)$  and its second intersection with $\CC$; therefore, it converges to $P_{2,M^\ga}$ at $\infty$ or it admits a limit cycle too. Finally consider the trajectory ${\bf T}^{1,2}_{st}$ which tends to $P_{1,M^\ga}$ at $\infty$. It stays in the negatively invariant region $\{(x,y):y>\frac{2}{p-1}x\text{ or }x>X_{1,M}\}$. If it stays in {\bf Q}, then the solution is defined on $\BBR$ and the trajectory is bounded from \rcor{p47}. So, either it converges to a fixed point or it has a limit cycle around $P_{2,M^\ga}$ when $t\to-\infty$. This is impossible since it would intersect ${\bf T}_{reg}$. Hence for $M=M^\ga$, ${\bf T}^{1,2}_{st}$ leaves {\bf Q} in finite time at $(x(\bar t),0)$ for some $x(\bar t)>0$. Then there exist $t_0>t_1>\bar t$ such that $y_t(t_0)=0$ and $y(t_0)$ is the maximum of
 $y$ on $(\bar t,\infty)$ and $(x(t_1),y(t_1))\in\CL$. Hence $x(t)\leq x(t_1)$ for $\bar t\leq t\leq t_1$.\smallskip

 \nind{\it Step 2}. Next consider any $M\in (-\gm^*(1),\overline M)$. In any case the trajectories ${\bf T}^{1,2}_{st}$ and ${\bf T}_{reg}$ are bounded as long as the stay in {\bf Q} from \rcor{p47}. Then we define $(x^M_{reg},y^M_{reg})$ as the farthest point on $\CL$ belonging to the closure $\overline {\bf T}_{reg}$ of ${\bf T}_{reg}$, and $(x^M_{st},y^M_{st})$ as the farthest point on $\CL$ belonging to $\overline{\bf T}_{st}$.\smallskip

  Let $A$ be the set of $M\in (-\gm^*(1),\overline M)$ such that  there is no ground state, let $B_1$ be the set of $M\in (-\gm^*(1),\overline M)$ such that there exists a ground state converging to $P_{1,M}$ at $\infty$ and let $B_2$ be the set of $M\in (-\gm^*(1),\overline M)$
such that there exists a ground state converging to $P_{2,M}$ or having a limit cycle at $\infty$. Then $(-\gm^*(1),\overline M)=A\cup B_1\cup B_2$. Clearly $\overline M\in A$ and $M^\ga\in B_2$, furthermore if $M\in A\cup B_1\cup B_2$ we have three possibilities.\smallskip

\nind $\bullet$ Any $M\in B_2$ has the same properties as $M^\ga$, hence  ${\bf T}_{reg}$ is a ground state which converges to $P^{2,M}$ or has a limit cycle around $P^{2,M}$; ${\bf T}^{1,3}_{unst}$ either converges to $P^{2,M}$ or has a limit cycle around $P^{2,M}$; and ${\bf T}^{1,2}_{st}$ intersects $\CL$ at a last value $t_1$. \smallskip

\nind $\bullet$ If $M\in A$, ${\bf T}_{reg}$ is not a ground state and it leaves {\bf Q} through the semi-axis $\{x=0, y>0\}$; ${\bf T}^{1,2}_{st}$ and ${\bf T}^{1,3}_{unst}$ are included in the region of {\bf Q} bordered by ${\bf T}_{reg}$ and three configurations are possible:\\
A-(i) either ${\bf T}^{1,2}_{st}$ has a limit cycle around $P_{2,M}$ and ${\bf T}^{1,3}_{unst}$ leaves {\bf Q}\\
A-(ii) or ${\bf T}^{1,3}_{unst}$ converges to $P_{2,M}$ or has  a limit cycle around.\\
A-(iii) or ${\bf T}^{1,2}_{st}={\bf T}^{1,3}_{unst}$ which means this trajectory is homoclinic with respect to $P_{1,M}$.\\
Note that $\overline M_\ge$ satisfies A-(i).\smallskip

As in the proof of \rth{psuperM<0} the mappings
$M\mapsto x^M_{reg}$ and $M\mapsto x^M_{1,2}$ are continuous. We set $g(M)=x^M_{reg}- x^M_{1,2}$. If $M\in A$, then
$g(M)<0$ and if $M\in B_2$, then $g(M)>0$. Since $g$ is continuous there exists $M\in (-\gm^*(1),-\overline\gm)$ such that  $g(M)=0$. Hence  ${\bf T}_{reg}={\bf T}^{1,2}_{st}$ and $B_1\neq\emptyset$. More precisely we can define $\overline\gm<\tilde\gm_{min}<\tilde\gm_{max}<\gm^*(1)$ such that $\tilde\gm_{min}\in B_1$, $\tilde\gm_{max}\in B_1$. If $\abs M<\tilde\gm_{min}$, $g(M)>0$ and there is no ground state. If $\abs M>\tilde\gm_{max}$, $g(M)<0$ and there exists a ground state $u$ such that $r^{\frac{2}{p-1}}u(r)\to X_{2,M}$  or such that $r^{\frac{2}{p-1}}u(r)$ is turning around $X_{2,M}$ as $r\to\infty$.\smallskip

Finally we consider the relative position of ${\bf T}^{1,2}_{st}$ and ${\bf T}^{1,3}_{unst}$ when $M\in [-\tilde\gm_{min}, -\overline\gm)$. In any of the three situations A-(i), A-(ii) and A-(iii), $\CL$ intersects ${\bf T}^{1,3}_{unst}$ at a first point of x-coordinate $x^M_{1,3}$ and ${\bf T}^{1,2}_{st}$ at a last point of x-coordinate $x^M_{1,2}$. We define  the continuous function $M\mapsto h(M)=x^M_{1,3}- x^M_{1,2}$. Then $h(\overline M_\ge)>0$ and $h(-\tilde\gm_{min})=x^M_{1,3}-x^M_{reg}<0$. Hence there exists at least one
$M\in (-\tilde\gm_{min}, -\overline\gm)$ where $h(M)=0$ and for such a $M$, A-(iii) holds. Then we define $\hat\gm_{min}$ and $\hat\gm_{max}$ in $(\overline\gm,\tilde\gm_{min})$  such that A-(iii) holds if $\abs M=-\hat\gm_{min}$ or if $\abs M=-\hat\gm_{max}$. Hence, if $-\hat\gm_{min}<M<-\overline\gm$, $h(M)>0$ and the trajectory ${\bf T}^{1,3}_{unst}$ starts from $P_{1,M}$ and converges to $P_{2,M}$ or has a limit cycle around $P_{2,M}$. If $-\hat\gm_{max}<M<-\tilde \gm_{min}$, $h(M)<0$ and the trajectory ${\bf T}^{1,2}_{st}$ starts from $(0,0)$ with the slope $N-2$ and converges to $P_{1,M}$ when $t\to\infty$.\qeda\medskip

\nind{\it Proofs of Theorem C and C'. }They are a consequence of \rprop{T11}, \rth{bar2points}, \rprop{bar2points*}  and \rth{psubM<0}.\medskip

\nind\Remark It is an open problem whether the cycles which may exist for some $M$ are unique or not. It is a numerical evidence that it holds if $M>0$, but unclear if $M<0$.
  \subsection{The case $M<0$, $N=2$ and $p>1$}
  A first difficulty in this case comes from the fact that there exist singular solutions $u$ with a logarithmic blow-up. The main difficulty comes from the equality of $\gm^*(2)$ and $\overline\gm$. Hence $-\overline\gm$ is no longer a weak source as in the case $N>2$.
  \bth{en2} Assume $N=2$, $p>1$. There exist positive numbers $\tilde\gm_{min}$ and $\tilde\gm_{max}$
such that $-\gm^*(2)<\tilde\gm_{min}\leq \tilde\gm_{max}<\gm^*(1)$ with the following properties:\\
1- For $-\tilde\gm_{min}<M<-\gm^*(2)$ there exists a ground state.\\
2- for $M<-\tilde\gm_{max}$ there exists a ground state $u$ either such that $u(r)\sim U_{2,M}(r)$ or ondulating around $U_{2,M}(r)$ when $r\to\infty$.\\
 3- $M=-\tilde\gm_{min}$  there exists a ground state $u$ such that $u(r)\sim U_{1,M}(r)$ when $r\to\infty$.
  \es
  \Proof If $M=\overline M=-\gm^*(2)$ there exists no ground state from \rprop{cross}. By continuity this property is still valid for $\overline M^\ge=\overline M-\ge$ for $\ge>0$ small enough.
As in the proof of \rth{psubM<0} with $N\geq3$ we still denote by $A$ the set of $M\in (-\gm^*(1),-\gm^*(2))$ such that there is no ground state. We define in a similar way the set $B_1$ and $B_2$. The previous situation
is still valid with the only difference that $\overline M_\ge$ does not satisfies A-(i) but A-(ii): indeed from \cite[Th. 8.2, Lemma 8.7]{Kuz}, see Appendix,  for $\ge<\ge_0$ small enough there is no cycle around $P_{2,M}$ which is a sink by \rprop{2points}. Thus ${\bf T}^{1,3}_{unst}$ converges to $P_2,M$ when $t\to\infty$ and ${\bf T}^{1,2}_{st}$ converges to $0$ when $t\to-\infty$, and since there is no ground state it satisfies $x(t)\sim ce^{\frac{2t}{p-1}}\abs t$ when ${t\to-\infty}$ for some $c>0$.
Hence the function $g:M\mapsto g(M)=x^M_{reg}-x^M_{1,2}$ defined as in the proof of \rth{psubM<0} shares the same properties and the conclusion follows. $\phantom{----}$ $\phantom{----}$ \qeda
   \medskip

   \nind\Remark We conjecture that there is no cycle when $N=2$. If it is true, then for any $-\tilde\gm<M<-\gm^*$, ${\bf T}^{1,2}_{st}$ converges to $(0,0)$ as $t\to-\infty$. Equivalently there exists a positive solution $u$ of $(\ref{I-1})$ with a logarithmic blow-up at $r=0$ and such that $u(r)\sim U_{1,M}(r)$ as $r\to\infty$. Hence there exist also a positive solution $u$ of $(\ref{I-1})$ such that
   $u(r)\sim U_{1,M}(r)$ as $r\to 0$ and $u(r)\sim U_{2,M}(r)$ as $r\to \infty$.
     \subsection{The case $M<0$, $N=1$ and $p>1$}
     In the case $N=1$, the equation is invariant under the translation group $\CT_\ga[u](.)=u(.+\ga)$ for $\ga=0$ and any ground state is symmetric with respect to its vertex.

   \bth{un} Let $N=1$. Then there exists a ground state $u$ if and only if $M\leq-\gm^*(1)$. It satisfies $u(r)\sim U_{2,M}(r)$ as $r\to\infty$.
   Furthermore, if $M<-\gm^*(1)$ there exists a positive singular solution $u$ which satisfies $u(r)\sim U_{1,M}(r)$ as $r\to 0$ and $u(r)\sim U_{2,M}(r)$ as $r\to\infty$.
   \es
   \Proof The existence when $M<-\gm^*(1)$  is proved in \rlemma{m1} but the proof therein is not valid when $M=-\gm^*(1)$ in which case a second beautiful construction due to Chipot and Weissler \cite{ChWe} applies: if  $M\leq-\gm^*(1)$ there exist singular
   solutions $U_{1,M}$ and $U_{2,M}$. If ${\bf T}_{reg}$ is not a ground state, the corresponding solution $u$ vanishes at
   $r=r_0>0$. Hence there  exists a translation of $u$, say $r\mapsto u(r-c)$ which is tangent to $U_{1,M}$ which is impossible. If  $M<-\gm^*(1)$ estimate $({\ref{M19}})$ implies $({\ref{M21}})$ which in turn implies that $r^{\frac{2}{p-1}}u(r)$ cannot converge to $X_{1,M}$. Notice that there exists no cycle in the phase plane $(x,y)$ otherwise the corresponding solution $u$ would be singular and ondulating hence a translation of it say $x\mapsto u(x+c)$ which is now singular at $x=-c$ and defined for $x>-c$ could be made tangent somewhere to $U_{1,M}$ (or $U_{2,M}$) which is impossible. Therefore $r^{\frac{2}{p-1}}u(r)$ converge to $X_{2,M}$ as $r\to\infty$.  \smallskip

   \nind In order to prove that there exists a heteroclinic connecting $P_{1,M}$ to $P_{2,M}$ and since ${\bf T}_{reg}$ converges to $P_{1,M}$, there exists a smallest $\gt$ such that $x_{reg}(\gt)=X_{1,M}$ and the vector field ${\bf H}$ is directed to the right on the segment $J=\{(x,y): x=X_{1,M}, y_{reg}(\gt)\leq y\leq Y_{1,M}\}$,  we have three possibilities:\\
   (i) either $x_{reg}(t)\to X_{1,M}$ monotonically. In that case the region bordered by the segment $J$,  the portion of $\CL$ between $P_{1,M}$ and $P_{2,M}$ and the portion of trajectory ${\bf T}_{reg}$ for $t>\gt$ is positively invariant. Since ${\bf T}^{1,3}_{unst}$ belongs to this region, it converges to $P_{2,M}$ when $t\to\infty$.\\
   (ii) either  ${\bf T}_{reg}$ has a first intersection with $\CL$ at a point $(x(t_1),y(t_1))$ with $x(t_1)>X_{2,M}$. Then it enters successively the region {\bf D} where $x_t<0$ and $y_t>0$ and the region  the region {\bf A} where $x_t<0$ and $y_t<0$ and finally intersects $\CL$ between $P_{1,M}$ and $P_{2,M}$ at some point $(x(t_2),y(t_2))$, or converges to $P_{2,M}$ monotonically, in which case we set $t_2=\infty$. In that case the region bordered by the segment $J$,  the portion of $\CL$ between $P_{1,M}$ and $(x(t_2),y(t_2))$ and the portion of trajectory ${\bf T}_{reg}$ for $\gt<t<t_2$ is positively invariant. Since ${\bf T}^{1,3}_{unst}$ belongs to this region we conclude as in case (i).
\qeda
\subsection{Appendix: Non-existence of cycle in the case $N=2$}
The difficulty comes from the fact that when $M=\overline M=-\overline\gm=-\gm^*(2)$, the linearized system at $P_{2,\overline M}$ has zero as a double eigenvalue. Following Kuznetsov's notations \cite[Lemma 8.7]{Kuz} we consider the system associated to $(\ref{I-3-2})$, with two extra parameters $\overrightarrow{\ga}=(\ga_1,\ga_2)$ called the bifurcation parameters,
\bel{A1} \BA {lll}
x_t=v\\
v_t=\ga_2+\myfrac{4}{p-1}v-\myfrac{4}{(p-1)^2}x-x^p+(\overline\gm+\ga_1)\left(\myfrac{2}{p-1}x-v\right)^{\frac{2p}{p+1}}:=g(x,v,\ga_1,\ga_2).
\EA
\ee
We recall that $\overline\gm=(p+1)p^{-\frac{p}{p+1}}$ and set $\overrightarrow x=(x,v)$,  $g(x,v,\ga_1,\ga_2)=g(\overrightarrow x,\overrightarrow{\ga})$ and
$$   y_0=Y_{\overline M,2}=\frac{1}{p^{\frac{1}{p-1}}}\left(\frac{2}{p-1}\right)^{\frac{p+1}{p-1}}=\,,\; x_0=X_{\overline M,2}=\frac{p-1}{2}y_0
$$
 linearize $(\ref{A1})$ at $\tilde P_0=(x_0,0)$, with $\overrightarrow{\ga}$ fixed, we obtain the new system
 \bel{A2} \BA {lll}
\overline x_t=v\\[2mm]
v_t=\ga_2+\ga_1y_0^{\frac{2p}{p+1}}+\myfrac{2\ga_1p}{p^2-1}y_0^{\frac{p-1}{p+1}}\overline x
-\frac{2\ga_1p}{p+1}y_0^{\frac{p-1}{p+1}}v-\ga_1y_0^{\frac{2p}{p+1}}+\ga_2\\[4mm]
\phantom{v_t}
+\myfrac{1}{2}\left(g_{xx}(P_0)\overline x^2+2g_{xv}(P_0)\overline xv+g_{vv}(P_0)v^2\right)+R(\overrightarrow x,\overrightarrow{\ga}).
\EA
\ee
In order to agree with Kuznetsov's notations, we write $(\ref{A2})$ under the form
 \bel{A3} \BA {lll}
y_{1\,t}=y_2\\[2mm]
y_{2\,t}=g_{00}+g_{10}y_1+g_{01}y_2
\myfrac{1}{2}\left(g_{20}y_1^2+2g_{11}y_1y_2+g_{02}y_2^2\right)+R(\overrightarrow y,\overrightarrow{\ga})
\EA
\ee
where $R(\overrightarrow y,\overrightarrow{\ga})=O(|\overrightarrow y|^3)$ and
$$\BA{lll}g_{00}=\ga_2+\ga_1y_0^{\frac{2p}{p+1}}\,,\;g_{1,0}=\myfrac{4p\ga_1}{p^2-1}y_0^{\frac{p-1}{p+1}}\,,\;g_{01}=-\myfrac{2p\ga_1}{p+1}y_0^{\frac{p-1}{p+1}}\,,\;g_{11}=-\myfrac{4p}{(p+1)^2}\left(\myfrac {p+1}{p^{\frac{p}{p+1}}}+\ga_1\right)y_0^{-\frac{2}{p+1}}\\[4mm]
g_{20}=-\myfrac{4p}{(p+1)^2}\left(\myfrac {p+1}{p^{\frac{p}{p+1}}}-\myfrac{2\ga_1}{p-1}\right)y_0^{-\frac{2}{p+1}}\,,\;g_{02}=-\myfrac{2p(p-1)}{(p+1)^2}\left(\myfrac {p+1}{p^{\frac{p}{p+1}}}+\ga_1\right)y_0^{-\frac{2}{p+1}}.
\EA$$
Note that, if $\ga_1=0$, the three coefficients $g_{11}$, $g_{20}$ and $g_{02}$ are negative. \\
Following Kuznetsov proof, we perform several changes of variables:\\
1- Setting $y_1=v_1+\gd$, $v_2=y_2$ where $\gd=\gd(\overrightarrow\ga)$, we can get rid of the coefficient of $v_2$ in the second equation and obtain
 \bel{A4} \BA {lll}
v_{1\,t}=v_2\\[2mm]
v_{2\,t}=h_{00}+h_{10}y_1+
\myfrac{1}{2}\left(h_{20}v_1^2+2h_{11}v_1v_2+h_{02}v_2^2\right)+Q(\overrightarrow v,\overrightarrow{\ga})
\EA
\ee
where $Q(\overrightarrow v,\overrightarrow{\ga})=O(|\overrightarrow v|^3)$ and
$$\gd=-\myfrac{g_{01}}{g_{11}(0)}(1+o(1))\quad\text{as }\;\overrightarrow{\ga}\to 0,
$$
with $g_{01}$ (resp. $g_{00}$ and $g_{10}$) stands for $g_{01}(\overrightarrow{\ga})$ (resp. $g_{00}(\overrightarrow{\ga})$ and $g_{10}(\overrightarrow{\ga})$) and
$$h_{00}=g_{00}-\myfrac{g_{01}g_{10}}{g_{11}(0)}+o(1)\,,\; h_{10}=g_{10}-\myfrac{g_{01}g_{20}}{g_{11}(0)}+o(1)
\,\text{ and }\;h_{20}=g_{20}\,,\;h_{11}=g_{11}\,,\;h_{02}=g_{02}.
$$
2- Time scaling, $\frac{dt}{d\gt}=1+\gth v_1(t)$ where $\gth=\gth(\overrightarrow{\ga})$, and $\xi_1=v_1$, $\xi_2=(1+\gth v_1))v_2$ in order to get rid of the coefficients of $v_2^2$ in the equations for $\xi_2$, then
 \bel{A5} \BA {lll}
\xi_{1\,\gt}=\xi_2\\[2mm]
\xi_{2\,\gt}=f_{00}+f_{10}\xi_1
\myfrac{1}{2}hf_{20}\xi_1^2+f_{11}\xi_1\xi_2+h_{02}v_2^2+P(\overrightarrow \xi,\overrightarrow{\ga})
\EA
\ee
where $P(\overrightarrow \xi,\overrightarrow{\ga})=O(|\overrightarrow \xi|^3)$,
$$f_{00}=h_{00}(1+o(1))\,,\; f_{10}=\left(h_{10}-\myfrac{h_{00}h_{02}}{2}\right)(1+o(1))\,,\; f_{20}=\left(h_{20}-h_{10}h_{02}\right)(1+o(1)),
$$
and $f_{11}=h_{11}(1+o(1))$. We rewrite the equation of $\xi_{2\,t}$ under the form
 \bel{A6} \BA {lll}
\xi_{2\,t}=\gm_1(\overrightarrow{\ga})+\gm_2(\overrightarrow{\ga})\xi_1+A(\overrightarrow{\ga})\xi_1^2+B(\overrightarrow{\ga})\xi_1\xi_2
+O(|\overrightarrow \xi|^3),
\EA
\ee
where, for $|\overrightarrow \ga|$ small enough,
$$\gm_1=h_{00}=\left(g_{00}-\myfrac{g_{10}g_{01}}{g_{11}(0)}\right)(1+o(1)))=\ga_1y_0^{\frac{2p}{p+1}}(1+o(1)),
$$
$$\gm_2=\left(h_{10}-\myfrac{h_{02}h_{00}}{2}\right)(1+o(1)))=\myfrac{2(p^2+2p-1)}{p^2-1}y_0^{\frac{p-1}{p+1}}(1+o(1)),
$$
$$A(\overrightarrow{\ga})=-\myfrac{2p}{p+1}\left(\myfrac1{p}\right)^{\frac{p}{p+1}}y_0^{-\frac{2}{p+1}}(1+o(1)),
$$
and
$$B(\overrightarrow{\ga})=-\myfrac{4}{p+1}\left(\myfrac1{p}\right)^{\frac{p}{p+1}}y_0^{-\frac{2}{p+1}}(1+o(1)).
$$
At end we change again time and put
$$t=\abs{\myfrac{B}{A}}\gt\,,\; \eta_1=\myfrac{B^2}{A}\xi_1\,\text{ and }\;\eta_2=sgn\left(\myfrac{B}{A}\right)\myfrac{B^3}{A^2}\xi_2,
$$
in order to see that $\overrightarrow\eta=(\eta_1,\eta_2)$ verifies
 \bel{A7} \BA {lll}
\eta_{1\,t}=\eta_2\\[2mm]
\eta_{2\,t}=\gb_1(\overrightarrow{\ga})+\gb_2(\overrightarrow{\ga})\eta_1+\eta_1^2+sgn\left(\myfrac{B(0)}{A(0)}\right)\eta_1\eta_2+O(|\overrightarrow \eta|^3).
\EA
\ee
In our situation
$$sgn\left(\myfrac{B(0)}{A(0)}\right)=1.
$$
After further computations and simplifications we obtain, with $\ga_2=0$ and eliminating  the terms which contain $\ga_1^2$,
 \bel{A8} \BA {lll}
\gb_1=-\myfrac{64}{p^2-1}\ga_1y_0^{\frac{p-1}{p+1}}(1+o(1))\,\text{ and }\;\gb_2=\myfrac{8p^2+2p-1}{p^2-1}\ga_1y_0^{\frac{p-1}{p+1}}(1+o(1)).
\EA
\ee
Therefore, the discriminant $\gb^2_2-4\gb_1$ of the polynomial $\CP(\eta_1)=\gb_1+\gb_2\eta_1+\eta_1^2$
is positive for $\ga_1>0$.

By \cite[Lemma 8.7]{Kuz} there is no cycle in the region of the plane $(\gb_1,\gb_2)$ located in the second quadrant, which is our case since $\ga_1>0$. Furthermore, the equilibrium $(\gb_1,\gb_2)=(0,0)$,  which has a double zero eigenvalue has one stable trajectory converging when $t\to\infty$ and one unstable trajectory converging when $t\to-\infty$. Hence it is a saddle point.

\end{document}